\documentclass[12pt]{article}

\usepackage{epic}
\usepackage{eepic}
\usepackage{epsf} 
\usepackage{ amssymb, amscd}
\usepackage{amsmath}
\usepackage{color}
\numberwithin{equation}{section}
\usepackage{epsfig}
\usepackage{graphicx}
\usepackage{tikz}
\usetikzlibrary{matrix}

\def\cA            {{\mathcal{A}}}
\def\cB            {{\mathcal{B}}}
\def\cC            {{\mathcal{C}}}
\def\cD            {{\mathcal{D}}}
\def\cE            {{\mathcal{E}}}

\def\cG            {{\mathcal{G}}}
\def\cH            {{\mathcal{H}}}
\def\cJ            {{\mathcal{J}}}
\def\cK            {{\mathcal{K}}}
\def\cL            {{\mathcal{L}}}
\def\cM            {{\mathcal{M}}}

\def\cO            {{\mathcal{O}}}
\def\cP            {{\mathcal{P}}}
\def\cQ            {{\mathcal{Q}}}

\def\cS            {{\mathcal{S}}}
\def\cT            {{\mathcal{T}}}
\def\cV            {{\mathcal{V}}}

\def\cY            {{\mathcal{Y}}}
\def\cZ            {{\mathcal{Z}}}

\def\bbC           {\mathbb{C}}

\def\bbQ           {\mathbb{Q}}
\def\bbR           {\mathbb{R}}
\def\bbT           {\mathbb{T}}
\def\bbZ           {\mathbb{Z}}

\def\ran           {\rangle}
\def\lan           {\langle}

\def\la            {\lambda}

\def\ran           {\rangle}

\def\i{{\rm i}}
\def\sdprod{{\times\!\vrule height5pt depth0pt width0.4pt\,}}
\def\boxit#1{\vbox{\hrule\hbox{\vrule{#1}\vrule}\hrule}}
\def\stimes{\,\,{\boxit{$\times$}}\,\,}
\def\eps{\epsilon}
\def\bfe           {{\bf1}}

\title{Tambara-Yamagami, loop groups, bundles and KK-theory}
\author{
{\sc David E.\ Evans}\\
 {\footnotesize School of Mathematics, Cardiff University,}\\
 {\footnotesize Senghennydd Road, Cardiff CF24 4AG, Wales, U.K.}\\
 {\footnotesize e-mail: {\tt EvansDE@cf.ac.uk}}\\ \\
 {\sc Terry  Gannon }\\
 {\footnotesize Department of Mathematics, University of Alberta,}\\
{\footnotesize Edmonton, Alberta, Canada T6G 2G1}\\
{\footnotesize e-mail: {\tt tgannon@math.ualberta.ca}} }

\begin{document}
\maketitle

\begin{abstract} 
This paper is part of a sequence interpreting quantities of conformal field theories $K$-theoretically. Here we give geometric constructions of the associated module categories (modular invariants, nimreps, etc). In particular,  we give a $KK$-theory interpretation of all modular invariants for the
loop groups of tori, as well as most known modular invariants of loop groups. In addition, we find unexpectedly that   the  Tambara-Yamagami fusion category has an elegant description as bundles over a groupoid, and use that to interpret its  module categories as  $KK$-elements. We establish reconstruction for  the doubles of all Tambara-Yamagami categories, generalizing work of Bischoff to even-order groups. We conclude by relating the 
modular group representations coming from finite groups and loop groups  to the Chern character and to the Fourier-Mukai transform respectively.

\end{abstract}

{\footnotesize
\tableofcontents
}

\section{Introduction}
The applicability of $K$-theory to conformal field theory {(CFT)}, e.g. to D-brane charges  
\cite{MiMo} or to the  fusion rings \cite{FHT}, is clear. In a series of papers \cite{Ev,EG1,EG3,EG7,EG8}, the authors have been extending the range of these applications.
This paper is part of that sequence, in that much of its inspiration derives from that story.

This paper brings together three themes. One is the elegance of using $K$-theory and $KK$-theory to capture fusion categories and module categories.   The fusion ring, describing the product structure of the primary fields or sectors, is the most elementary structure of chiral CFT (which also contains a modular tensor category). Analogously, the most
elementary object in the full CFT (which is a module category over that modular tensor category) is the modular invariant partition function (\textit{modular invariant}
for short), which describes how the two chiral halves glue together. It is an integral matrix indexed by the
primaries (the preferred basis in the fusion ring, namely the irreducible representations of the chiral algebra or vertex operator algebra or conformal net). 
In particular, when the fusion ring  has a natural $K$-group realisation $K(X)$, then the modular invariant 
 is
a linear map between $K$-groups, i.e. is an element of $KK(X,X)$. We expect in these cases that the modular invariants can be expressed in a natural way as very special $KK$-elements. The most important manifestation of this idea will be to realise them as spectral triples (see \cite{FrTe} for a preliminary step in this direction). As a first step however, in this
paper we realise these $KK$-elements with correspondences.
The result, as we'll see, is the most elegant expression for those modular invariants that we have seen.

We see this strategy implicitly for the case of finite groups in \cite{EG7} (developed further and more explicitly in \cite{EG8}). The fusion ring of  the Drinfeld double $\cD(G)$ of a finite group $G$ is naturally identified with the equivariant $K$-group $K^0_G(G)$, where $G$ acts on itself by conjugation.   Finite group doubles describe the representation theory of finite group orbifolds of a chiral algebra with a trivial representation theory. By abstract considerations \cite{Ost}, the module categories of $\cD(G)$  
are parametrised up to equivalence by a subgroup $H\le G\times G$ and a 2-cocycle class $[\psi]\in H^2_H(\mathrm{pt}
;\bbT)$ (here and elsewhere, $\bbT$ denotes the unit circle in $\bbC$). However, it is difficult to identify the corresponding modular invariant matrix (or indeed most other aspects of the full CFT). But $KK$-theory provides an elegant answer. We explain in \cite{EG8} that this matrix is given by
the correspondence (in the sense of \cite{CoSk},\cite{EM})
\begin{equation}\begin{tikzpicture}\label{modinvgen1}
  \matrix (m) [matrix of math nodes,row sep=2em,column sep=2em,minimum width=2em]
  {
    & (H/\!/H^{adj},\beta^\psi) &  \\ G/\!/ G^{adj}&&G/\!/G^{adj}\\ };
  \path[-stealth]
    (m-1-2) edge node [left] {$p_L\ $} (m-2-1)
            edge node [right] {$\ p_R$} (m-2-3);
\end{tikzpicture}\end{equation}
where $G/\!/G^{adj}$ refers to the groupoid with $G$ acting adjointly on itself, and
where $p_{L,R}$ are the obvious coordinate maps and $\beta^\psi$ is a certain line bundle
depending on $\psi$. From this the matrix entries and other
information of the modular invariant (e.g. its type 1 parents) can be immediately read off. Indeed, \eqref{modinvgen1} describes the geometric realisation of the module category parametrised by $(H,\psi)$, as bundles over groupoids.

The loop groups are an important source of theories, and they have a similar $K$-theoretic description. In seminal work, Freed-Hopkins-Teleman \cite{FHT,FHTi,FHTii,FHTlg} identified
the fusion ring of the loop group $LG$ at level $k$, for any compact simple connected simply-connected Lie group $G$, with the twisted equivariant
$K$-group $^\tau K^d_G(G)$ for some $\tau\in Z^3_G(G;\bbZ)$ depending on $k$,
where $d$ is the dimension of $G$ and $G$ acts on itself by conjugation. In this paper we identify the $KK$-description of almost all module categories of the loop groups.

Another class of simple examples are the loop groups of tori  $\bbR^d/\bbZ^d$.  Here, the possible chiral data (e.g. modular tensor categories)
are parametrised by $d$-dimensional even lattices $L$. For example, the fusion ring will be
the group ring $\bbZ[G]$ of the finite abelian group $G=L^*/L$, where $L^*$ is the dual of $L$, and the inner product $\langle,\rangle$ on $L$ identifies $G$ with its irreps $\widehat{G}$. We give many different parametrisations of their modular invariants, each with their associated $KK$-element. Here is one such parametrisation: pairs $(H,[\psi])$, where $H\le G$ and $[\psi]\in H^2_H(\mathrm{pt};\bbT)$. Compare this classification with that of the modular categories for finite group doubles given a couple paragraphs earlier. Using $H$ and $\psi$, there is a second homomorphism $G\to\widehat{G}$, which we call $\epsilon$, and the modular invariant is the comparison of the two:
\begin{equation}\begin{tikzpicture}\label{modinvtorigen}
  \matrix (m) [matrix of math nodes,row sep=2em,column sep=2em,minimum width=2em]
  {
    K(G)&&K(G)\\ & K(\widehat{G})&   \\  };
  \path[-stealth]
    (m-1-1) edge node [left] {$\langle,\rangle\ $} (m-2-2) 
       (m-1-3)     edge node [right] {$\ \epsilon$} (m-2-2);
\end{tikzpicture}\end{equation}

Any modular tensor category comes with a unitary matrix representation of SL$_2(\bbZ)$ (hence their name). Building on the discussions mentioned above,  we conclude the paper with a $K$-theoretic reinterpretation of modular data for these finite group, loop group, and toroidal modular tensor categories.

The examples mentioned so far are quite classical, being completely group or Lie theoretic. In this paper we unexpectedly discover a nonclassical family possessing a $K$-theoretic description, and this leads to our second theme: 
the Tambara-Yamagami fusion categories \cite{TY}. They are the simplest quadratic categories, i.e. quadratic extensions of  fusion categories with fusion ring $\bbZ[G]$. The quadratic categories have attracted attention recently as many of them (most famously the Haagerup \cite{iz3,EG2})  seem to be exotic in the sense that they have no known direct constructions starting from classical structures like finite groups or Lie algebras.
The Tambara-Yamagami fusion ring is spanned by $[\alpha_g]$ for $g\in G$ ($G$ a finite additive abelian group), and $[\rho]$, and obeying the relations \begin{equation}\label{TYfusions}[\alpha_g][\alpha_h]=[\alpha_{g+h}]\,,\ [\alpha_g][\rho]=[\rho]=[\rho][\alpha_g]\,,\  [\rho]^2=\sum_{g\in G}
[\alpha_g]\end{equation}
 In any case, we realise the full Tambara-Yamagami categories (including associators) naturally as categories over groupoids.  We also realise their module categories as $KK$-elements. It would be very interesting to apply similar techniques to explore which other exotic fusion categories have natural $K$-theoretic descriptions.

In order to discuss modular invariants, we need to bring in modular tensor categories. Given a fusion category, one can always obtain a modular tensor category through the double construction. Moreover, we show   here that the Tambara-Yamagami categories have a $\bbZ_2$-crossed braiding, and  when the order $|G|$ is odd,  their $\bbZ_2$-equivariantisations are modular.  These modular tensor categories are closely related to, but rather simpler than, the doubles. We show the fusion rings of both the doubles and $\bbZ_2$-equivariantisations are also captured by $K$-theory. Being metaplectic categories, they have enormous numbers of modular invariants, but few of these will be sufferable (i.e. correspond to module categories or full CFTs). We show how to realise some of these sufferable modular invariants by $KK$-elements. 

 Theme 3 concerns a
 foundational question in the theory, namely reconstruction: whether or not it is possible to realize all (unitary) modular tensor categories as the category of modules of an RCFT chiral algebra, i.e.\ of a rational vertex operator algebra or local conformal net of factors. This is trivial for the torus and loop group examples, and in \cite{EG7} we show this is true for any twisted double of a finite group.  In this paper we establish reconstruction for both the doubles and  $\bbZ_2$-equivariantisations of Tambara-Yamagami categories. This generalises the work of Marcel Bischoff  \cite{Bis}, who established reconstruction for Tambara-Yamagami when $|G|$ is odd. Bischoff's work for odd $|G|$ was independent and simultaneous to ours --- we announced our reconstruction for even and odd $|G|$ in a  talk by the second author at the Isaac Newton Institute, available online \cite{GanINI}.

In the process of doing this we prove that any \textit{pointed} modular tensor category (i.e. one whose simple objects are all invertible) can be reconstructed as a lattice theory. This was an explicit assumption in Bischoff \cite{Bis}.
More generally, a weakly-integral fusion category is one where all dimensions-squared are integers. It is tempting to guess that any weakly-integral modular tensor category can be reconstructed as a group orbifold of a lattice theory. This generalises what we now know to be true for doubles and $\bbZ_2$-equivariantisations of the Tambara-Yamagami categories.

One reason for being especially interested in Tambara-Yamagami categories is because of their possible relevance to reconstruction for the doubles of the Haagerup-Izumi series of subfactors. More precisely, the modular data of the double of the Haagerup can be most easily recovered as the grafting (see \cite{EG7}) of the $\bbZ_2$-equivariantisations of a $\bbZ_3\times\bbZ_3$ and  $\bbZ_{13}$ Tambara-Yamagami category. We describe this in section 5.6.

\section{Background}

\subsection{$KK$-theoretic background {via correspondences}}

The $KK$-groups can be defined as follows. For $C^*$-algebras $A,B,E$, the extensions
\begin{equation}\label{**}
0 \rightarrow \cK \otimes B \rightarrow E \rightarrow A
\rightarrow 0
\end{equation}
together with suspensions, where $\cK$ is the algebra of compact operators of a separable infinite-dimensional Hilbert space, yield the Kasparov groups
$KK_\star(A,B)$ (p. 118 of \cite{EKaw}).  Now,
the Universal Coefficient Theorem (see e.g. section 23.1 in \cite{Black}) says there is a short exact sequence
\begin{equation}\label{UCT}
0\rightarrow\mathrm{Ext}^1_\bbZ(K_\star(A),K_{\star-1}(B))\rightarrow KK_\star(A,B)
\rightarrow\mathrm{Hom}(K_\star(A),K_\star(B))\rightarrow 0\,.\end{equation}
Therefore if either $K_\star(A)$
or $K_\star(B)$ is  a free finitely generated $\bbZ$-module ---  which will always the case for us
--- then 
$KK_\star(A,B)\cong \mathrm{Hom}(K_\star(A),K_\star(B))$. 

For spaces (the language we find more convenient), this becomes $KK_\star(X,\mathrm{pt})=K_\star(X)$, 
$K^\star(\mathrm{pt},Y)=K^\star(Y)$, and $KK_\star(X,Y)=\mathrm{Hom}(K^\star(X),K^\star(Y))$.
 In particular, taking $X = Y$
to be the object giving the fusion ring, a modular
invariant is  an element of End$(K(X))$ and so
gives rise to an element of $KK(X,X)$. 

Hence a modular
invariant gives rise to very special $KK$-elements, as
do other structures in CFT such as sigma-restriction and alpha-induction. The intersection or Kasparov product of $KK$-elements then corresponds to
the usual matrix multiplication of modular invariants. And so on.

Elements of $KK$-groups can be described through spectral triples, Fredholm modules and Dirac operators \cite{Connes book}. However
the approach we take in this paper is to follow the correspondence description due to Connes and Skandalis \cite{CoSk}
for topological manifolds, in particular manifolds and foliations, and most recently developed by Emerson and Meyer
\cite{EM} for groupoid equivariant theory on manifolds.

Let $X,Y$ be smooth manifolds. A \textit{correspondence}  (defined in \cite{CoSk}, refined by \cite{EM}; see also \cite{BMRS}
for a gentle introduction) is given by the diagram
$$\begin{tikzpicture}
  \matrix (m) [matrix of math nodes,row sep=2em,column sep=2em,minimum width=2em]
  {
    & (Z,E) &  \\ X& &Y\\ };
  \path[-stealth]
    (m-1-2) edge node [left] {$f\ $} (m-2-1)
            edge node [right] {$\ b$} (m-2-3);
\end{tikzpicture}$$
where $Z$ is a smooth manifold, 
$E$ is a complex vector bundle over $Z$, the forward map $f : Z\rightarrow X$ is smooth and
proper, and the backward map $b : Z \rightarrow Y$ is $K$-oriented. Any correspondence naturally defines a class 
\begin{equation}\label{kkcorr}b_! (f^* (-)\otimes E)\end{equation}
in the bivariant $K$-theory group $KK(X, Y ) := {KK}(C_0 (X), C_0 (Y ))$, where $f^*$ is the pullback and $b_!$ is the pushforward.  The
collection of all correspondences for $X, Y$ forms an additive category under disjoint union:
\begin{equation}
{(Z_1,E_1)+(Z_2,E_2)=(Z_1\sqcup Z_2,E_1\sqcup E_2)}\,.\label{bundlesum}
\end{equation}
{Quotienting by appropriate notions of cobordism, direct sum and vector bundle
modification recovers the $KK$-theory group $KK(X, Y )$.} {In this quotient, a special case of \eqref{bundlesum} is $[Z,E_1]+[Z,E_2]=[Z,E_1\oplus E_2]$.}

{In} the correspondence picture,  the intersection product on $KK$-theory is  the bilinear associative map
$$\otimes_M:{KK}(X,M)\times{KK}(M,Y)\rightarrow {KK}(X,Y)$$
defined by the pullback,  sending two correspondences 
$$\begin{tikzpicture}
  \matrix (m) [matrix of math nodes,row sep=2em,column sep=2em,minimum width=2em]
  {
    & (Z_1,E_1) & &(Z_2,E_2)& \\ X& & M&&Y\\ };
  \path[-stealth]
    (m-1-2) edge node [left] {$f\ $} (m-2-1) edge node [right] {$\ b_M$} (m-2-3)
     (m-1-4) edge node [left] {$f_M\ $} (m-2-3) edge node [right] {$\ b$} (m-2-5);
\end{tikzpicture}$$
to the single correspondence {$(Z,E)=(Z_1,E_1)\otimes_M(Z_2,E_2)$} with $Z=Z_1\times_M Z_2=\{(z_1,z_2)\in Z_1\times Z_2\,:\,
b_M(z_1)=f_M(z_2)\}$,
the fibred product, and $E=E_1\otimes_M E_2$ is the restriction of the  $Z_1\times Z_2$-bundle $E_1\otimes E_2$ to
 $Z_1\times_M Z_2$. {The maps $f'$ and $b'$ from $(Z,E)$ are then just the obvious restrictions.}
 This definition requires a transversality condition on the
two maps $f_M$ and $b_M$ in order to ensure that the fibred product $Z$ is a smooth manifold.

Equivariant $KK$-theory on topological spaces is too narrow for us, as we will have different groups
acting on different algebras or spaces. In general we need $KK$-theory on groupoids. We describe 
next how to modify the correspondence picture on topological spaces and manifolds to correspondences between groupoids
and hence elements of their $KK$-groups.

A  groupoid is a category whose morphisms $f$ have both left and right inverses.
 Groupoid theory is summarised for instance in Appendix A of \cite{FHTi}. When a group $G$ 
 acts on a set $X$, we will write $X/\!/G$ for the corresponding groupoid.
A morphism (or map) between groupoids is a functor between the corresponding categories. 
{In the loop group setting, topology must be considered, and Lie groupoids used;
in this case all maps are required to be smooth.}

We will be twisting our groupoids by 2- and 3-cocycles; in this framework, {the twist} does not
affect the underlying groupoid but does affect the {class of} bundles considered.
 An (untwisted) \textit{bundle} over a groupoid   is simply a functor from the groupoid to the category of
finite-dimensional vector spaces; the twist controls the projectivity of the groupoid action on the fibres. 
For example, untwisted bundles over pt$/\!/G$ ($G$ a finite group) correspond precisely to
representations of $G$; for $\psi\in Z^2_G(\mathrm{pt};\bbT)$, $\psi$-twisted bundles on pt$/\!/G$ 
correspond precisely to projective representations of $G$ with cocycle $\psi$. By a slight abuse of notation,
we will speak of e.g. bundles over  pt$/\!/_\psi G$ rather than $\psi$-twisted bundles over pt$/\!/G$. 
Given a groupoid $\cG$ and such a twist $\tau$, let $\cC_\tau(\cG)$ denote the category
whose objects are $\tau$-twisted bundles over $\cG$, and whose morphisms are natural transformations 
between those bundles. The corresponding $K$-group will be denoted ${}^\tau K^0(\cG)$ (or ${}^\tau
K^0_G(X)$ in the case of action groupoids $X/\!/_\tau G$). See section 2.2  for detailed examples.

The \textit{weak pullpack for groupoids}  is defined in section 2.5 of \cite{EM} and section 2 of  \cite{BHW}: given groupoid maps 
$$\begin{tikzpicture}
  \matrix (m) [matrix of math nodes,row sep=2em,column sep=2em,minimum width=1em]
  {
   T & & S \\ & X&\\ };
  \path[-stealth]
    (m-1-1) edge node [left] {$q\ $} (m-2-2)
     (m-1-3)       edge node [right] {$\ p$} (m-2-2);
\end{tikzpicture}$$
the pullback $P$ will have objects consisting of all triples $(s,t,\alpha)$ where $s$ and $t$ are
objects of groupoids $S$ and $T$ respectively, and $\alpha:p(s)\rightarrow q(t)$ is a morphism 
in $X$. A morphism in $P$ from $(s,t,\alpha)$ to $(s',t',\alpha')$ consists of a morphism
$f:s\rightarrow s'$ in $S$ and a morphism $g:t\rightarrow t'$ in $T$ such that the diagram
$$\begin{tikzpicture}
  \matrix (m) [matrix of math nodes,row sep=2em,column sep=2em,minimum width=1em]
  {
   p(s) & q(t)  \\ p(s') &q(t')\\ };
  \path[-stealth]
    (m-1-1) edge node [left] {$p(f)\ $} (m-2-1)  
    (m-1-1) edge node [above] {$\alpha $} (m-1-2)  
     (m-2-1)       edge node [below] {$\ \alpha'$} (m-2-2)
     (m-1-2)       edge node [right] {$\ q(g)$} (m-2-2) ;
\end{tikzpicture}$$
commutes. $P$ inherits its (Lie) groupoid structure from $S$ and $T$.
{A \textit{weak} pullback $P$ satisfies the weak universality property, in the following sense.}
Suppose $R$ is a groupoid, and $F:R\rightarrow S$ and $G:R\rightarrow T$ are groupoid
morphisms. We say that the diagram
\begin{equation}\begin{tikzpicture}
  \matrix (m) [matrix of math nodes,row sep=2em,column sep=2em,minimum width=1em]
  {
   &R&\\T && S  \\ &X&\\ };
  \path[-stealth]
    (m-1-2) edge node [left] {$G\ $} (m-2-1)  
    (m-1-2) edge node [right] {$\ F $} (m-2-3)  
     (m-2-1)       edge node [left] {$q\ $} (m-3-2)
     (m-2-3)       edge node [right] {$\ p$} (m-3-2) ;
\end{tikzpicture}\label{commut}\end{equation}
\textit{commutes} if for each object $r$ of $R$, there is a morphism $\alpha\in\mathrm{Hom}_X(
p(F(r)),q(G(r)))$ of $X$, and for each morphism $h$ of $R$ we have $\alpha'p(F(h))=q(G(h))\alpha$
(where $\alpha,\alpha'$ are associated to the head resp. tail of $h$).
Weak universality means that the map $R\rightarrow P$ exists but is not necessarily unique. 
It is elementary to verify {that $P$ constructed above} is a weak pullback in this sense.

The pullback $E$ of bundles $E_S$ and $E_T$ over groupoids $S$ and $T$ respectively,
is now defined as before: it is the functor defined by $E(s,t,\alpha)=E_S(s)\otimes_\bbC E_T(t)$
and $E(f,g)=E_S(f)\otimes E_T(g)$, restricted to the pullback $P$ defined last paragraph.
{Composition
of correspondences in this more general picture is again} defined using the pullback.
Equivalence of correspondences is defined in section 2.2 of \cite{EM}.

\subsection{Warm-up examples}

We are interested in interpreting, in as natural a way as possible, the sufferable modular invariants
(i.e. those realised by module categories, defined in section 2.3)  {as $KK$-elements. To get comfortable with these $KK$-groups,  this subsection gives examples, including
matrix units, for the $KK$-rings of finite groupoids. By \textit{matrix units} $E_{i,j}$ in a matrix ring $M_{d\times d}(\bbZ)$
we mean the standard bases (matrices with all 0's except for a 1 in the $ij$-entry). These obey
$E_{i,j}E_{k,l}=\delta_{jk}E_{i,l}$.

Begin with} $KK(\mathrm{pt},\mathrm{pt})=\bbZ$. We realise the number $n\ge 0$ through the
correspondence 
$$ \begin{tikzpicture}
  \matrix (m) [matrix of math nodes,row sep=2em,column sep=2em,minimum width=2em]
  {
    & \sqcup_{j=1}^n\{\mathrm{pt}\} &  \\ \mathrm{pt}& &\mathrm{pt}\\ };
  \path[-stealth]
    (m-1-2) edge node [left] {$f\ $} (m-2-1)
            edge node [right] {$\ b$} (m-2-3);
\end{tikzpicture}$$
Composition of these correspondences amounts to multiplying the numbers.
Equivalently, we could have replaced the discrete union $\sqcup_{j=1}^n\{\mathrm{pt}\}$ here with the bundle $\bbC^n$ over a single point: {in this case,  $f^*$  associates to $k\in K(\mathrm{pt})$} the vector space $\bbC^{kn}$, while  $b_!$ takes dimension. {Composition of the correspondences for $\bbC^m$ and $\bbC^n$ is manifestly the correspondence for $\bbC^m\otimes
\bbC^n\cong\bbC^{mn}$.}
 As always, negative elements of $KK$ come from taking virtual or $K$-theoretic bundles.

Now let $G$ be any finite group. We can identify the representation ring  $K_G^\star(\mathrm{pt})=R_G$ with $KK(\mathrm{pt},\mathrm{pt}/\!/G)$.
 We can realise any given representation $V$ as the correspondence
\begin{equation} \begin{tikzpicture}
  \matrix (m) [matrix of math nodes,row sep=2em,column sep=2em,minimum width=2em]
  {
    & (\mathrm{pt}/\!/G,V) &  \\ \mathrm{pt}& &\mathrm{pt}/\!/G\\ };
  \path[-stealth]
    (m-1-2) edge node [left] {$\pi\ $} (m-2-1)
        edge node [right] {\ id} (m-2-3);
\end{tikzpicture}\label{MUgp1}\end{equation}
where $\pi$ sends $G$ to {1}. Here $\pi^*$ associates $n$ to the trivial $G$-representation $\bbC^n$ and id$_!:R_G\rightarrow
R_G$ is the identity, {so the correspondence \eqref{MUgp1} corresponds to the class $n\mapsto \bbC^n\otimes V$ in $KK(\mathrm{pt},\mathrm{pt}/\!/G)$.} Instead of the bundle $ (\mathrm{pt}/\!/G,V)$, we could instead use the disjoint union  
$\sqcup_j({\mathrm{pt}}/\!/G,\rho_j)$ where $V=\oplus_j\rho_j$ is a sum of {\textit{irreps} ($=$irreducible representations)}.

Consider now $KK_G(\mathrm{pt},\mathrm{pt})=KK(\mathrm{pt}/\!/G,\mathrm{pt}/\!/G)$, {which we identify with
End$( R_G)$.} {Let us} construct using correspondences
the matrix units for this ring. For any irreducible $G$-modules  $V,W$ define $E_{W,V}$ to
be the element of $KK_G(\mathrm{pt},\mathrm{pt})$ given by the product
$$ \begin{tikzpicture}
  \matrix (m) [matrix of math nodes,row sep=2em,column sep=2em,minimum width=2em]
  {
    & (\mathrm{pt}/\!/G,V^*) && (\mathrm{pt}/\!/G,W)& \\ \mathrm{pt}/\!/G&&\mathrm{pt}& &\mathrm{pt}/\!/G\\ };
  \path[-stealth]
    (m-1-2) edge node [left] {id\ \ } (m-2-1)
            edge node [right] {$\ \pi_R$} (m-2-3)
             (m-1-4) edge node [left] {$\pi_L\ $} (m-2-3)
            edge node [right] {\ id} (m-2-5)            ;
\end{tikzpicture}$$
Pulling this back gives the correspondence
\begin{equation} \begin{tikzpicture}
  \matrix (m) [matrix of math nodes,row sep=2em,column sep=2em,minimum width=2em]
  {
 & (\mathrm{pt}/\!/(G\times G),V^*\otimes W)& \\  \mathrm{pt}/\!/G&&\mathrm{pt}/\!/G\\ };
  \path[-stealth]
    (m-1-2) edge node [left] {{$\pi_L'\ \ $}} (m-2-1)
            edge node [right] {{$\ \pi_R'$}} (m-2-3)
             ;
\label{MUgp2}\end{tikzpicture}\end{equation}
where the morphisms are the obvious left/right projections, and $V^*$ denotes the contragredient of $V$. 
In $V^*\otimes W$ {in \eqref{MUgp2},} the left $G$ acts on $V^*$ and the right acts on $W$. Then $\pi_L^{\prime\,*}$ lifts $V'\in R_G$ to $V'\otimes 1$ (where the
left $G$ acts on $V'$ and the right acts trivially), while $\pi'_{R!}$ sends the $G\times G$-rep $V_L\otimes V_R$ to the $G$-rep {$\langle V_L,1\rangle V_R$}.
Here and elsewhere we {let $\langle V,\rho_j\rangle$}  denote the multiplicity of irrep $\rho_j$ in 
$G$-rep $V$, and extend the definition
{to other $\rho$} by $\langle V,\oplus_jn_j\rho_j\rangle =\sum_jn_j\langle V,\rho_j\rangle$.  From this we obtain that $E_{W,V}$ sends $V'\in R_G$ to $\langle V',V\rangle W$, and
so for irreps $V,W\in\mathrm{Irr}(G)$, $E_{W,V}$ is 
the $(W,V)$ matrix unit, as the notation suggests. 

Note that the inputs for correspondences are on the left, whereas inputs for matrices are written on the right. Hence in compositions, the orders of operators in the correspondence picture is the reverse of that of matrices. 

This describes the matrix units for any finite groupoid  $\cG$, by working locally in terms of orbits and stabilisers. For example, the $K$-theory for such a $\cG$ is $K^0(\cG)=\oplus_oR_{G(o)}$
where $o$ runs over all orbits in $\cG$ and $G(o)$ is the stabiliser (as an abstract group)  of $o$ in $\cG$. We can identify this with $KK(\mathrm{pt},\cG)$. Then \eqref{MUgp1} generalises to $\cG$
as follows. Fix some representative $x$ of orbit $o$, and let $G_x\cong G(o)$ be its stabiliser. 
Then to any $G_x$-module $V$, we obtain an element of $KK(\mathrm{pt},\cG)$  by the correspondence
\begin{equation} \begin{tikzpicture}
  \matrix (m) [matrix of math nodes,row sep=2em,column sep=2em,minimum width=2em]
  {
    & {(x/\!/G_x,V)} &  \\ \mathrm{pt}& &\cG\\ };
  \path[-stealth]
    (m-1-2) edge node [left] {$\pi\ $} (m-2-1)
        edge node [right] {\ $\iota$} (m-2-3);
\end{tikzpicture}\label{MUdiscr1}\end{equation}
where $\pi$ sends {$x$ to pt and $G_x$ to 1}, and where $\iota$ {embeds $x/\!/G_x$ into} 
the component $o$ of $\cG$. {As $o$ runs over the different components of $\cG$ and
$V$ runs over all simple $G_x$-modules, we recover the canonical basis of $KK(\mathrm{pt},\cG)$.}
The $KK$-group $KK(\cG,\cG)$ becomes the direct sum $\oplus_{o,o'}\mathrm{Hom}(R_{G(o)},R_{G(o')})$, one summand for each orbit $o$ (where Hom here denotes the space of linear maps between vector spaces $R_{G(o)}$ and $R_{G(o')}$).

Our main examples in this paper are groupoids over finite groups. However, in section 6, as well as isolated places in sections 4 and 7, we consider cases where $G$ is a compact connected Lie group. This involves a more elaborate treatment, which we will describe in those subsections.

\subsection{Chiral and full data of a CFT}

The chiral data of a conformal field theory (CFT) can be identified with a \textit{vertex operator algebra} (VOA) or, what should be the same thing (at least in the unitary setting), a \textit{conformal net of factors} on $S^1$. The only place we refer to these is section 4.2, where we discuss lattice theories, and section 5.4, where we discuss their $\bbZ_2$-orbifolds. For an introduction to their complicated theory, see e.g. \cite{LL,CKLW}.

A \textit{fusion category} is a semisimple abelian rigid monoidal category with irreducible  tensor unit.
A \textit{modular tensor category} $\cC$ is a  fusion category with a 
braiding which is maximally non-degenerate in a certain sense. The category of representations Rep$(\cV)$ of a rational VOA or conformal net $\cV$ (such as the lattice theories and their $\bbZ_2$-orbifolds) will be a modular tensor category.

 Let $\Phi$ denote the (finite) set of {isomorphism classes $\lambda$ of} irreducible objects in $\cC$; we call these \textit{primaries}.
The Grothendieck ring of $\cC$ is called the \textit{Verlinde} or \textit{fusion ring} Fus, and has 
basis $\Phi$.  A modular tensor category has an associated unitary representation (unique up to a cube root of 1) of
the modular group SL$_2(\bbZ)$ on the complexification $\bbC\otimes_\bbZ\mathrm{Fus}$, called the \textit{modular data}.
Now, SL$_2(\bbZ)$ is generated by  $\left({0\atop 1}{ -1\atop  \,0}\right)$ and $\left({1\ 1\atop 0\ 1}
\right)$, so this representation is uniquely determined by the matrices $S,T\in M_{\Phi\times\Phi}(\bbC)$ corresponding to those generators. $T$ is a diagonal matrix, whereas $S$ determines the structure constants $N_{\lambda,\mu}^\nu$ of Fus through Verlinde's formula
\begin{equation}\label{Fusl}N_{\lambda,\mu}^\nu=\sum_\kappa S_{\lambda,\kappa}\,\frac{S_{\mu,\kappa}}{S_{\mathbf{0},\kappa}}\,\overline{S_{\nu,\kappa}}\end{equation}
where here and elsewhere we use $\mathbf{0}$ to denote the isomorphism class of the tensor unit.

A full CFT associated to a given rational VOA or conformal net $\cV$, consists of two local extensions $\cV_+$ and $\cV_-$, and a braided equivalence Rep$(\cV_+)\to \mathrm{Rep}(\cV_-)$. These extensions $\cV_+,\cV_-$ are called the type 1 parents. 
Write  $b_\pm$ for the restriction (branching rules) from   $\cV_\pm$-modules to $\cV$. Then 
associated to the full CFT we have the matrix
\begin{equation}\cZ=b_-\sigma b_+^t\,,\label{modinvsigma}\end{equation}
 where  $\sigma$ is a permutation matrix corresponding to the braided equivalence mentioned above.

\medskip\noindent\textbf{Definition 1.} \textit{A matrix $\cZ=(
\cZ_{\lambda,\mu})_{\lambda,\mu\in\Phi}$ is called a} modular invariant
\textit{if $\cZ S=S\cZ$, $\cZ T=T\cZ$, each entry $\cZ_{\lambda,\mu}$
is a nonnegative integer, and $\cZ_{\mathbf{0},\mathbf{0}}=1$.}\medskip

For example, $\cZ=I$ is always a modular invariant. The matrix $\cZ$ in \eqref{modinvsigma} coming from a full CFT is a modular invariant. Conversely, a given modular invariant may come from many, one, or no full CFTs. When it comes from at least one, we call it \textit{sufferable}. For instance, $\cZ=I$ is always sufferable. The modular invariant is a combinatorial shadow cast by the full CFT. Our real interest of course is in
this full structure, although in this paper we only really use this full structure in sections 4.4, 5.4, 5.5 and 6.2.

This full structure was first captured mathematically in the language of factors  \cite{BE1,BE4}, and then axiomatised  in the language of module categories  \cite{Ost1}. A \textit{module category} $\cM$ over a fusion category $\cC$
consists of a bifunctor $\otimes:\cC\stimes\cM\rightarrow\cM$ (corresponding to the \textit{nimrep}) together with compatible associativity and unit
isomorphisms. 
There are obvious notions of equivalence and direct sums
of module categories, and of indecomposable module categories. Given any module category $\cM$ over $\cC$,  there is an algebra $A$ in $\cC$ such that the category Mod$_\cC(A)$ of right $A$-modules in $\cC$ is equivalent to $\cM$.

Now consider $\cC$ a modular tensor category. The category Bimod$_\cC(A)$ of $A$-$A$-bimodules in $\cC$ is a fusion category called the \textit{full system}. There are tensor functors $\alpha_\pm$ from $\cC$ into Bimod$_\cC(A)$ called the \textit{alpha-inductions}.  The modular invariant \eqref{modinvsigma} associated to the module category $\cM$ is \cite{BEK1}\begin{equation}\label{modinvalpha}\cZ_{\lambda,\mu}=\mathrm{dim}\,\mathrm{Hom}_{\mathrm{Bimod}_\cC(A)}( \alpha_+(\lambda),\alpha_-(\mu))\end{equation}
The intersection of the images $\alpha_\pm(\cC)$ is called the \textit{ambichiral system}, the modular tensor category shared by the two type 1 parents of $\cM$. By \textit{sigma-restriction} we mean the forgetful functor from the full system to $\cC$; restricted to the ambichiral system, it recovers the branching rules familiar to CFT.

Type 1 module categories (i.e. those of pure extension type) in a modular tensor category $\cC$ correspond to rigid commutative algebras $A$ with trivial twist. In this case, $\alpha_\pm$ are equivalent (though non-equal); they can be recast as a tensor functor from $\cC$ to Mod$_\cC(A)$, with the latter now regarded as the full system. The local (or dyslectic) $A$-modules form the ambichiral system. The modular invariant will be block-diagonal. The other extreme are the type 2 module categories (i.e. those of pure automorphism type), where the modular invariant is a permutation matrix. As \eqref{modinvsigma} suggests, any (indecomposable) module category will be a combination of two type 1's linked by a type 2.
 
 Any full rational CFT is expected to correspond in this way to
a module category. This framework captures its boundary data, defect lines, spaces of conformal
blocks in arbitrary genus, etc (see e.g. \cite{RFFS}).

\subsection{{Two lemmas}}

Throughout this paper, we let $\widehat{G}$ denote the 1-dimensional representations of a group $G$. For $G$ abelian, we call $\gamma:G\times G\to\bbT$ a \textit{pairing} (or bicharacter) if $$\gamma(gh,k)=\gamma(g,k)\,\gamma(h,k)\,,\ \ \gamma(g,hk)=\gamma(g,h)\,\gamma(g,k)\ \ \forall g,h,k\in G\,.$$  A pairing is called \textit{non-degenerate} if
 $\lan g,h\ran=1$ for all $h\in G$ implies $g=0$. Equivalently, a pairing is a group homomorphism $G\to\widehat{G}$, and is non-degenerate when this is an isomorphism. We call a pairing \textit{symmetric} if $\gamma(g,h)=\gamma(h,g)$ $\forall g,h\in G$, and \textit{alternating} if $\gamma(g,h)=\overline{\gamma(h,g)}$ $\forall g,h\in G$. Symmetric non-degenerate pairings, which exist for any finite abelian group, are central to section 4.

Given any (not necessarily abelian) finite group $G$ and 2-cocycle $\psi\in Z_G^2(\mathrm{pt};\bbT)$, write 
$\beta^{[\psi]}$
for the element of $\prod_{g\in G}\widehat{C_G(g)}$ with components $\beta^{[\psi]}_g$ {defined by}
 \begin{equation}\label{betag}\beta^{[\psi]}_g(h):=\psi(g,h)\,\overline{\psi(h,g)}=\overline{
\beta^{[\psi]}_h(g)}\end{equation}
for all $h\in C_G(g)$, the centraliser of $g$ in $G$. Indeed, the 2-cocycle condition for $\psi$ directly yields
\begin{equation} \beta^{[\psi]}_g(hk)=\beta^{[\psi]}_g(h)\,\beta^{[\psi]}_g(k)\,,\label{bichar2} \end{equation}
for all $g\in G$ and $h,k\in C_G(g)$, so in particular $\beta^{[\psi]}_g$ lies in $\widehat{C_G(g)}$. 
{It is now easy to verify that}  this map $[\psi]\mapsto \beta^{[\psi]}$ yields a well-defined group homomorphism $H^2_G(\mathrm{pt};\bbT)\rightarrow \prod_{g\in G}\widehat{C_G(g)}$. 
When $G$ is abelian, $\beta^{[\psi]}$ is an alternating pairing. The second sentence of part (b) is implicit
 in the proof of Theorem 2 in  \cite{Hugh}, but our proof --- an immediate consequence of the much more
 general part (a) --- seems to be new.

\medskip\noindent\textbf{Lemma 1.} {\textit{Let $G$ be any finite group.}}\smallskip

\begin{itemize} \item[(a)] \textit{{This homomorphism {$H^2_G(\mathrm{pt};\bbT)\rightarrow \prod_{g\in G}\widehat{C_G(g)}$},
$[\psi]\mapsto \beta^{[\psi]}$, is one-to-one.}}

\item[(b)] \textit{Suppose in addition that $G$ is abelian. Then the map $[\psi]\mapsto \beta^{[\psi]}$
is an isomorphism from $H^2_G(\mathrm{pt};\bbT)$ to the group $\cA\cP(G)$ of alternating pairings $\gamma$ on $G$ ($\gamma(g,h)=\beta^{[\psi]}_g(h)$). 
Moreover, the class  $[\psi]\in H^2_G(\mathrm{pt};\bbT)$ contains a cocycle $\gamma\in\cA\cP(G)$,
iff $[\psi]$ is a square, i.e. iff $[\psi]=[\psi']^2$ for some $[\psi']\in H^2_G(\mathrm{pt};\bbT)$, in which
case $\gamma=\beta^{[\psi']}$ is that cocycle.}
\end{itemize}

\noindent\textit{Proof.} {We begin by proving part (a) for $p$-groups. Let $p$ be any prime. Suppose for
contradiction that there exists a $p$-group $G$ and a nontrivial class  $[\psi]\in H_G^2(\mathrm{pt};\bbT)$ such that $\beta^\psi(g,h)=1$ whenever $g,h\in G$ commute. Without loss of generality we may
assume $G$ has minimal order amongst such counterexamples. Being a $p$-group, $G$ has nontrivial
centre $Z=Z(G)$. Because $Z\le C_G(g)$ for all $g\in G$, we can restrict $\beta^{[\psi]}\in\prod_g\widehat{C_G(g)}$
to a group homomorphism $\widetilde{\beta^{[\psi]}}:G\rightarrow \widehat{Z}$, namely by $\widetilde{\beta^{[\psi]}}(g)(z)
=\beta_g^{[\psi]}(z)$. But $\widetilde{\beta^{[\psi]}}$ is identically 1, by hypothesis. By Theorem 4.3 of \cite{Rea}, there exists a 2-cocycle $\psi'$ cohomologous to $\psi$, such that $\psi'(zg,z'g')=\psi'(g,g')$ for all
$g,g'\in G$, $z,z'\in Z$. This means $\psi'$ descends to a 2-cocycle on $G/Z$. But $\beta^{[\psi']}=\beta^{[\psi]}$, so $\psi'\in Z^2_{G/Z}(\mathrm{pt};\bbT)$ has $\beta^{[\psi']}$ identically 1. Because $G$ is by
hypothesis a minimal counterexample (and $Z\ne 1$), this means $[\beta']=1$ as a 2-cocycle of 
$G/Z$, which means it is  also trivial as a  2-cocycle of $G$. Hence $[\psi]$ itself is trivial, as desired.

Now we can do the general case of (a). Let $G$ be any finite group, and $[\psi]\in H^2_G(\mathrm{pt};\bbT)$
is nontrivial. Pick any prime $p$ dividing the order of $[\psi]$. Then $p$ will divide $|G|$, so let $P$
denote a Sylow $p$-subgroup of $G$. By transfer (see e.g. Theorem III.10.3 of \cite{Brown}), the
restriction res$^G_P$ maps the $p$-primary part of $H^2_G(\mathrm{pt};\bbT)$ isomorphically onto
the $G$-invariant elements of $H^2_P(\mathrm{pt};\bbT)$. In particular, the restriction of $[\psi]$ to
$P$ is nontrivial. By the argument of the previous paragraph, this means the restriction of $\beta^{[\psi]}$
to $P$ will be nontrivial, hence $\beta^{[\psi]}$ certainly can't be identically 1 on $G$, and we're done.}

\smallskip
{Now turn to the proof of part (b).} {Suppose $G$ is abelian. We may write $G=\langle g_1\rangle
\times \cdots\times \langle g_t\rangle\cong\bbZ_{n_1}
\times\cdots\times\bbZ_{n_t}$ where $n_t|n_{t-1}|\cdots|n_1$ are the orders of $g_t,\ldots,g_1$ resp.
Then} $H^2_G(\mathrm{pt};\bbT)
\cong \prod_{i=2}^t\bbZ_{n_i}^{i-1}$, thanks to K\"unneth. This  is isomorphic to the group $\cA\cP(G)$, as $\gamma\in\cA\cP(G)$ is uniquely determined by its values $\gamma(g_i,g_j)=\overline{\gamma(g_j,g_i)}$
for $i<j$,
where $g_i\in G$ is a generator of the subgroup $\bbZ_{n_i}$, and that value will be
an $n_j$th root of unity (since $n_j|n_i$). But for any 2-cocycle $\psi$, $\gamma^\psi(g,h):=\beta^{[\psi]}_g(h)$ lies in $\cA\cP(G)$, thanks to 
\eqref{bichar2}. Thus
the map $[\psi]\mapsto\beta^{[\psi]}$, which is injective by part (a), must indeed be an isomorphism onto $\cA\cP(G)$.

In particular, any pairing $\gamma\in\cA\cP(G)$ is a normalised 2-cocycle. We need to
 determine when $\gamma,\gamma'\in\cA\cP(G)$
 lie in the same cohomology class. {That is, we must identify the  functions} $f:G\rightarrow\bbT$ 
satisfying $$ f(g^2)=f(g)^2\,,\ f(gh)^2=f(g)^2f(h)^2\,,\ f(ghk)\,f(g)\,f(h)\,f(k)=
f(gh)\,f(gk)\,f(hk)$$
and $f(g_i)=1$ (we can assume $f(g_i)=1$, as dividing $f$ by the character
$\phi\in\widehat{G}$ matching the values $f(g_i)$ {does not change the associated 2-coboundary}). Then $f^2\in\widehat{G}$ and is
trivial on all generators. Thus $f:G\rightarrow \{\pm 1\}$.

It is elementary to verify that such an $f$ is uniquely determined by the values $f_{ij}:=f(g_ig_j)$
for $i<j$: indeed, $f(g^2h)=f(h)$ and $f(g_{l_1}\cdots g_{l_k})=\prod_{1\le i<j\le k}
f_{l_il_j}$. Moreover, {for any $1\le i<j\le t$ define $f^{(ij)}:G\rightarrow \bbT$ by $f^{(ij)}_{ij}=-1$ and
all other $f^{(ij)}_{i'j'}=+1$; then} it is easy to verify that $f^{(ij)}$
is a solution. This means that the number of solutions
$f$ is precisely $2^{\left({L\atop 2}\right)}$. It also means that the number of $\gamma\in\cA\cP(G)$ lying in
cohomologically distinct classes is precisely equal to the number of classes in $H^2_G(\mathrm{pt};\bbT)$ which are perfect squares. That each one of these classes does intersect $\cA\cP(G)$, is clear from the map $\gamma^\psi(g,h)$ defined above. Indeed, it is trivial that $\gamma^\psi\in\cA\cP(G)$ for any 2-cocycle $\psi$, and that 
if $[\psi]=[\psi']$, then $\gamma^\psi=\gamma^{\psi'}$ as functions on $G\times G$. Likewise, it is elementary that, for finite abelian $G$, 2-cocycles $\psi(g,h)$ and $\overline{\psi(h,g)}$ lie in the same class.

It suffices to show that any 2-cocycle $\psi$ with $\gamma^\psi$ identically 1, 
must be coboundary. We have $\psi(g,h)=\psi(h,g)$, and any cohomologous 
cocycle will likewise be symmetric. Since $H^2_{\bbZ_n}(\mathrm{pt};\bbT)=1$, we may assume
$\psi(g_i^a,g_i^b)=1$ for all $i,a,b$ (this only involves the values $f(g_i^c)$, so we can do
each $\langle g_i\rangle$ separately). Now, assume there is a $k\ge 1$ such that $\psi(g_{1}^{a_1}\cdots g_{{k}}^{a_{k}},
g_{1}^{b_1}\cdots g_{{k}}^{b_{k}})=1$ for all $a_i,b_j$. Indeed, $k=1$ works. Induce on $k$. By 
defining $f(g_{1}^{a_1}\cdots g_{k}^{a_k} g_{k+1}^b)$
appropriately,  we can require $\psi(g_{1}^{a_1}\cdots g_{k}^{a_k},
g_{k+1}^b)=1$. The cocycle condition on $\psi$,
for $x=g_{1}^{a_1}\cdots g_{k}^{a_k},y=g_{1}^{b_1}\cdots g_{{k}}^{b_{k}},z=g_{k+1}^c$ then
says $\psi(g_{1}^{a_1}\cdots g_{k}^{a_k},g_{1}^{b_1}\cdots g_{{k}}^{b_{k}}g_{k+1}^c)=1$,
while the choice $x=g_{1}^{a_1}\cdots g_{k}^{a_k},y=g_{k+1}^b,z=g_{k+1}^c$  says
$\psi(g_{1}^{a_1}\cdots g_{k}^{a_k}g_{k+1}^b,g_{k+1}^c)=1$. From these, the choice
$x=g_{1}^{a_1}\cdots g_{k}^{a_k},y=g_{k+1}^b,z=g_{1}^{c_1}\cdots g_{{k+1}}^{c_{k+1}}$ 
then establishes the induction hypothesis for $k+1$. Thus if $\gamma^\psi$ is identically 1,
 $\psi$ must be coboundary. 
 \textit{QED to Lemma 1}

\medskip The following
simple observation will also be useful.

\medskip\noindent{\textbf{Lemma 2.} \textit{Let $J_1,J_2$ be finite abelian groups, and $\epsilon:J_1\rightarrow\widehat{J_2}$
be a pairing. Then there exists a subgroup $J_0\le J_2$ such that $\phi\in\widehat{J_2} $ is in the image of $\epsilon$ iff
ker$\,\phi$ contains $J_0$. Moreover, $J_1/\mathrm{ker}\,\epsilon\cong J_2/J_0$.}}\medskip 
 
\noindent\textit{Proof.}  Fix any generators $h_i$ resp.\ $h_j'$ of $J_1$ resp.\ $J_2$. Define the (rectangular) matrix $E$ with entries $E_{j,i}\in\bbQ$  defined  by $\epsilon(h_i)(h'_j)
 =e^{2\pi i E_{j,i}}$. The entries are only uniquely determined mod 1, but that isn't important. Then, using Smith normal form over $\bbZ$, there exist integer matrices $P,Q$ invertible over $\bbZ$ such that $PEQ$ is
 diagonal, with entries $r_1,r_2,\ldots,r_t\in\bbQ$. What matters are the denominators $n_i$ of the $r_i$ (when the $r_i$ are written in reduced form). 
 Then the image of  $\epsilon$  is $\cong \bbZ_{n_1}\times \bbZ_{n_2}\times\cdots$. To see the order $n_k$ generator, let $u$ be the $k$th column of $Q$, and $v$ the $k$th row of $P$. Then $\epsilon(\sum_iu_ih_i)(\sum_jv_jh'_j)=e^{2\pi i r_k}$, which has order $n_k$ as desired. 

Define  $\tilde{\epsilon}:J_2\rightarrow\widehat{J_1}$ by
 $\tilde{\epsilon}(g)(h)=\epsilon(h)(g)$. The treatment for $\tilde{\epsilon}$ is similar. Its matrix $\tilde{E}$, and its Smith normal form, is the transpose of that for $E$. Hence the image of $\tilde{\epsilon}$ is also $\cong \bbZ_{n_1}\times \bbZ_{n_2}\times\cdots$, for the same reason.
Therefore $J_2/\mathrm{ker}\,\tilde{\epsilon}\cong J_1/\mathrm{ker}\,\epsilon$.  We find that $J_0=\mathrm{ker}\,\tilde{\epsilon}$, the common kernel of all the
 $\epsilon(g)$'s, works. \textit{QED to Lemma 2}

\subsection{Categories of bundles over groupoids}

Consider a groupoid $X/\!/G$, where both $G$ and $X$ are finite. By a bundle $V$ over $X/\!/G$, or a  
 $G$-equivariant vector bundle over $X$, we mean a choice of vector space $V_x$ over each $x\in X$, and an action of $G$ on the  total space  $V=\oplus_{x\in X}V_x$ such that $g(V_x)=V_{g.x}$ and
 $(g_1g_2)v=g_1(g_2v)$.  We write  $V\in K^0_G(X)$. Here, $ K^1_G(X)=0$. Equivalently, a bundle  $V\in K^0_G(X)$ is a choice of vector space $V_x$ for each 
 orbit representative $x\in X$, such that $V_x$ carries a representation of the stabiliser of $x$ in $G$. 

Given a 2-cocycle $\psi\in Z_G(\mathrm{pt};\bbT)$, by a $\psi$-twisted bundle $V$ over $X/\!/ G$ we mean  $(g_1g_2)v=\psi(g_1,g_2)g_1(g_2v)$, or equivalently to each orbit representative $x\in X$, $V_x$ carries a projective representation of the stabiliser of $x$ in $G$, with cocycle given by $\psi$. 
By a slight abuse of notation, we speak of bundles over $X/\!/_\psi G$
rather than $\psi$-twisted bundles over $X/\!/G$. We write  $V\in{}^\psi K^0_G(X)$. Again, ${}^\psi K^1_G(X)=0$.

 Direct sums of bundles are defined in the obvious way. A bundle is indecomposable iff $V_x=0$ for all but one orbit, and for that orbit,  $V_x$ is an irrep of the stabiliser.

A morphism $V \to W$ between bundles $V,W\in{}^\psi K^0_G(X)$ is a set of linear maps $\oplus_{x\in\cO}V_x \to \oplus_{x\in\cO}W_x$ for each orbit $\cO\subset X$ which commute with the $G$ action. 
We write Bun$(X/\!/_\psi G)$ for the category of bundles over the groupoid $X/\!/_\psi G$. When a twist $\psi$  is identically 1, or a group is trivial, we usually drop it from the notation.

Defining a multiplicative structure on the bundles is much more delicate. Here is a basic approach. First, we have the natural identification $K^0_G(X)\times K^0_G(X)\rightarrow K^0_{G\times G}(X\times X)$, because $K^1$ vanishes ($X$ is finite). Embedding the diagonal subgroup $\Delta_G=\{(g,g)\}$ in $G\times G$ gives the restriction $K^0_{G\times G}(X\times X)\rightarrow K^0_{\Delta_G}(X\times X)$. Now, suppose we have a map $M:X\times X\to X$
which is $G$-equivariant  in the sense that $M(g.x,g.y)=g.M(x,y)$. Suppose also we choose any subset $Y\subset X\times X$ stable under the action of  $\Delta_G$. Then this gives a wrong-way map $M_!:K^0_{\Delta_G}(X\times X)\rightarrow K^0_G(X)$ in $K$-theory, sending some bundle $V$ on $X\times X$ to the bundle on $X$ whose fibre at $x\in X$ is $\oplus_{(x',x'')\in Y\cap M^{-1}(x)}V_{(x',x'')}$. Composing all three maps gives a product   $K^0_G(X)\otimes K^0_G(X)\rightarrow K^0_{G}(X)$ of bundles. We give examples of this in section 5.5.

More generally, consider any $X/\!/G$ and any $X'/\!/_\psi H$.  Suppose we have a group homomorphism $\phi:H\to G$ and  a function $M':X\times X'\to X'$ which is $H$-equivariant in the sense that $M'(h.x,h.y)=h.M'(x,y)$ (where $H$ acts on $X$ through $\phi$). Take $Y'$ to be any subset of $X\times X'$ which is stable under the diagonal $H$ action. Then the  bundles over $X/\!/G$ multiply those over $X'/\!/_\psi H$, in an obvious modification of the previous argument.  We get  a product   $K^0_G(X)\otimes {}^\psi K^0_H(X')\rightarrow {}^\psi K^0_{H}(X')$.

Most of the axioms of fusion category and module category are automatically satisfied by these categories of bundles. The exception is associativity of tensor  products. More specifically, for Bun$(X/\!/ G)$ to be a fusion category, we require 
for each triple $X,Y,Z$ of bundles  a specific choice of \textit{associator} $a_{X,Y,Z}:(X\times Y)\otimes Z\to X\otimes(Y\otimes Z)$ in our category, which are required to satisfy
$$(I_X\otimes a_{Y,Z,W})a_{X,Y\otimes Z,W}(a_{X,Y,Z}\otimes I_W)=a_{X,Y,Z\times W} a_{X\otimes Y,Z,W}$$
For Bun$(X'/\!/_\psi H)$ to be a module category over Bun$(X/\!/ G)$, we require an analogous condition. In both cases, it suffices to consider indecomposable bundles. Certainly this has no chance to work unless the $G$-orbits of $M(M(x,y),z)$ and $M(x,M(y,z))$ coincide (with a similar condition for $M'$). Situations where we do get fusion and module categories through this approach are given in section 5.5.

For a simple example, let $G$ be finite (and not necessarily abelian). Then the fusion category Rep$(G)$ is realised as bundles over the groupoid $\mathrm{pt}/\!/G$, where $G$ fixes pt. Here, $M(\mathrm{pt}\times\mathrm{pt})=\mathrm{pt}$ and $Y=\mathrm{pt}\times\mathrm{pt}$. Indecomposable bundles correspond to irreps $V_\phi$ of $G$. The tensor product of bundles $V_\phi\otimes V_\rho$ is simply their usual tensor product $V_{\phi\otimes\rho}$ as $G$-representations, as the $K$-theoretic description mimics the usual coproduct. Associators $a_{\phi,\rho,\chi}$ can be taken to be 1. The indecomposable module categories of Rep$(G)$ are in  bijection with pairs $(H,[\psi])$ where $H\le G$ and $[\psi]\in H^2_H(\mathrm{pt};\bbT)$. These categories can be realised as bundles over the groupoid $\mathrm{pt}/\!/_\psi H$. Here, $M'(\mathrm{pt}\times\mathrm{pt})=\mathrm{pt}$ and $Y'=\mathrm{pt}\times\mathrm{pt}$. We find that for $\phi\in \mathrm{Irr}(G)$ and $\chi\in\mathrm{Irr}_\psi(H)$, $\phi\otimes\chi=\phi|_H\chi$. Again, the associators are all trivial.

\section{Simple-currents}

Amongst the most accessible and important modular invariants, are the simple-current modular invariants. These are always sufferable. In this section we describe these, for any
choice of modular data, in a uniform way that makes the subsequent $KK$ analysis
more straightforward. The result is remarkably similar to the classification of the module categories 
of finite group doubles by Ostrik \cite{Ost}: for the untwisted double of $G$, these are in bijection with pairs $(H,[\psi])$ of subgroups  $H\le G\times G$ and classes $[\psi]\in H^2_H(\mathrm{pt};\bbT)$. There is some overlap of our Theorem 1 with \cite{KrSch}, although our proof
is unrelated.

For an arbitrary modular tensor category, 
it is natural (though probably naive) to hope for some correspondence as in section 2, between  module categories and some sort of cohomological data. This is known to happen for the 1-dimensional lattice theory $\sqrt{2n}\bbZ$ \cite{BE6}: there the modular invariants are in one-to-one correspondence
with the subgroups of $\bbZ_n$. 
In this subsection we generalise this characterisation to all simple-current modular invariants
for all modular data. (Recall that in that lattice case, all modular invariants are simple-current ones.)

\subsection{Simple-currents and modular data}

As always, we let $\widehat{G}$ denote the 1-dimensional representations. By definition, a \textit{simple-current} is an invertible object in a modular tensor category, or equivalently
any primary with Perron-Frobenius dimension 1. The simple-currents form a finite abelian group $\cJ$ of the modular tensor category.
The multiplication in this group is restriction of the fusion product. The terminology
simple-current, and with it the symbol $J$, are rather obscure and come from the
CFT literature. The basic properties of simple-currents we need are
proved in \cite{Ganmd}.

To each simple-current $j\in\cJ$ is associated two things: a permutation of the primaries
$\Phi$ given by the fusion product, and a grading 
of the fusion ring coming from Verlinde's formula. More precisely, we
have a faithful group homomorphism $\cJ\rightarrow \mathrm{perm}(\Phi)$, and 
a surjective grading $\cQ:\Phi\rightarrow \widehat{\cJ}$:
for any primaries $a,b,c\in\Phi$, if $c$ appears in the fusion product of $a$ and $b$,
then $\cQ_a\cQ_b=\cQ_c$. {In particular,
\begin{equation}\label{Qprod}\cQ_{ja}=\cQ_j\cQ_a\,.\end{equation}}
Moreover, we have the reciprocity \cite{Ganmd}
\begin{equation}\label{jsym}\cQ_j(j')=\cQ_{j'}(j)\ \forall j,j'\in\cJ\,.\end{equation}

Simple-currents respect modular data:
\begin{eqnarray} S_{ja,b}/S_{ja,\mathbf{0}}=S_{b,ja}/S_{\mathbf{0},ja}=\cQ_b(j)\,S_{a,b}/S_{a,\mathbf{0}}\,,\label{Sj}\\
T_{ja,ja}\overline{T_{a,a}}=\overline{\cQ_a(j)}T_{j,j}\overline{T_{\mathbf{0},\mathbf{0}}}\,,\label{TTj}\\
 \left(T_{j,j}\overline{T_{\mathbf{0},\mathbf{0}}}\right)^{2}=\overline{\cQ_j(j)}\,,\label{Tj}\end{eqnarray}
for all simple-currents $j\in \cJ$ and primaries $a,b\in\Phi$. \eqref{Sj} can be taken as the
definition of $\cQ_b(j)$. When the order of $j$ is odd,  $T_{j,j}\overline{T_{\mathbf{0},\mathbf{0}}}$ also has odd order. When the modular data comes from a so-called \textit{unitary} CFT
(the case of interest in this paper), the denominators $S_{ja,\mathbf{0}}$ and $S_{a,\mathbf{0}}$
in \eqref{Sj} are equal and can be dropped; in nonunitary theories such as the Yang-Lee model,
they contribute a sign.

For example, the group $\cJ(G/\!/G)$ of simple-currents for the double of a finite group $G$ consists of
all pairs $(z,\psi)$ where $z$ is in the centre of $G$ and $\psi$ is a 1-dimensional
representation of $G=C_G(z)$: 
$\cJ\cong Z(G)\times \widehat{G/[G,G]}$. The action on $\Phi$ is $(z,\psi)(g,\chi)=(zg,\psi\chi)$.
The grading is $\cQ_{(g,\chi)}(z,\psi)=\chi(z)\psi(g)/\chi(1)$, and $T_{(z,\psi),(z,\psi)}
=\overline{\psi(z)}$.

For another example, every primary for the loop groups of tori twisted by
an even lattice $L$, are simple-currents:
$\cJ=\Phi=L^*/L$, with group operation being addition, $\cQ_{[\lambda]}([\mu])=e^{-2\pi\i
\lambda\cdot\mu}$ and $T_{[\lambda],[\lambda]}\overline{T_{[0],[0]}}=e^{\pi\i\lambda\cdot
\lambda}$. For loop groups of compact simply-connected groups $G$, the group of simple-currents is isomorphic to the centre $\cZ(G)$ (with one source of exceptions: $E_8$ at level 2).

Incidentally, all abelian groups can arise as a group of simple-currents, even among the fusion rings 
of simply-connected loop groups:  e.g.
 $Z(SU(n_1)\times\cdots \times SU(n_t))\cong \bbZ_{n_1}\times\cdots\times\bbZ_{n_t}$.

For each $j,j'\in\cJ$, write $q(j)=T_{j,j}\overline{T_{\mathbf{0},\mathbf{0}}}$ and $\langle j,j'\rangle=Q_j(j')$. Together, \eqref{Qprod} and \eqref{jsym} say $\langle,\rangle$ is a symmetric pairing. \eqref{TTj},\eqref{Tj} then say $q(-j)=q(j)$ and $\langle j,j'\rangle=q(j)q(j')\overline{q(j+j')}$. We define a quadratic form in section 4.1; this $q(j)$ satisfies all properties of a quadratic form except possibly the non-degeneracy of $\langle,\rangle$.

\medskip\noindent\textbf{Proposition 1.} \cite{Ganmd} \textit{Fix any modular data. Let $\cJ$ be its group of simple-currents, and let $\cZ$ be any modular invariant. Let $j,j'\in\cJ$ be such that $\cZ_{j,j'}\ne 0$. Then
 $\cZ_{ja,j'b}=\cZ_{a,b}$ for all $a,b\in\Phi$, and  $\cZ_{a,b}\ne 0$ only when $\cQ_j(b)=\cQ_{j'}(a)$.}
 \medskip

\subsection{Simple-current modular invariants}

Call $j\in\cJ$ \textit{quaternionic} if $\cQ_j(j)$ has the same order in $\bbC^\times$ as $j$ has in
$\cJ$, and this order is {even}. Equivalently (using \eqref{Tj}), $j\in \cJ$ is quaternionic iff
$(T_{j,j}\overline{T_{\mathbf{0},\mathbf{0}}})^{\mathrm{ord}(j)}= -1$ (that power is always either $\pm 1$). We want to avoid these, for this reason (though they do play a role in supersymmetric theories):

\medskip\noindent\textbf{Lemma 3.} \textit{Let $J$ be a subgroup of the group of simple-currents. Then $J$ contains no  quaternionic elements, iff there is a pairing $\epsilon$ for $J$ satisfying
\begin{eqnarray}\label{Qeps}\cQ_j(j')\,\epsilon_j(j')\,\epsilon_{j'}(j)=1\ \ \ \forall j,j'\in J\,,\\
\label{Teps}\epsilon_j(j)=T_{j,j}\overline{T_{\mathbf{0},\mathbf{0}}}\ \ \ \forall j\in J\,.\end{eqnarray}
The set of all such $\epsilon$ is in bijection with $H^2_J(\mathrm{pt};\bbT)$ (provided $J$ contains no quaternionic elements).}\medskip

\noindent\textit{Proof.} One direction is clear: if $j\in J$ is quaternionic, then the order of $T_{j,j}\overline{T_{\mathbf{0},\mathbf{0}}}$ doesn't divide that of $j$, so $\epsilon$ cannot be a pairing. 

Now suppose $J$ contains no quaternionic elements. Write $J=\langle h_1\rangle\times\cdots\times\langle h_s\rangle$, where $n_s|n_{s-1}|\cdots|n_1$, the orders of $h_s,h_{s-1},\ldots,h_1$. Write $\langle j,j'\rangle=\cQ_j(j')$ and $q(j)=T_{j,j}\overline{T_{\mathbf{0},\mathbf{0}}}$ as in section 3.1. For any $j,j'\in\cJ$ write $j=\sum_ij_ih_i$ and $j'=\sum_ij'_ih_i$, for $j_i,j'_i\in\bbZ$. Define
\begin{equation}\epsilon^J_j(j')=\prod_{i=1}^sq(h_i)^{j_ij'_i}\prod_{i=1}^{s-1}\langle \pi_i (j),j'_ih_i\rangle\label{epsJ}\end{equation}
where $\pi_i(j)=\sum_{k>i}j_kh_k$. This is manifestly a pairing, thanks to the absence of quaternionic elements in $J$. Moreover, it satisfies \eqref{Qeps},\eqref{Teps} by the properties of $\langle,\rangle$ and $q(j)$ mentioned in section 3.1. So there is at least one solution.

Note that $\epsilon$ is a pairing on $J$ satisfying \eqref{Qeps},\eqref{Teps}, iff $\epsilon\,\overline{\epsilon^J}$ is an alternating pairing on $J$. By Lemma 1(b), we get the bijection  $\epsilon=\beta^{[\psi]}\epsilon^J$ with $H^2_J(\mathrm{pt};\bbT)$.  \textit{QED to Lemma 3}\medskip

More precisely, the set of all such $\epsilon$ forms a torsor over $H^2_J(\mathrm{pt};\bbT)$ --- the choice of base-point $\epsilon^J$ is arbitrary. The relevance of the equations \eqref{Qeps},\eqref{Teps} will become clear in the proof of Theorem 1.

For $a\in\Phi$, write stab$(a)=\{j\in\cJ\,:\,ja=a\}$.
Given a subgroup $J\le\cJ$, we call $a\in \Phi$ \textit{$J$-free} 
if stab$(a)=\{0\}$, i.e. if $\|J a\|=|J|$. 

 Call the modular data \textit{sufficiently nonzero} if {for any $\psi,\psi'\in\widehat{\cJ}$ there exists
$a,b\in\Phi$ with $\cQ_a=\psi$, $\cQ_b=\psi'$, and $S_{a,b}\ne 0$. The modular data is sufficiently nonzero
if  for any $a\in\Phi$ and $\psi\in\widehat{\cJ}$
 with stab$_\cJ(a)\le \mathrm{ker}\,\psi$,
 there exists a $b\in \Phi$ $\cJ$-free such that $\cQ_b=\psi$ and $S_{a,b}\ne 0$.
(If instead stab$_\cJ(a)$ was not a subgroup of ker$\,\psi$, then any $b\in\Phi$ with $\cQ_b=\psi$
would necessarily have $S_{a,b}=0$, thanks to \eqref{Sj}.) To see why this implies the sufficiently nonzero hypothesis, apply it twice: firstly to $a'=\mathbf{0}$ to get a  $\cJ$-free $b'\in\Phi$ with $\cQ_{b'}=\psi$, and secondly
to $a=b'$ to get $b\in\Phi$ with the desired properties. }

This notion {of sufficiently nonzero is needed for the strongest results. As stated, it} is  stronger than we need, but {it is hard to find any examples of modular data which is not sufficiently nonzero.}
For example, the only $SU(n)$ level $k$ modular data which is not sufficiently nonzero,
is $SU(2)$ level 2. Torus modular data is always sufficiently nonzero, as are the untwisted doubles
of abelian groups {and all (unitary or nonunitary) minimal models except for the Ising}.

In the following theorem, we parametrise modular invariants associated to simple-currents,
by subgroups $J$ of simple-currents and 2-cocycles $\psi$. More precisely, the modular invariants (in fact, module categories) associated to $J$ form a torsor over $H^2_J(\mathrm{pt};\bbT)$.
 
Parts (a) and (b) in the following theorem are the important cases for us;  {their proof is essentially the same as that of} part (c),
which perhaps is of independent interest. It cannot be a coincidence  that $H^2_J(\mathrm{pt};\bbT)\cong H^3_J(\mathrm{pt}
;\bbZ)$ and $ H^1_J(\mathrm{pt};\bbZ_2)\cong J/2J$ are the standard twists for the $K$-group $K_J^\star(\mathrm{pt})$.
Compare Example 9.7.2 of \cite{EGNO}, which gives module categories over Vect$^\omega_G$, and \cite{KrSch},
which concerns simple-current modular invariants in rational conformal field theory.

\medskip\noindent\textbf{Theorem 1(a)} \textit{Let $J$ be any subgroup of $\cJ$ containing
no quaternionic elements, and choose any class $[\psi]\in H^2_{J}(\mathrm{pt};\bbT)$. Define the pairing $\epsilon=\beta^{[\psi]}\epsilon^J$ on $J$, where $\beta^{[\psi]}$ resp.\ $\epsilon^J$ is defined in \eqref{bichar2} resp.\ \eqref{epsJ}. Let $J_0\le J$ {be as in 
Lemma 2.} Then  there is a modular invariant $\cZ={\cZ(J,[\psi])}=\cZ(J,\epsilon)$
whose nonzero entries are}
 \begin{equation}\label{scinvform}
 \cZ_{a,ja}=\frac{|J_0|}{\|J_0a\|}\end{equation}
\textit{for all $j\in J$ and $a\in\Phi$ satisfying $\cQ_a|_J=\epsilon_j$. }

\smallskip\noindent\textbf{(b)} \textit{Any modular invariant constructed in (a) satisfies
\begin{equation}\cZ_{a,b}\ne 0\Rightarrow
b\in \cJ a\,.\label{scinvdef}\end{equation}
Conversely, suppose the modular data is sufficiently nonzero.
Then any modular invariant satisfying \eqref{scinvdef} equals $\cZ(J,[\psi])$ for
 one and only one pair $(J,[\psi])$.}

\smallskip\noindent\textbf{(c)} \textit{There is an analogous classification of all matrices
$\cZ$ with nonnegative integer entries and $\cZ_{\mathbf{0},\mathbf{0}}=1$, which commute with $S$, and triples 
$(J,[\psi],[\phi])$ where $J$ and $\psi$ are as in (a), and $[\phi]\in {H^1_J(\mathrm{pt};\bbZ_2)\cong} J/J^2$ (interpret $[\phi]$ as a homomorphism $J\rightarrow \{\pm 1\}$, and replace $q(h_i)$ in \eqref{epsJ} with $q(h_i)\,[\phi](h_i)$).
}\medskip

\medskip\noindent\textit{Proof.} Begin with part (a). For any subgroup $J$ and pairing  
$\epsilon$ as in the theorem, let $J_L=J_0$ be as in Lemma 2, and
 define $J_R=\mathrm{ker}\,\epsilon$. 
We will first show that \eqref{Qeps},\eqref{Teps} imply that the matrix  $\cZ=\cZ(J,\epsilon)$ defined by \eqref{scinvform} commutes with
$S$ and $T$. Commutation with $T$ is trivial: by  \eqref{TTj}, it is equivalent to \eqref{Teps}.  

{To show $\cZ$ commutes with $S$, first observe using \eqref{Qeps} that $J_L$, the common
kernel of all $\epsilon_j$, equals  $\{j'\in J\,:\,\cQ_{j'}|_J=\overline{\epsilon_{j'}}\}$. Moreover, 
if $j\in J$ and $ja=a$ for some $a\in\Phi$, then $\cQ_j$ is identically 1 thanks to \eqref{Qprod}, in which case
$j\in J_L$ iff $j\in J_R$. Thus for all $a\in\Phi$,  $\|J_La\|=\|J_Ra\|$.}
 
{We compute directly from  \eqref{Sj} that}
\begin{eqnarray} (S\cZ)_{a,b}={|J_L|}\left\{\begin{matrix} S_{a,j'b}&\mathrm{if}
\ \cQ_a(J_L)=1\ \mathrm{and}\ \cQ_{j'b}|_J=\overline{\epsilon_{j'}}\ \mathrm{for}\ j'\in J\\
0&\mathrm{otherwise}\end{matrix}\right.\label{SZ}\\  \label{ZS}
(\cZ S)_{a,b}={|J_L|}\left\{\begin{matrix} S_{ja,b}&\mathrm{if}
\ \cQ_b(J_R)=1\ \mathrm{and}\ \cQ_{a}|_J={\epsilon_{j}}\ \mathrm{for}\ j\in J\\
0&\mathrm{otherwise}\end{matrix}\right.\,,
\end{eqnarray}
{since $\|J_Ra\|=\|J_La\|$. From Lemma 2 and \eqref{Qeps} we see that $\cQ_a(J_L)=1$
iff there is a $j\in J$ such that $\cQ_{a}|_J={\epsilon_{j}}$; likewise, $\cQ_b(J_R)=1$
iff there is a $j'\in J$ such that $\cQ_{j'b}|_J=\overline{\epsilon_{j'}}$.} 
{Hence, the commutation
$(S\cZ)_{a,b}=(\cZ S)_{a,b}$ is the trivial $0=0$, unless $\cQ_a(J_L)=\cQ_b(J_R)=1$ in which case it reduces to
$S_{a,j'b}=S_{ja,b}$ for $\cQ_a|_J=\epsilon_j$ and $\cQ_{j'b}|_J=\overline{\epsilon_{j'}}$, which 
holds thanks to \eqref{Qeps},\eqref{Sj}.}

{Now turn to part (b).} 
Consider any modular invariant $\cZ$ satisfying \eqref{scinvdef}, {and assume the sufficiently zero
hypothesis}.
Define $J_L=\{j\in\cJ\,:\,\cZ_{j,\mathbf{0}}\ne 0\}$ and $J_R=\{j\in\cJ\,:\,\cZ_{\mathbf{0},j}\ne 0\}$.
Then it is elementary to see \cite{Ganmd} that $\cZ_{ja,j'b}=\cZ_{a,b}$ for all $j\in J_L,j'\in J_R$ and $a,b\in \Phi$.
This implies $J_L,J_R$ are both groups, and $\cZ_{j,j'}=1$ for all $j\in J_L,j'\in J_R$.
Moreover \cite{Ganmd}, $\cZ_{a,b}\ne 0$ implies $\cQ_a(J_L)=\cQ_b(J_R)=1$, and conversely,
$\cQ_a(J_L)=1$ implies there exist $b\in\Phi$ such that $\cZ_{a,b}\ne 0$.

Note that $T\cZ=\cZ T$ implies that if $\cZ_{a,ha}\ne 0$, then {by \eqref{TTj}}
\begin{equation}T_{h,h}\overline{T_{\mathbf{0},\mathbf{0}}}=\cQ_a(h)\ \,.\label{MTTM1}\end{equation} 
Choose any $a\in\Phi$ with $\cQ_a(J_L)=1$; evaluating $S\cZ=\cZ S$ at $(a,\mathbf{0})$ gives
\begin{equation}\sum_{[h]\in \cJ/J_R}\cZ_{a,ha}=\frac{|J_L|}{\|J_Ra\|}\,.\label{normalise}\end{equation}
This implies {(taking $a=\mathbf{0}$) that} $|J_L|=|J_R|$. It also implies that if such an $a$ is $J_R$-free,
then there exists a unique class $[h_a]\in\cJ/J_R$ such that $\cZ_{a,b}=\delta_{J_R b,J_Rh_aa}$
for all $b$.

 Define 
 \begin{equation}
\label{Jdefgen}J=\langle \{j\in \cJ\,:\,\cZ_{a,ja}\ne 0\ \mathrm{for\ some}\ J_R\mathrm{-free}\ a\in\Phi\}\rangle\,.
\end{equation}
Suppose for contradiction that for some $a\in\Phi$ (not necessarily $J_R$-free), there are $h,h'\in \cJ$ such that both  $\cZ_{a,ha}\ne 0$
and $\cZ_{a,h'a}\ne 0$, for $h'{h}^{-1}\not\in
\langle J_R,\mathrm{stab}_Ja\rangle$. Then there exists some character $\psi$ of ${J/
\langle J_R,\mathrm{stab}_Ja\rangle}$ such that $\psi(h')\ne \psi(h)$. By the
sufficiently nonzero hypothesis, there exists some $\cJ$-free $b\in\Phi$ such that
$\cQ_b|_J=\psi$. For this choice we find by the triangle inequality and \eqref{normalise}
$|(\cZ S)_{a,b}|<|J_L|\,{|S_{ab}|}$, which contradicts $(S\cZ)_{a,b}=|J_L|\,{S_{ab}}$
(since $b$ is $J$-free). Thus for all $a\in\Phi$
with $\cQ_a(J_L)=1$, there is  a unique class $[h_a]\in\cJ/J_R$ such that $\cZ_{a,b}=\frac{|J_L|}{\|J_Ra\|}\delta_{J_R b,J_Rh_aa}$
for all $b$. The argument with left and right interchanged now forces $\|J_La\|=\|J_Ra\|$
for such $a$.
Now $\cZ S=S\cZ$ reduces to $S_{h_aa,h_bb}=S_{a,b}$ whenever $\cQ_a(J_L)=1=\cQ_b(J_L)$.
As usual the sufficiently nonzero hypothesis then leads to
\begin{equation}\label{QQQ}\cQ_{h_a}(h_b)\,\cQ_a(h_b)\,\cQ_b(h_a)=1\,,\end{equation}
 for these $a,b$. {For each $j\in J$, \eqref{Jdefgen} says there exists an $a\in\Phi$ such that
 $\cZ_{a,ja}\ne 0$. Define a map $\epsilon:J\rightarrow\widehat{J}$ by $\epsilon_j=\cQ_a|_J$. This assignment
 $\epsilon_j$ is well-defined (i.e. independent of $a$), because we can use \eqref{jsym} to rewrite \eqref{QQQ}
 as $\cQ_a(h_b)=\overline{\cQ_{h_bb}(h_a)}$, which tells us that if $h_{a_1}J_R=h_{a_2}J_R$, then $\cQ_{a_1}
 =\cQ_{a_2}$. The same equation also tells us this $\epsilon$ is a homomorphism.}
 
{For this choice of $\epsilon$, $\cZ$ is given by \eqref{scinvform}, and \eqref{QQQ},\eqref{MTTM1}
reduce to \eqref{Qeps},\eqref{Teps}. It suffices to show that any $\epsilon$ satisfying \eqref{Qeps},\eqref{Teps}
come from a 2-cocycle $[\psi]$. For this purpose,  let $\gamma(j,j')=\epsilon_j(j')
/\epsilon'_j(j')$ be the quotient of two solutions to \eqref{Qeps},\eqref{Teps}. Then $\gamma \in\cA\cP(G)$, so by Lemma 1(b) does indeed come from a 2-cocycle, and we're done.}

{Now suppose $\cZ(J,\epsilon)=\cZ(J',\epsilon')$. If $J\ne J'$, then there exists $j\in J$, $j\not\in J'$
(or the other way around). By the sufficiently nonzero hypothesis, there exists a $\langle j\rangle$-free $a\in\Phi$
such that $\cQ_a=\epsilon_j$. Then $\cZ(J,\epsilon)_{a,ja}\ne 0$ but $\cZ(J',\epsilon')_{a,ja}= 0$, a contradiction.
Thus $J=J'$. If now $\epsilon_j\ne\epsilon'_j$ for some $j\in J$, then for the same $a$, $\cZ(J,\epsilon)_{a,ja}\ne 0$ but $\cZ(J',\epsilon')_{a,ja}= 0$, another contradiction. This means $\cZ(J,\epsilon)=\cZ(J',\epsilon')$ forces
$J=J',\epsilon=\epsilon'$ as desired. Of course, Lemma 3 tells us $\epsilon=\epsilon'$ is equivalent to
$[\psi]=[\psi']$.}
This completes the proof of part (b).

Finally, turn to part (c). Dropping the condition $\cZ T=T\cZ$ amounts to dropping \eqref{Teps}.
Note that \eqref{Qeps} for $j=j'$ reduces to the square of \eqref{Teps}, thanks to \eqref{Tj}. 
This means that the values of $\epsilon_j(j)$ are determined only up to signs. These signs
are arbitrary, provided $\epsilon$ is a pairing. The possibilities for these signs is then
captured by a homomorphism $J\rightarrow\{\pm 1\}$, i.e. by a class in $H^1_J(\mathrm{pt};\bbZ_2)$.
\textit{QED to Theorem 1}\medskip

By definition, a \textit{simple-current modular
invariant} $\cZ$ is any modular invariant of the form \eqref{scinvform}. The modular invariants \eqref{scinvform} all correspond
to at least one module category (see \cite{FRS}).
Different pairs $J,\phi$ can correspond to the same modular invariant
(e.g. $SU(2)$ at level 2, $J=1$ and $J=\bbZ_2$ both correspond to the identity matrix).
But as explained in the proof, the correspondence is one-to-one if, for every $\psi\in\widehat{\cJ}$, there is a
$\cJ$-free $a\in\Phi$ with $\cQ_a=\psi$. 

For example, consider an (untwisted) finite group double $G/\!/G$. There, the group of simple-currents is
$Z(G)\times \widehat{G/[G,G]}$. No simple-current $j=(z,\psi)$ can be quaternionic here. 
If we take $Z'$ to be a subgroup of $Z(G)$, and an arbitrary $\phi\in Z^2_{Z'}(\mathrm{pt};
\bbT)$, then in Ostrik's parametrisation this corresponds to subgroup $H=\Delta_G(1\times
Z')$ and $\phi$ lifts to $\phi^H$ on $H$ by defining $\phi^H((g_1,g_1z_1),(g_2,g_2z_2))=\phi(z_1,z_2)$ (this satisfies the cocycle condition because $Z'$ is in the centre). On the other hand, if $Z'$
is a subgroup of $\widehat{G/[G,G]}$, and $\phi\in Z^2_{Z'}(\mathrm{pt};\bbT)$, then 
 $H=\Delta_{G'}$ where $Z'=\widehat{G/G'}$, i.e. $G'=\cap_{\psi\in Z'}\mathrm{ker}\,\psi$. 
Incidentally,  when $G$ is abelian and the 3-cocycle $\omega$ is trivial,
Theorem 1 classifies all modular invariants, and we can explicitly recover Ostrik's parametrisation.

The fusion rules obeyed by modular invariants (more precisely, module categories) was first explored in \cite{EP2}. Here we can be very explicit.

 \medskip\noindent\textbf{Corollary 1.} \textit{Given any two modular invariants $\cZ=\cZ(J,[\psi]),\cZ'=\cZ(J',
 [\psi'])$, their matrix product satisfies $\cZ\cZ^{\prime\, tr}=n\cZ(J'',[\psi''])$ for some admissible $(J'',[\psi''])$, where $n=\|\{a\in \Phi\,:\,\mathrm{both}\ \cZ_{\mathbf{0},a}\ne
 0\ \mathrm{and}\ \cZ_{\mathbf{0},a}'\ne 0\}\|$.}
 
\medskip \noindent\textit{Proof.} Let $J_R=\{j\in\cJ\,:\,\cZ_{\mathbf{0},j}\ne 0\}$ and $J_R'=\{j\in\cJ\,:\,
\cZ'_{j,\mathbf{0}}\ne 0\}$ as usual. Then Proposition 1 says $\cZ_{a,jb}=\cZ_{a,b}$ and $\cZ'_{a,j'b}=
\cZ'_{a,b}$ for all $j\in J_R$,
$j'\in J_R'$, and $a,b\in\Phi$. In particular, $J_R,J_R'$ are groups, and $\cZ_{\mathbf{0},j}=1=\cZ'_{\mathbf{0},j'}$ for all $j\in J_R$,
$j'\in J_R'$. 

We compute $(\cZ\cZ^{\prime\, tr})_{\mathbf{0},\mathbf{0}}=\sum_a \cZ_{\mathbf{0},a}\cZ'_{\mathbf{0},a}=|J_R\cap J'_R|$, which equals what we call $n$ above. Moreover, for any
 $a,b\in\Phi$ with $(\cZ\cZ^{\prime\, tr})_{a,b}\ne 0$, there must exist $j_0\in J$, $j'_0\in J'$ such that $\cZ_{a,j_0a}\ne 0$ and
 $\cZ'_{b,j_0a}\ne 0$ where $b=j_0^{\prime-1}j_0a$. In this case, $j_0\in J$ and $j_0'\in J'$ satisfy 
 $\cQ_a=\epsilon_{j_0}$ on $J$ and $\cQ_b=\epsilon'_{j_0'}$ on $J'$. Using \eqref{Qeps}, we can rewrite the latter as $\cQ_a=\overline{\cQ_{j_0}}(\star)\,\overline{\eps'_\star(j'_0)}$ on $J'$.

We learn from the proof of Theorem 1(b) that $J''$ will be the set of all $j\in \cJ$ such that $(\cZ\cZ^{\prime\, tr})_{a,ja}\ne 0$ for some $a$. The previous paragraph shows that such $j$ can always be written as $j_0j_0^{\prime-1}$ for some $j_0\in J,j_0'\in J'$, i.e. $J''\subseteq JJ'$.
So $j_0j_0^{\prime-1}\in J''$, for $j_0\in J$ and $j_0'\in J'$, if there is some $a\in \Phi$ with $\cQ_a$ on 
$J$ resp. $J'$ given as above. Those formulas for $\cQ_a$ can be taken to define $a$. But there is a consistency relation: for all $j\in J\cap J'$,  $\epsilon_{j_0}(j)=\overline{\cQ_{j_0}}(j)\,\overline{\eps'_j(j'_0)}$. Putting this together, we get that
$$J''=\{j_0j_0^{\prime -1}:j_0\in J,\ j_0'\in J',\ \epsilon_\star(j_0)|_{J\cap J'}=\eps'_\star(j_0')|_{J\cap J'}\}\,.$$
Moreover, for any $j_0j_0^{\prime-1},jj^{\prime-1}\in J''$, $$\epsilon''_{j_0j_0^{\prime-1}}(jj^{\prime-1})=
\eps_{j_0}(j)\,\cQ_{j_0}(j')\,\eps'_{j'}(j_0')\,.$$

 Write $J_{R,a}$  and $J'_{R,a}$ for the groups containing $j\in J_R$ resp.\ $j'\in J'_R$ which fix $a$.
 Using \eqref{scinvform} and Proposition 1,
$$(\cZ\cZ^{\prime\, tr})_{a,b}=\|J_{R,a}\|\|J'_{R,a}\|(|J_R\cap J'_R|/\|J_{R,a}\cap J'_{R,a}\|)=n\|J_{R,a}J'_{R,a}\|\,.$$
 Thus all entries of $\cZ\cZ^{\prime\, tr}$ are divisible by $n$, as desired. Thus $\cZ\cZ^{\prime\, tr}=
n\cZ(J'',\epsilon'')$ for $J'',\epsilon''$ defined above. As in the proof of Theorem 1(b), such an $\epsilon''$ corresponds to a unique class $[\psi'']$.
\textit{QED to Corollary 1}\medskip

(This Corollary should be compared with section 3.3 of \cite{KrSch}.) For the special case where $J=J'$ and $\psi=\psi'$, we get $J''=J$ and $\psi''=1$. Note that in general, $\cZ(J,\eps)^{tr}=\cZ(J,\eps^{tr})$, where $\eps^{tr}_j(j')=\eps_{j'}(j)$.

\section{Pointed modular tensor categories}

\subsection{Quadratic forms and pointed categories}\medskip

As before,  $\bbT$ denotes the complex numbers of modulus 1 and $\widehat{G}$ denotes the 1-dimensional irreps of $G$. Write $\xi_m=e^{2\pi i/m}$.

Let $G$ be a finite abelian group (written additively). Recall the definition of non-degenerate symmetric pairing from section 2.4. We will denote these by $\langle\cdot,\cdot\rangle$. 
 Unless otherwise stated, all symmetric pairings in this section are assumed to be non-degenerate.
A \textit{quadratic form} on $G$ is a map $q:G\rightarrow\bbT$ obeying $q(-g)=q(g)$ for all $g\in G$, {such that} the formula
\begin{equation}\label{quadbil}\lan g,h\ran:=q(g)\,q(h)\,\overline{q(g+h)}\end{equation}
defines a non-degenerate symmetric pairing $\lan\cdot,\cdot\ran$ on $G$. We call two quadratic forms $q,q'$ on $G$ \textit{equivalent}, if there is a group automorphism $\alpha$ of $G$ such that $q\circ\alpha=q'$. Lemma 4 explains to what extent the correspondence $\langle,\rangle\leftrightarrow q$ is a bijection.

\medskip\noindent\textbf{Lemma 4.} \textit{Let $\langle\cdot,\cdot\rangle$ be any  non-degenerate symmetric pairing  on $G$.}

\begin{itemize}

\item[(a)] \textit{$\langle\cdot,\cdot\rangle$ has precisely $2^l$ associated quadratic forms, where $2^l$ is the order of the largest elementary 2-subgroup $\bbZ_2\times\cdots\times\bbZ_2$ of $G$.}

\item[(b)] \textit{At most two of those quadratic forms will be inequivalent. All $2^l$ of them will be equivalent, unless there is an order 4 element $g\in G$ with $\langle g,g\rangle^2=-1$.}

\end{itemize}

We prove this at the end of the subsection. We use Lemma 4 in section 5.2.

Given any quadratic form $q$ on $G$, define the matrices
\begin{align}\label{WeilT} T_{g,h}=&\,xq(g)\,\delta_{g,h}\,,\\ \label{WeilS} S_{g,h}=&\,\sqrt{|G|}{}^{-1}\lan g,h\ran\end{align} 
where 
$x^{-1}$ is any cube-root of $\sqrt{|G|}^{-1}\sum_{k\in G}q(k)$. This defines a representation of SL$_2(\bbZ)$ called the \textit{Weil representation}, through $\left({0\atop 1}{-1\atop 0}\right)\mapsto S$ and $\left({1\atop 0}{1\atop 1}\right)\mapsto T$.

The \textit{pointed} modular tensor categories (i.e. those whose simple objects are all simple-currents) are described in section 8.4 of \cite{EGNO}. In particular, 
given a quadratic form $q$ on $G$, there exists a skeletal pointed unitary modular tensor category $\cC(q,G)$ with simple objects $g\in G$ and modular data given by the Weil representation.  Its fusion ring is isomorphic to the group ring $\bbZ[G]$. Here, $x=e^{-\pi ic/12}$ where $c$ is the central charge (unique up to a shift by $8\bbZ$). Any pointed modular tensor category is braided tensor equivalent to some $\cC(q,G)$. Moreover, $\cC(q,G)$ is braided tensor equivalent to $\cC(q',G')$ iff there is a group isomorphism $\varphi:G\rightarrow G'$ {(since the fusion rings are isomorphic)} such that $q=q'\circ \varphi$ (since the $T$-matrices must match). Thus pointed modular tensor categories are identified up to braided tensor equivalence by their modular data.

Let $q_1$ resp.\ $q_2$ be quadratic forms on $G_1$ resp.\ $G_2$. Then $q_1\times q_2$ defined by $(q_1\times q_2)(g_1,g_2)=q_1(g_1)\,q_2(g_2)$ is a quadratic form on $G_1\times G_2$. Moreover, $\cC(q_1\times q_2,G_1\times G_2)$ is braided tensor equivalent to the Deligne product
$\cC(q_1,G_1)\stimes \cC(q_2,G_2)$. An $x$ for $q_1\times q_2$ is the product of  $x$'s for $q_1$ and $q_2$.

In e.g. \cite{Nik} we learn that any quadratic form is a direct product of indecomposable ones. These are:

\begin{itemize}\item[\textbf{type} $p^k_s$:] $G=\bbZ_{p^k}$ for $p$ an odd prime and $k\in\bbZ_{\ge 1}$, and  $s=\pm 1$. Then  $q(\ell)=\xi_{p^k}^{m\ell^2}$ for $\ell\in G$, where $2m$ is/is not a quadratic residue mod $p$, for $s=+1$ resp. $-1$. Here, $x^{3}=s^k\epsilon_{p^k}$ where $\epsilon_n=1$ resp. $-i$ for $n\equiv 1$ resp. $-1$ (mod 4). 

\item[\textbf{type}  $2^k_m$:]  $G=\bbZ_{2^k}$, and $m=\pm 1,\pm 3$, with $q(\ell)=\xi_{2^{k+1}}^{m\ell^2}$. Here,
$x^{3}=\epsilon_m\xi_8^{-m}$ where $\epsilon_m=-1$ if both $k$ odd and $m=\pm 3$, otherwise $\epsilon_m=+1$. 

\item[\textbf{type}  $2^k2^k_{\ i}$:] $G=\bbZ_{2^k}\times \bbZ_{2^k}$ and $q(\ell,m)=\xi_{2^k}^{\ell m}$.
Here, $x^{3}=1$.

\item[\textbf{type}  $2^k2^k_{\ ii}$:] $G=\bbZ_{2^k}\times \bbZ_{2^k}$  and $q(\ell,m)=\xi_{2^k}^{\ell^2+\ell m+m^2}$. Here, $x^3=(-1)^{k}$. 
\end{itemize}

\noindent\textit{Proof of Lemma 4.} 
Write $G\cong \bbZ_{m_1}\times\cdots\times\bbZ_{m_n}$ and let $g_i$ generate the $\bbZ_{m_i}$ factor. When $m_i$ is odd, put $q(g_i)=\langle g_i,g_i\rangle^{(m_i-1)/2}$; when $m_i$ is even, let $q(g_i)=\overline{\sqrt{\langle g_i,g_i\rangle}}$, for either choice of square-root. Now define
$$q(\sum_ik_ig_i)=\prod_iq(g_i)^{k_i^2}\prod_{i<j}\overline{\langle g_i,g_j\rangle}^{k_ik_j}$$
Then $q$ is well-defined, and is a quadratic form on $G$ realising the given $\langle,\rangle$.
 
Suppose $q'$ is another quadratic form realising the same symmetric pairing. Define $\psi(g)=q'(g)\overline{q(g)}$. Then \eqref{quadbil} implies $\psi\in \widehat{G}$. Moreover, $\psi(g)^2=1$  for all $g\in G$. Conversely, given any $\psi\in\widehat{G}$ satisfying $\psi^2=1$, $q'(g):=\psi(g)q(g)$ will be a quadratic form realising $\langle,\rangle$. Thus there are precisely $2^l$ different $q'$ realising $\langle,\rangle$, where $l$ is the number of even $m_i$. These correspond precisely to the different square-roots in the previous paragraph. 

Now turn to part (b). We know that  $q,q'$ realise the same symmetric pairing iff $q'=\psi q$ for some $\psi\in\widehat{G}$ obeying $\psi^2=1$. It thus suffices to consider 2-groups $G$. 

Suppose first that $q$ is \textbf{type} $2^k2^k_{\ i}$. Write $\psi_{a,b}(\ell,m)=(-1)^{\ell a+mb}$. When $k\ge 2$, the quadratic forms $q$ and $\psi_{a,b}q$ are equivalent, through the automorphism $\alpha(\ell,m)=(\ell+bm2^{k-1},m+a\ell 2^{k-1})$. \textbf{Type} $2^k2^k_{\ ii}$ when $k\ge 2$ is treated identically. When $q$ is \textbf{type} $2^12^1_{\ i}$ or  $2^12^1_{\ ii}$, $\psi_{a,b}q$ will also be either  \textbf{type} $2^12^1_{\ i}$ or  $2^12^1_{\ ii}$.

Finally, consider \textbf{type $2^k_m$} and $\psi$ is the unique nontrivial choice. Then $\psi q$ will be of \textbf{type} $2^k_{(2^k+1)m}$. But $2^k+1$ is a perfect square mod $2^{k+1}$ (namely it equals $(2^{k-1}+1)^2$) when $k\ge 3$. Trivially, \textbf{type} $2^1_m$ coincides with \textbf{type} $2^1_{m+4}$. However, \textbf{type} $2^2_m$ is inequivalent to \textbf{type} $2^2_{m+4}$ (to see this, note that the $T$ matrices in the associated modular data \eqref{WeilT} will be different).

It thus suffices to consider products of \textbf{type} $2^2_m$ quadratic forms for various $m$'s. Note that \textbf{type} $2^2_m\times 2^2_{m'}$ is equivalent to \textbf{type} $2^2_{m+4}\times 2^2_{m'+4}$, through automorphism $\alpha(\ell,\ell')=(\ell+2\ell',\ell'+2\ell)$. Using this, any product of \textbf{type} $2^2_m$'s is equivalent to a product of \textbf{type} $2^2_{\pm 1}$'s  together with at most one of  \textbf{type} $2^2_{\pm 3}$.
\textit{QED to Lemma 4}\medskip

\subsection{Lattices and reconstruction for pointed categories}\medskip

We follow the notation and terminology of \cite{CS}. An even positive-definite lattice $L$ is a free $\bbZ$-module of finite rank together with a positive-definite symmetric bilinear form $u\cdot v$ on $L$, such that $u\cdot u\in 2\bbZ$ for all $u\in L$. By the {dual} $L^*$ we mean $\{v\in\bbR\otimes_\bbZ L:v\cdot L\subseteq\bbZ\}$. Then for $L$ even, $G_L:=L^*/L$ is a finite abelian group, and for any classes $[u],[v]\in L^*/L$ the norm $u\cdot u$ is well-defined mod 2, and  $u\cdot v$ is well-defined mod 1. Hence $q_L([u]):=e^{\pi i u\cdot u}$ is a (well-defined) non-degenerate quadratic form on $G_L$. 

Given an even positive-definite lattice $L$, we get a rational vertex operator algebra $\cV(L)$ (see e.g. sections 6.4,6.5 of \cite{LL}) and rational local conformal net $\cA(L)$ \cite{DX}, and the category of modules for both is braided tensor equivalent to $\cC(q_L,G_L)$. For later convenience, we'll sketch some of this. As a vector space, $\cV(L)$ is spanned by combinations of the form 
$$h_{-n_1}^{(1)}h_{-n_2}^{(2)}\cdots h_{-n_m}^{(m)}\mathbf{e}^v$$
where $h^{(i)}\in \bbC\otimes_\bbZ L\cong \bbC^{\mathrm{dim}\ L}$ (corresponding to the Heisenberg part of the VOA), $n_i\in\bbZ_{>0}$,  and $\mathbf{e}^v$ (for $v\in L$) is the standard basis of the group ring $\bbC[L]$. We may demand $n_1\ge n_2\ge\cdots\ge n_m$ and that each $h^{(i)}$ is taken from a basis of $L$. The conformal weight of that combination is $\sum_i n_i+v\cdot v/2$.

The main purpose of this subsection is to show that  every category $\cC(q,G)$ can be realised through \textit{positive-definite} even lattices in this way. Though often stated, a complete explicit proof seems to be lacking in the literature (positive-definiteness is crucial for the complete rationality of the VOA and conformal net).

Given any classes $[u_i]\in G_L$ satisfying $u_i\cdot u_i\in 2\bbZ$ and $u_i\cdot u_j\in\bbZ$, write $G=\langle [u_1],[u_2],\ldots\rangle$. Then
$L[G]:=\cup_{[v]\in G}[v]$ is an even positive-definite lattice satisfying $L[G]/L=G$ and $|L[G]^*/L[G]|=|L^*/L|/|G|^2$. This is called the gluing construction in \cite{CS}.

The following result should be well-known --- see e.g.\ \cite{GaLa}.

\medskip\noindent\textbf{Lemma 5.} \textit{Let $L$ be an even positive-definite lattice. Then there exists an even \textit{self-dual} positive-definite lattice $\Lambda$, such that $L$ is isomorphic to a sublattice of $\Lambda$. Then $L^\perp:=\{v\in\Lambda\,|\,v\cdot L=0\}$ is also an even positive-definite lattice, with a group isomorphism $\varphi:G_L\rightarrow G_{L^{\perp}}$ satisfying $q_{L^\perp}(\varphi([u]))=\overline{q_L([u])}$ for all $[u]\in G_L$.}\medskip

Here, `isomorphic lattice' means the isomorphism preserves the inner product. In fact $\Lambda$ can be taken so that the orthogonal direct sum $L\oplus L\oplus L\oplus L\oplus L\oplus L\oplus L\oplus L$ is finite index in it, though which particular $\Lambda$ is chosen is not important to us. {This Lemma will be used shortly in the proof of Theorem 2.}

The main purpose of this subsection is reconstruction for pointed modular tensor categories (see also Corollary 1.10.2 in \cite{Nik}):

\medskip\noindent\textbf{Theorem 2.} \textit{Let $\cC$ be a pointed modular tensor category. Then  there is an even positive-definite lattice $L$ such that the lattice VOA $\cV(L)$ and  conformal net $\cA(L)$ both have  their category of modules Mod$(\cV(L))$ and Rep$(\cA(L))$ braided tensor equivalent to $\cC$.}\medskip

Conversely, it is certainly not the case that only lattice VOAs have pointed Mod$(\cV(L))$. For example, the cyclic orbifold $\cV=(V^{\natural}\otimes V^\natural)^{\bbZ_2}$ of the Moonshine module  has category of modules Mod$(\cV)$ tensor equivalent to the Drinfeld double $\cD(\bbZ_2)$ \cite{EG7}, which is pointed, but its associated Lie algebra $\cV_1$ is trivial so it  cannot be a lattice VOA (a lattice VOA $\cV(L)$ has associated Lie algebra $\cV(L)_1\supseteq\bbC^{\mathrm{dim}\,L}$). 

\smallskip The remainder of this subsection is devoted to a proof of Theorem 2. It suffices to show that any non-degenerate quadratic form $q$ on  an abelian group $G$ can be realised as $q_L$ on $G_L$ for some even positive-definite lattice $L$. We can restrict to the 4 types of quadratic forms listed last subsection, as $(q,G)\cong(q_L,G_L)$ and $(q',G')\cong(q_{L'},G_{L'})$ implies $(q\times q',G\times G')\cong (q_{L\oplus L'},G_{L\oplus L'})$ where $L\oplus L'$ denotes orthogonal direct sum.

We begin by realising  quadratic forms of \textbf{type} $p^k_s$ ($p$ an odd prime). Consider first primes $p\equiv -1$ (mod 4) (so $-1$ is not a quadratic residue of $p$). Note that the root lattice  $L=A_{p^k-1}$
has $G_L\cong\bbZ_{p^k}$, with generator $[1]$ satisfying $[1]\cdot[1]\equiv \frac{p^k-1}{p^k}$ (mod 2),
so this realises $m=(p^k-1)/2$, or equivalently $m=-2$. Letting $L^\perp$ be as in Lemma 5, we can also 
realise $m=+2$. Hence we have realised both $s=\pm $ (again, $-1$ is not a quadratic residue of $p$), and we're done \textbf{type} $p^k_s$ for those primes.

Next consider $p\equiv 1$ (mod 4). Choose any prime $p'\equiv -1$ (mod 4) with Legendre symbol $\left(\frac{p'}{p}\right)=s$   (there are infinitely many such $p'$).  Then by quadratic reciprocity, $\left(\frac{p}{p'}\right)=s$. By the previous paragraph (because $p'\equiv-1$ (mod 4)), there is an even positive-definite lattice $L'$ with $G_{L'}\cong\bbZ_{p'}$ and with generator $[\gamma]$ satisfying $\gamma\cdot\gamma\equiv -2p^k/p'$ (mod 2). Take $L$ to be the lattice gluing $(\sqrt{2p'p^k}\bbZ\oplus L'\oplus \sqrt{2}\bbZ)\langle[p'p^k,0,1],[2p^k,\gamma,0]\rangle$, where $[p'p^k,0,1]$ and $[2p^k,\gamma,0]$ denote the cosets containing the dual lattice vectors $(\frac{p'p^k}{\sqrt{2p'p^k}},0,\frac{1}{\sqrt{2}})$ resp. $(\frac{2p^k}{\sqrt{2p'p^k}},\gamma,0)$. Then $L$ is even  (since $[p'p^k,0,1]\cdot[p'p^k,0,1]\equiv \frac{p'p^k}{2}+\frac{1}{2}\equiv 0$ (mod 2), $[2p^k,\gamma,0]\cdot[2p^k,\gamma,0]\equiv \frac{2p^k}{p'}-\frac{2p^k}{p'}\equiv 0$ (mod 2), and $[p'p^k,0,1]\cdot[2p^k,\gamma,0]\equiv 0$ (mod 1)) and positive-definite (since $\sqrt{2p'p^k}\bbZ\oplus L'\oplus \sqrt{2}\bbZ$ manifestly is). Moreover, $G_L$ has order $(2p'p^k)(p')(2)/(2p')^2=p^k$ and contains $[2p',0,0]$ (since $[2p',0,0]\cdot[p'p^k,0,1]\equiv 0\equiv[2p',0,0]\cdot[2p^k,\gamma,0]$ (mod 1)), so $G_L\cong\bbZ_{p^k}$ with generator $[2p',0,0]$ having $[2p',0,0]\cdot[2p',0,0]\equiv \frac{2p'}{p^k}$ (mod 2). Thus $L$  realises \textbf{type} $p^k_s$.

\textbf{Type} $2^k_m$ is now easy. First, $L=\sqrt{2^k}\bbZ$ works for $m=+1$.  To do $m=3$, let $L'$ be a \textbf{type} $3^1_s$ lattice for $G=\bbZ_3$ with generator $[\gamma]$ satisfying $[\gamma]\cdot[\gamma]\equiv \frac{2^{k+1}}{3}$ (mod 2) (explicitly, we can take $L'$ to be the root lattice $A_2$ resp. $E_6$ if $k$ is even resp. odd). Then
$L=(\sqrt{3\cdot 2^k}\bbZ\oplus L')\langle[2^k,\gamma]\rangle$ is even positive-definite, and $L^*/L\cong\bbZ_{2^k}$ has generator $[3,0]$
having $[3,0]\cdot[3,0]\equiv \frac{3}{2^{k}}$ (mod 2), i.e. realising $m=3$.
The $L^\perp$ argument from Lemma 5 then takes care of $m=-1$ and $m=-3$.

Now look at \textbf{type} $2^k2^k_{\ i}$. Choose  any $L'$ from case (ii) with $G_{L'}\cong \bbZ_{2^k}$ and $m=-1$, and let $[\gamma]\in G_{L'}$ be a generator with $\gamma\cdot\gamma\equiv \frac{-1}{ 2^{k}}$ (mod 2). Let $L=( 2^{k}\bbZ\oplus 2^{k}\bbZ\oplus L'\oplus L')[1,1,\gamma,\gamma]$. Then $L$ is positive-definite and even, and $|G_L|=2^{4k}/(2^k)^2=2^{2k}$, with generators $a=[1,0,0,\gamma]=:(1,0)$ and $b=[0,-1,0,0,\gamma]=:(0,1)$, so $G_L\cong\bbZ_{2^k}\times\bbZ_{2^k}$, and we recover the desired quadratic form.

\textbf{Type} $2^k2^k_{\ ii}$ is similar. Choose any $L'$ from \textbf{type} $2^k_{-3}$, and let $[\gamma']\in G_{L'}$ be a generator with $\gamma'\cdot\gamma'\equiv \frac{-3}{2^k}$ (mod 2).
Let $L=( 2^{k}\bbZ\oplus 2^{k}\bbZ\oplus 2^{k}\bbZ\oplus L')[1,1,1,\gamma']$. The desired generators are $[1,1,0,0]=:(1,0)$ and $[1,0,1,0]=:(0,1)$, and we have recovered \textbf{type} $2^k2^k_{ii}$.

Incidentally, since $\cV(L)$ has central charge $d$ equal to the dimension of $L$, in all cases we obtain that $x^3=e^{-\pi i d/4}$
(this recovers Milgram's signature theorem).

\subsection{Chiral data associated to tori}

Let $T$ be a $d$-dimensional torus.  Section 2.2 of \cite{EG1} describe
the possible twists $H^3_T(T;\bbZ)$ and $H^1_T(T;\bbZ_2)$, where $T$ acts trivially on itself, while section 5.1 of \cite{FHLT}
identifies the transgressed ones. In particular, transgressed (non-degenerate) twists $\tau_L$ correspond 
bijectively to \textit{even positive-definite $d$-dimensional lattices} $L$. 

Fix a positive-definite inner product $ u\cdot v\in\bbR$ on $\bbR^d$. We interpret these $L$ as $d$-dimensional $\bbZ$-submodules
of $\bbR^d$, on which $ L\cdot L\subseteq\bbZ$ and $ u\cdot u\in 2\bbZ$ $\forall u\in L$. As mentioned last subsection,  the dual $L^*$ contains $L$. If we realise the torus $T$ by $T_L=\bbR^d/L$, then we get the natural
identifications $L=\mathrm{Hom}(\bbT,T)$ and $L^*=\mathrm{Hom}(T,\bbT)$.

Recall that $^{\tau_L} K_T^d(T)$  is naturally isomorphic to the group ring of $L^*/L$, which is
isomorphic as a ring to $K_{L^*/L}^0(\mathrm{pt})$ (on the other hand, $^{\tau_L} K_T^{d+1}(T)=K_{L^*/L}^1(\mathrm{pt})=0$). This can be understood through the calculation
$$^\kappa K_T(T)=R_T\otimes_{R_{T^*}}K_{T^*}(T)\,,$$
where $T^*$ is regarded as the quotient $T/L^*$, so that $R_{T^*}$ is naturally a subring of
$R_T$,
together with the identification $K_{T^*}(T)\cong \bbZ$ (on which $R_{T^*}$ acts by dimension)
coming from \cite{FHTii}, and the natural identification of $R_T\otimes_{R_{T^*}}\bbZ$
with $R_{T}/R_{T^*}=R_{L^*/L}$.

As mentioned last subsection, the modular data corresponding to lattice $L$ is given by \eqref{WeilT},\eqref{WeilS}, where $G=L^*/L$, $q([v])=e^{\pi i v\cdot v}$ and $\langle [u],[v]\rangle=\cQ_{[u]}([v])=e^{-2\pi i u\cdot v}$. The tensor unit
is $[0]$. In terms of VOAs and conformal nets, the torus theory $T_L$ is captured by the lattice VOA $\cV(L)$ and the lattice conformal net.

The question of which lattices have quaternionic simples, is  a little subtle. For example, the root lattice $L=D_8$ has $L^*/L\cong\bbZ_2\times\bbZ_2$, but none of these 4 simples is quaternionic. On the other hand, the root lattice $L=A_1\oplus E_7$ also has  $L^*/L\cong\bbZ_2\times\bbZ_2$, but now 2 of its 4 simples are quaternionic.

There are two frameworks here, related by the short exact sequence
$$1\rightarrow L^*/L\rightarrow T_L\rightarrow T_{L^*}\rightarrow 1\,.$$
We will work first with the finite group framework, then switch to the toroidal one. These two frameworks are elegantly related by \cite{FHLT}
into $C=(\bbR^d\times L^*)/L$, where $L$ embeds diagonally in $\bbR^d$ and $L^*$ (their $\Pi$ resp.\ $\Lambda$
are our $L$ resp.\ $L^*$).
We would expect that much of our work in this section can be rephrased in that language.

Let $G$ be a finite abelian group, and $H$ any subgroup. Then $\widehat{H}$ is naturally a quotient of $\widehat{G}$: the projection $\widehat{G}\to\widehat{H}$ is just restriction to $H$. Let $K\le\widehat{G}$ be the kernel of that projection. Now suppose $G$ has  a non-degenerate pairing $\gamma:G\to \widehat{G}$ (not necessarily symmetric). Write $H^\perp=\gamma^{-1}(K)$:  i.e.\begin{equation}H^\perp=\{k\in G\,|\, \gamma(k)(h)=1\ \forall h\in H\}\label{Hperp}\end{equation}
Then $H\cong G/H^\perp$, so in particular $|H^\perp|=|G|/|H|$. Of course we also have $\widehat{K}\cong \widehat{\hat{G}}/\widehat{\hat{H}}$, i.e. $H^\perp\cong G/H$.  When $\gamma$ is in  fact symmetric, then $(H^\perp)^\perp=H$.

Another little fact which can be useful: given lattices $M\supseteq L$ of the same dimension, there is a natural group isomorphism $M/L\cong \widehat{L^*/M^*}$, by $m\mapsto e^{2\pi\i \lambda\cdot m}$.

\subsection{A Galois correspondence and Jones tower for lattices}

 Let $L\subseteq M$ be arbitrary
lattices of equal dimension, and write $G=M/L$.
There is an elementary  correspondence between intermediate lattices $L\subseteq D\subseteq M$,
and subgroups $H\le G$, through $H=D/L$ and $D=\cup_{[v]\in H}[v]$. Write $M_H$ for the sublattice corresponding to subgroup $H$.  In the language of section 4.2, $D$ is the gluing $L[H]$. 

When there is a non-degenerate pairing $\gamma$ on $G$, this becomes a
Galois (contravariant) correspondence using \eqref{Hperp}. More precisely, given a subgroup $H\le G$, define $M^H$ to be the lattice $M_{H^\perp}$.  Then $L\subseteq M^H\subseteq M$, $M^H/L\cong H^\perp$ and $M/M^H\cong
(M/L)/(M^H/L)=G/H^\perp\cong H$. Moreover, when $H\le K\le G$, we have $M^H\supseteq M^K$ and $M^H/M^K\cong H^\perp/K^\perp\cong K/H$.

This has a direct  Galois interpretation. For any $x\in\bbC\otimes L$, define a map on the group algebra  $\bbC[M]$ by $x.\mathbf{e}^v=e^{2\pi i x\cdot v}\mathbf{e}^v$. Note that $x\in M^*$ acts trivially. Through the pairing, we  identify $G=M/L$ with $\widehat{G}=L^*/M^*$, $g\in G\mapsto [x_g]\in L^*/M^*$. Then $G$ acts on the group algebra through $x_g$. Now choose any subgroup $H\le G$. Define $K\le\widehat{G}$ as at the end of section 4.3, and interpret $K\le L^*/M^*$ as usual. Then the fixed-point subalgebra $\bbC[M]^K$ is naturally identified with $\bbC[M_H]$. Choosing a different $K$ would recover $\bbC[M^H]$ spot-on.

Indeed, this action of $G$ on lattice group algebras plays a fundamental role in VOA theory.
Let $L$ be any even positive-definite lattice of dimension $d$. Write $G_L:=L^*/L$ and for each $[v]\in G_L$ put $q_L([v])=e^{\pi i v\cdot v}$.  Recall the sketch of lattice VOAs offered in section 4.2.
The automorphism group of the VOA $\cV(L)$ (c.f. \cite{DoNa}) is generated by isometries of the lattice, as well as automorphisms of the form $\alpha_x$ for $x\in\bbC\otimes L$. An isometry $\sigma$ of $L$ sends  $\mathbf{e}^v$ to $\mathbf{e}^{\sigma(v)}$ and each $h^{(i)}$ to $\sigma(h^{(i)})$, whereas $\alpha_x$ fixes the $h^{(i)}$ and sends $\mathbf{e}^v$ to $x.\mathbf{e}^v$ as above. The corresponding orbifolds behave very differently --- in particular, orbifolding by isometries $\sigma$ generally gives a nonlattice theory (see section 5.4), whereas orbifolding by $x\in\bbQ\otimes_\bbZ L$ always yields a lattice VOA.
 We need both kinds of orbifolds in section 5, but in this subsection we consider only the latter.

\medskip\noindent\textbf{Lemma 6.} \textit{Let $L$ be an even positive-definite lattice, and let $L_0$ be a sublattice of $L$ of finite index. Then the lattice VOA  $\cV({L_0})$ is a VOA orbifold of $\cV(L)$, by a group $G\cong L/L_0$.}\medskip

\noindent\textit{Proof.} The automorphism $\alpha_x$ is trivial iff $x\in L^*$ so we are really only interested in $x$ mod $L^*$.
When $x\in\bbQ\otimes_\bbZ L$, $\alpha_x$ is finite order. 
Choose any finite set of $x_i'\in\bbQ\otimes_\bbZ L$; then the group $G'$ generated by the automorphisms $\alpha_{x'_i}$ is finite, and the orbifold 
$\cV(L)^{G'}$ is a lattice theory $\cV({L'})$ where $L'=\{u\in L:x_i\cdot u\in\bbZ\ \forall i\}$.

In particular, choosing $x_i\in L^*_0$ so they generate $L^*_0/L^*=:G$, we obtain $\cV(L)^G=\cV({L_0})$. As mentioned in section 4.3,  $L^*_0/L^*\cong L/L_0$.  \textit{QED to Lemma 6}\medskip

Now let $n$ be the exponent of $G$, so $nG=0$. We have $nL\subseteq n M_H\subseteq nM\subseteq L$ and both $n M_H/nL
\cong H$ and $nM/nL\cong G$, so we get  intermediate  lattices $L^H$ and $L^G$ satisfying $nL\subseteq L^G\subseteq L^H\subset L$, $L/L^G\cong G$, 
$L/L^H\cong H$, and $L^H/L^G\cong H^\perp$.

We can also go in the other direction, though it is less natural. Find some $m\in\bbZ_{>0}$ for which $m^{-1}L\supset M$. Then $m^{-1}M/m^{-1}L\cong G$ so there is a quotient (hence subgroup $G'$) of the abelian group $m^{-1}M/M$ which is isomorphic to $G$. Let $L'$ be the lattice gluing $M[G']$, so $L'/M\cong G$. Then for each subgroup $H\le G$, we get as before the lattices $L^{\prime H},L'_H$ intermediate to $M\subset L'$, with the desired quotients. We use this construction in section 4.6.

Thus, given $L\subseteq M$, a non-degenerate pairing on $G=M/L$, and a subgroup $H\le G$, we get  a  tower $$\cdots\subset  M_{-2}\buildrel H\over\subset M_{-1}\buildrel H^\perp\over \subset
M_0\buildrel H\over\subset M_1\subset \cdots$$ such that $M_0=L,M_1=M^H$, $M_2=M$, $M_{2i}/M_{2i-1}\cong H^\perp$ and $M_{2i+1}/M_{2i}\cong H$. In the special case $H=G$, this is analogous to the Jones tower   of subfactors corresponding to the inclusions $A^G\subset A\subset A\sdprod G$. 

This tower of lattices lifts to a  tower of lattice VOAs, in the negative direction (in the positive direction, it terminates when $M_i$ is no longer even). Note that not all  VOAs can belong to such a tower --- e.g. any tower containing a Virasoro minimal model would necessarily be finite.

\subsection{The finite group framework}

Let $L$ be any even positive-definite lattice. The modular data associated to $L$ is described in section 4.3.
Note that for this modular data, any modular invariant automatically satisfies
\eqref{scinvdef}. Moreover, $S_{[a],[b]}$ never vanishes and the simple-current stabiliser of any 
$[a]\in\Phi$ is trivial. This means it is sufficiently nonzero, in the sense of section 3.2, and 
Theorem 1 exhausts all modular invariants.

\medskip\noindent\textbf{Theorem 3(a)}  \textit{Let $L$ be even and write $G=L^*/L$. Let $q$ and $\langle,\rangle$ denote the associated quadratic form and symmetric pairing on $G$. The following are in bijection:}

\begin{itemize}

\item[(i)] \textit{the modular invariants $M$ for $L$;}

\item[(ii)] \textit{all pairs $(J,[\phi])$ where $J\le G$ containing no quaternionic elements, and $[\phi]\in H^2_J(\mathrm{pt};\bbT)$;}

\item[(iii)] \textit{subgroups  $D_\pm$ of $G$ with $q|_{D_\pm}=1$,  together with an isomorphism $\sigma:D_+^\perp/D_+\to D_-^\perp/D_-$ satisfying $q(k)=q(\sigma(k))$ for all $k\in D_+^\perp$;}

\item[(iv)] \textit{subgroups  $Z$ of $G\times G$ with $Z=Z^\perp$ and $q_2|_Z=1$, where on $G\times G$ we use the non-degenerate pairing $\langle(g,h),(g',h')\rangle_2= \langle g,g'\rangle \overline{\langle g',h'\rangle}$ and associated quadratic form $q_2(g,h)=q(g)\overline{q(h)}$.}

\end{itemize}

 \medskip\noindent\textbf{(b)} \textit{There is a one-to-one correspondence between all subgroups
 $Z\le G\times G$ satisfying $Z^\perp=Z$, and all triples $(J,[\psi],\phi)$ where $J\le G$,
no elements of $J$ are quaternionic,  $\phi\in \widehat{J/2J}$ and
 $[\psi]\in H^2_J(\mathrm{pt};\bbT)$.}\medskip
 
The notation $D^\perp$  is defined in \eqref{Hperp} where $\gamma(g)(h)=\langle g,h\rangle$. In (iii), as explained in the proof, we have that $D_\pm^\perp\ge D_\pm$, and both the pairing and quadratic form of $G$ restrict to well-defined functions on $D_\pm^\perp/D_\pm$, so we use the same symbols. The parametrisation of (ii) is given in Theorem 1, and directly describes the nimrep as we'll see.  Part (iii) gives the type 1 parents and permutation $\sigma$ appearing in \eqref{modinvsigma}. Part (iv) gives a geometric interpretation of these modular invariants: they are the self-dual lattices lying between the indefinite lattice $L\oplus\sqrt{-1}L$ and its dual. As we will see in the proof, the relation between the subgroup $Z$ of (iv) and the matrix $\cZ$ of (i) is
\begin{equation}g,h\in Z\ \mathrm{iff}\ \cZ_{g,h}=1\,,\ \forall g,h\in G\label{ZZ}\end{equation}
(the other entries $\cZ_{g,h}$ all equal 0).
There is a similar 4-part description for (b), but we only list two.

\medskip\noindent\textit{Proof.} The equivalence of (i) and (ii) is given by Theorem 1. The equivalence of (i) and (iii) is given in Proposition 1(a) of \cite{EG3}. We will complete the proof of Theorem 3(a) by establishing the equivalence of (iii) and (iv).

First, consider any subgroup $D\subseteq G$ with $q(d)=1$ for all $d\in D$. Then \eqref{quadbil} implies that $D^\perp\supseteq D$, and also that the quadratic form $q$ (hence the pairing $\langle,\rangle$) is  well-defined on $D^\perp/D$. Hence any such $D$ with quadratic form $q$ defines a possible extension (type 1 module category) for $(G,q)$.  
 As discussed in  section 2.3, sigma-restriction is the branching rules from the modular tensor category of $(D,q)$ (with simples in $D^\perp/D$), to that of $(G,q)$. Sigma-restriction here is given by the correspondence
 $$\begin{tikzpicture}
  \matrix (m) [matrix of math nodes,row sep=2em,column sep=2em,minimum width=2em]
  {
    & \bbZ[D^\perp] &  \\ \bbZ[D^\perp/D]&&\bbZ[G]\\ };
  \path[-stealth]
    (m-1-2) edge node [left] {$\pi\ $} (m-2-1)
            edge node [right] {\ $\iota$} (m-2-3);
\end{tikzpicture}$$
where we identify the $K$-group $K(G)$ etc with the group ring $\bbZ[G]$, where $\pi$ denotes the obvious projection $D^\perp\to D^\perp/D$, and $\iota$ is inclusion $D^\perp\to G$.

Let's begin by showing that (iii) implies (iv). Let $D_\pm$ and $\sigma$ be as in (iii). Then the matrix product \eqref{modinvsigma} is captured by the composition
$$\begin{tikzpicture}
  \matrix (m) [matrix of math nodes,row sep=2em,column sep=2em,minimum width=2em]
  {
 &\bbZ[{D_+}^\perp]&  & & \bbZ[{D_-}^\perp] &  \\ \bbZ[G]&&  \bbZ[D_+^\perp/D_+]&  \bbZ[ {D_-}^\perp/D_-] & &\bbZ[G] \\ };
  \path[-stealth] (m-1-2) edge node [left] {$\iota_+$\ \ } (m-2-1) (m-1-2) edge node [right] {\ $\pi_+$} (m-2-3) (m-2-3) edge node [above] {$\sigma$} (m-2-4)
    (m-1-5) edge node [left] {$\pi_-$\ } (m-2-4) (m-1-5) edge node [right] {\ $\iota_-$} (m-2-6)   ;
\end{tikzpicture}$$
of correspondences.
This composition collapses (pulls back) to
$$\begin{tikzpicture}
  \matrix (m) [matrix of math nodes,row sep=2em,column sep=2em,minimum width=2em]
  {
    & \bbZ[Z] &  \\ \bbZ[G]&&\bbZ[G]\\ };
  \path[-stealth]
    (m-1-2) edge node [left] {$p_+\ $} (m-2-1)
            edge node [right] {$\ p_-$} (m-2-3);
\end{tikzpicture}$$
where $Z$ is the group $$0\rightarrow D_-\rightarrow Z\rightarrow D_+^\perp\rightarrow 0$$
consisting of all $(g_+,g_-)\in D_+^\perp\times D_-^\perp$ satisfying $\pi_-(g_-)=\sigma(\pi_+(g_+))$,
and where $p_\pm(g_+,g_-)=g_\pm$. This pull-back defines the modular invariant matrix $\cZ$ (as well as the group $Z$) and we see \eqref{ZZ} holds.
 
 Conversely, let $Z\le G\times G$ be as in (iv). Write $D_+=\{g\in G\,|\,(g,0)\in Z\}$ and $D_-=\{g\in G\,|\,(0,g)\in Z\}$. Then $q(D_\pm)=1$ because $q_2(Z)=1$. Hence $(g_+,g_-)\in Z$ implies $g_\pm\in D_\pm^\perp$ because $Z^\perp=Z$. Also, both $(g,h),(g,h')\in Z$ only when $h-h'\in D_-$. Thus $|Z|$ is at most $|D_\pm^\perp|\,|D_\mp|=\frac{|G|\,|D_\mp|}{|D_\pm|}$. But $Z^\perp=Z$ forces $|Z|=|G|$. Hence for every $g\in D_+^\perp$ there is an $h\in D_-^\perp$ such that $(g+d_+,h+d_-)\in Z$ for all $d_\pm\in D_\pm$, and this assignment $\sigma(g+D_+)=h+D_-$ defines an isomorphism $D_+^\perp/D_+\to D_-^\perp/D_-$. Of course $q_2(Z)=1$ requires that $q(g)=q(\sigma(g))$ for all $g\in D_+^\perp$.
 
The proof of part (b) is similar. \textit{QED to Theorem 3}\medskip

In terms of the parametrisation (iii), the full system is $\bbZ[G/D_-\times G/D_+^\perp]$.
Alpha-induction $\alpha_\pm:\bbZ[G]\rightarrow \bbZ[L^*/D_-\times L^*/D_+^\perp]$ 
are given by  $\alpha_+=(\sigma\circ\pi_+,\pi_{+}^\perp)$ and $\alpha_-=(\pi_-,0)$, where $\pi_{+}^\perp:G\rightarrow G/D_+^\perp$
is the obvious projection. In particular, these $\alpha_\pm$ are linear, and so respect tensor products, and through \eqref{modinvalpha} recovers the modular invariant $\cZ$ of \eqref{ZZ}. In terms of the parametrisation (iii), the Grothendieck group of the associated module category is $\bbZ[G/J]$, and the action of the fusion ring $\bbZ[G]$ on it (called the nimrep) is $g.[h]=[g+h]$. To see this, note that the exponents (i.e. the diagonal entries) of the modular invariant are the $j\in J^\perp$ (all with multiplicity 1), so the complete list of  eigenvalues of the nimrep action of $g$ are required to be $\langle g,j\rangle$ for all $j\in J^\perp$.

The relation between the $(J,[\psi])$ parametrisation of (ii),
and the $(D_\pm,\sigma)$ parametrisation of (iii), is as follows. Given any subgroup $J\le G$ and class $[\psi]\in H^2_J(\mathrm{pt};\bbT)$, define the homomorphism $\epsilon:J\rightarrow \widehat{J}$ as in Theorem 1. Let $\varphi:G/J^\perp\to \widehat{J} $ be the isomorphism $\varphi([g])(j)=\langle g,j\rangle$, and define $\varepsilon=\varphi^{-1}\circ\epsilon:
J\rightarrow G/J^\perp$. From the proof of Theorem 1(a), we find in this notation $D_-=J_R=\mathrm{ker}(\varepsilon)$ and $D_+=J_0=\{j\in J\,|\,\varepsilon(j)=-j+J^\perp\}$. By Lemma 2, 
the image of $\varepsilon$ is $D_+^\perp/J^\perp$, and thus $\varepsilon$ gives an isomorphism $J/D_-\rightarrow D^\perp_+/J^\perp$ is an isomorphism. Therefore, given any $a\in D_+^\perp$, there is a unique class $[j]\in J/D_-$ such that $a+J^\perp=\varepsilon(j)$. Define $\sigma(a)=a+j+D_-$. We know that $\sigma$ is constant on each $D_+$ class, and that it defines an isomorphism $D_+^\perp/D_+\cong D^\perp_-/D_-$.

For example, the identity modular invariant $\cZ_{[u],[v]}=\delta_{[v],[u]}$ corresponds to $J=0$. More generally, the type 1 modular invariant (pure extension type) with $D_+=D_-=D$ corresponds to
$J=D$ and $\varepsilon(j)=0$ (equivalently, $\epsilon(j)=1$). `Charge-conjugation' $\cZ_{[u],[v]}=\delta_{[v],[-u]}$ corresponds to $J=2G$ and $\varepsilon(j)=-j/2$.

\subsection{The toroidal framework}

As warm-up for section 6, we should also reinterpret the previous subsection using 
Dixmier-Douady bundles for the torus acting on the torus. 
Recall from \cite{EG1} the construction of this bundle.  As mentioned in section 4.3, a
transgressed (non-degenerate) 
twist for the trivial action of a torus on itself can be identified with the choice of an even $d$-dimensional lattice $L$.
It is convenient to identify the torus $T$ with $T_L=\bbR^d/L$. Let $\pi$ denote the regular representation of $T_L$ on the
Hilbert space $\cH=L^2(T_L)$. For any $\gamma\in L^*$ we
have a character (1-dimensional representation) $\chi_\gamma$ for $T_L$ defined by
$\chi_\gamma (t) ={\rm e}^{2\pi \i \,\gamma(t)}$; we can extend it linearly to a functional on $\bbR^d$. As in \cite{EG1,EG2}, 
choose any unitary $U_\gamma$ satisfying $U_\gamma \pi U_\gamma^*
=\chi_\gamma \pi$. The bundle  on $T_L$, with
fibres the compacts $\mathcal{K}=\mathcal{K}(\cH)$, is
defined using the gluing conditions $f(t)=U_{\ell}
f(t+\ell) U_{\ell}^*$, for all $\ell\in L$,
$t\in\bbR^d$. We denote these bundles symbolically by $T/\!/_LT$  or just $\cA_L$.

By the $K$-group  ${}^LK_0^T(T)$ we mean the $K$-theory of the C$^*$-algebra of $\cK$-valued sections of $\cA_L$. It is readily computed to be isomorphic to the ring  $\bbZ[L^*/L]$, which we recognise as the fusion ring for the lattice theory of section 4.3.

Consider first a type 1 modular invariant, corresponding to a pair
 $L\subset D$ of even lattices (i.e. they are both transgressed) and identity map $\sigma$.
Sigma-restriction here is clear: it corresponds to the bundle map $\cA_L\rightarrow 
\cA_D$ arising from the natural projection $T_L\rightarrow T_D$ sending $x+L$ to $x+D$.
The associated modular invariant $bb^t$, as an element in $KK(\cA_L,\cA_L)={}^{\tau_L}KK_{T}(T,T)$, can be pulled-back to give the diamond 
$$\begin{tikzpicture}
  \matrix (m) [matrix of math nodes,row sep=2em,column sep=2em,minimum width=2em]
  {
    & \cA_{D'} &  \\ \cA_L&&\cA_L\\ &\cA_D&\\ };
  \path[-stealth]
    (m-1-2) edge node [left] {$p'\ $} (m-2-1)
            edge node [right] {$\ p'$} (m-2-3)  (m-2-1) edge node [right] {$\ p$} (m-3-2)
           (m-2-3) edge node [left] {$\ p$} (m-3-2);
\end{tikzpicture}$$
where $D'\subset L$ satisfies $L/D'\cong D/L$ (using the tower of
section 4.4), and $p$ is the obvious projection $T_{D'}/L=T_L$.

The other extreme,  type 2, when $D_\pm=L$ and  $\sigma$ is nontrivial, requires that we construct a new bundle.
Representations of $T_L$ can be regarded as projective representations of the
quotient $T_{L^*}=T_L/L^*$, with twists given by $\psi\in L^*/L$. This allows us to decompose
the compacts $\cK=\cK(\cH(T_L))$ into the finite sum $\oplus_\psi \cK^\psi$. We get
a vector bundle over $T_{L}$, with components parametrised by $\psi\in L^*/L$,
with fibres $\cK^\psi$, equivariant with respect to $T_{L^*}$, with sections 
in the $\psi$-component satisfying $f_\psi(t+\lambda)=\psi(\lambda)\mathrm{Ad}\,U_\lambda\,
f_\psi(t+\lambda)$, $\forall \lambda\in L^*$, $\psi\in L^*/L$. Call this bundle $\cB_L$ or $T/\!/_{L^*}T$. We get a natural map
$\iota:\cA_L\rightarrow \cB_L$, which is basically the identity on the base with a projection $\pi:T_L\rightarrow T_{L*}$ of the groups
of morphisms. Then the modular invariant corresponds to
$$\begin{tikzpicture}
  \matrix (m) [matrix of math nodes,row sep=2em,column sep=2em,minimum width=2em]
  {
 \cA_L&&&&\cA_L\\    & \cB_{L} & & \cB_{L}&\\ };
  \path[-stealth]
    (m-1-1) edge node [left] {$\iota\ $} (m-2-2)
         (m-1-5)   edge node [right] {$\ \iota$} (m-2-4)
     (m-2-2)   edge node [above] {$\sigma$} (m-2-4)        ;
\end{tikzpicture}$$
where $\sigma$ permutes the components of the bundle.

The general modular invariant can be understood through the combination $b_-\sigma b_+^t$ as
the composition 
$$\begin{tikzpicture}
{\small  \matrix (m) [matrix of math nodes,row sep=2em,column sep=2em,minimum width=2em]
  {
    \cA_L&&&&&\cA_{L}\\ &\cA_{D_+}&&&\cA_{D_-}&\\  &&\cB_{D_+} &\cB_{D_-}&&\\ };
  \path[-stealth] (m-1-1) edge node [left] {$p_+\ \ $} (m-2-2) (m-1-6) edge node [left] {$p_-\ $} (m-2-5) 
    (m-2-2) edge node [left] {$\iota_+\ $} (m-3-3) (m-3-3) edge node [above] {$\sigma$} (m-3-4)
(m-2-5) edge node [left] {$\iota_-\ $} (m-3-4)   ;}
\end{tikzpicture}$$

It is elementary to translate the alpha-induction given in the previous subsection, into this language: e.g. the obvious projection
$\pi_\pm$ on morphisms (and identity on base space) sends the bundle corresponding to $T_{D_\pm}$ acting on $T_L$, to the bundle corresponding to $T_L$ acting
on $T_L$, and corresponds to the right-way map $\pi_\pm:D_\pm/L\rightarrow L^*/L$ (embedding) and wrong-way map
$\pi_{\pm!}: L^*/L\rightarrow D_\pm/L$ (projecting away cosets not in $D_\pm$) of $K$-groups.

\section{Tambara-Yamagami  categories}

The Tambara-Yamagami categories $\cT\cY_\pm(G,\lan,\ran)$ are parametrised by a finite abelian group $G$, a non-degenerate symmetric pairing $\lan\cdot,\cdot\ran$ on $G$, and a sign $\pm$.
 As shown in \cite{TY}, these are the only fusion categories whose primaries (equivalence classes of simples) are $[\alpha_g]$ ($g\in G$) and $[\rho]$, and which obey the fusions \eqref{TYfusions}.

These categories are pairwise distinct:  $\cT\cY_s(G,\lan,\ran)$ and $\cT\cY_{s'}(G',\lan,\ran')$ are tensor equivalent  iff $s=s'$ and $\lan g,h\ran=\lan \varphi(g),\varphi(h)\ran'$ $\forall g,h\in G$, where $\varphi:G\rightarrow G'$ is a group isomorphism. 
The Tambara-Yamagami categories are amongst the simplest examples of near-group categories. 

The sign $\pm$ is not at all mysterious. Any fusion category with a grading by some (not necessarily abelian) group $K$, can have its associativity isomorphisms twisted by a cocycle in $H^3(G;\bbT)$. This is the explanation for the $H^3$ twists \cite{DW} of the category of $G$-graded vector spaces Vect$_G$, and for the $\bbZ_n$ twists of the Kazhdan-Wenzl sl$(n)$ quantum group categories \cite{KazW}. Each Tambara-Yamagami category has a $\bbZ_2$ grading, where $\alpha_g$ has grading 0 and $\rho$ has grading 1. Thus it can be twisted by $H^2(\bbZ_2;\bbT)\cong\bbZ_2$. This is the source of the sign $\pm$.

\subsection{Tambara-Yamagami as a Potts model}

Izumi (Example 3.7 in \cite{iz}) realised each  $\cT\cY_\pm(G,\lan,\ran)$ as a system of endomorphisms on a Cuntz algebra (though in a somewhat ad hoc manner). In this subsection we recover this in a calculation-free way, using the Potts model. This subsection, apart from Theorem 4, is not used in the rest of the paper.

By the \textit{Cuntz algebra} $\cO_n$, we mean the universal C$^*$-algebra generated by $n$ generators $S_i$ (and their adjoints) satisfying {$S_i^*S_j=\delta_{i,j}$ and $\sum_jS_jS_j^*=1$}. Fix a finite abelian group $G$. We will label the generators of $\cO_{|G|}$ by $S_g$ for $g\in G$ (so $S_g^*S_h=\delta_{g,h}$ etc).

An excellent way to realise fusion categories is using endomorphisms on an algebra $A$. Objects are algebra endomorphisms; given endomorphisms $\rho$ and $\tau$, 
by Hom$(\rho, \tau )$ we mean the space of intertwiners 
 $\{ t \in A \,|\, t\rho(a) = \tau(a) t,\, \forall a \in A\}$. Tensor product becomes composition, and direct sum is realised using e.g. Cuntz generators. When the category is unitary and $A$ is e.g.\ a C$^*$-algebra (the situation considered here), then the endomorphisms should be $*$-endomorphisms, and the dual of an object is given by its adjoint. See e.g.\ \cite{iz3,EG5} for details.

 In particular, we say that a system of $*$-algebra endomorphisms $\alpha_g,\rho$ on $\cO_{|G|}$ realises the fusion category $\cT\cY_\pm(G,\lan,\ran)$, if the compositions $\alpha_g\circ\alpha_h$, $\alpha_g\circ \rho$, and $\rho\circ \alpha_g$, equal
 $\alpha_{g+h}$, $\rho$, and $\rho$ resp., up to conjugations by unitaries of $\cO_{|G|}$, and $\rho(\rho(x))=\sum_gS_g\alpha_g(x)S_g^*$.

\medskip\noindent\textbf{Theorem 4.} \cite{iz3} Each $\cT\cY_s(G,\lan,\ran)$ can be realised by the following system $\alpha_g.\rho$ of endomorphisms on $\cO_{|G|}$: for each $g\in G$, define $\alpha_g(S_h)=S_{g+h}$ and 
\begin{equation}\label{rho}\rho(S_h)=U(h)\frac{s}{\sqrt{n}}\sum_{k\in G}S_kU(h)^*\end{equation}
where $U(h)=\sum_k\lan h,k\ran S_kS_k^*$ is unitary. Then $$\alpha_g\alpha_h=\alpha_{g+h}\,, \ U(g)U(h)=U(g+h)\,,\ U(g)^*=U(-g)\,,\ \alpha_g(U(h))=\overline{\lan g,h\ran}U(h)\,,$$ $$\alpha_g\circ\rho=\rho\,,\ \rho\circ\alpha_g=Ad(U(g))\circ\rho\,,\ \rho(U(g))=\sum_h S_{h-g}S_h^*\,,\ \rho^2(S_h)=\sum_gS_gS_{g+h}S_g^*$$

The subtle aspect here of course is the endomorphism $\rho$.
In this subsection we see that this can be derived from the high temperature - low temperature duality in the Potts model just as the Ising fusion category  is related to Krammers-Wannier duality in the Ising model \cite{Ev2}. This can be regarded as a conceptual explanation for the  $\rho$ formula \eqref{rho}.

The Ising mode  can be
generalized to a $Q$-state standard Potts model, where $Q=|G|$. The transfer matrix formalism for the $Q$-state Potts 
Hamiltonian
\begin{equation}
\label{pottshamilt}
H(\sigma)=-{\scriptstyle\sum}_{i,j\ n\cdot
n}J\delta(\sigma_i,\sigma_j)
\end{equation}
on the two dimensional lattice $\Bbb{Z}^2$ leads to the
following algebraic set up.
For each $g\in G$, define $Q\times Q$ matrices $U_g=\mathrm{diag}_{k\in G}\lan g,k\ran$ and $V_g$ by $(V_g)_{h,k}=\delta_{k,g+h}$. Then both assignments $g\mapsto U_g,V_g$ define $Q$-dimensional unitary representations of $G$. Moreover, $V_hU_g=\lan g,h\ran U_gV_h$
and together the $U_g,V_g$ generate the $Q \times Q$ complex matrix algebra
$M_Q = \mathrm{End}({\mathbb C}^Q)$. For each $g\in G$, define a sequence of unitaries $W^g_i$ in $M_Q\otimes
M_Q\otimes M_Q\otimes\cdots$ by
\begin{equation}
\label{defV}
W^g_{2i+1}=1\otimes\cdots\otimes1\otimes U_{-g}\otimes
U_g\otimes1\otimes\cdots
\end{equation}
with $U_g$ appearing as the $(i+1)$th factor, and
\begin{equation}
\label{defW}
W^g_{2i}=1\otimes\cdots\otimes1\otimes V_g\otimes1\otimes\cdots
\end{equation}
with $V_g$ appearing as the $i$th factor. Then
\begin{equation}
\label{ccr}
(W^g_j)^Q=1\,,\quad W^g_iW^h_{i+1}=\overline{\lan g,h\ran} W^h_{i+1}W^g_i\,,\quad
W^g_iW^h_j=W^h_jW^g_i\,,\quad|i-j|>1\,.
\end{equation}
The symmetric group $S_Q$ acts on $\Bbb{C}^Q$ by permuting
basis vectors and so induces a product action on the UHF C$^*$-algebra 
$F_Q=\otimes_{\Bbb{N}} M_Q$. In particular, for each $g\in G$ there is a
$G$-action on $F_Q$ by $\bigotimes_i{\rm
Ad}(V_g)=\prod_i{\rm Ad}(W^g_{2i})$. The
Temperley-Lieb\index{Temperley--Lieb} operators in $F_Q$ are
the spectral projections of $W_i^g$ corresponding to eigenvalue
1, namely
\begin{equation}
\label{espectral}
e_i=\frac{1}{Q}\left(\sum_g W^g_i\right)\,.
\end{equation}
In particular
\begin{eqnarray}
e_{2i-1}&=&1\otimes\cdots \otimes 1\otimes \pi\otimes1\otimes\cdots
\label{eodd}\\
e_{2i}&=&1\otimes\cdots\otimes 1\otimes \pi'\otimes1\otimes\cdots
\label{eeven}
\end{eqnarray}
where $\pi$ (occupying the $i$th and $(i+1)$th spots) and $\pi'$ (occupying the $i$th spot) are projections in $M_Q\otimes M_Q$ and $M_Q$
respectively given by
\begin{equation}
\label{g&f}
\pi=\sum_{i=1}^QE_{ii}\otimes E_{ii}\,,\qquad
\pi'=\sum_{i,j=1}^QE_{ij}/Q
\end{equation}
if $\{E_{ij}: i,j=1,2,\dots,Q\}$ are matrix units for
$M_Q$. In the Ising case:
\begin{equation}
\label{ising_e}
e_{2i}=(1+\sigma_z^i)/2\,,\quad
e_{2i+1}=(1+\sigma_x^i\sigma_x^{i+1})/2\,.
\end{equation}
The transfer matrix for the two dimensional $Q$-state Potts 
model\index{Potts model} is then (cf. Chapter 12 of \cite{Ba2})
described as
\begin{equation}
\label{transfer}
V=\exp \left(K_2\sum e_{2i+1}\right)\,,\qquad
W=\exp \left(K_1^*\sum e_{2i}\right)
\end{equation}
with Kramers-Wannier duality being
the shift $\rho : e_i\to e_{i+1}$. Here the family $\{e_i\}$
satisfy the relations that non-nearest neighbours commute:
\begin{equation}
\label{TL1}
e_ie_j=e_je_i\qquad|i-j|>1
\end{equation}
whilst nearest neighbours satisfy
\begin{equation}
\label{TL2}
e_ie_{i\pm1}e_i=\tau e_i
\end{equation}
where $\tau^{-1}=Q$.
The family
$\{e_i\}$ is
$S_Q$ invariant, and the $\{W_i^g\}$ generate the $G$
fixed point algebra of
$F_Q$.
 This formulation is related to that of generalized
Clifford algebras as follows. The formulae
\begin{equation}
\label{J-W}
\Gamma^g_i=W^g_1W^g_2\cdots W^g_i\,,\qquad
W^g_i=(\Gamma_{i-1}^g)^{-1}\Gamma^g_i
\end{equation}
is the generalized Jordan-Wigner transformation required
between (\ref{ccr}) and
\begin{equation}
\label{clifford}
(\Gamma_i^g)^Q=1\,,\quad\Gamma^g_i\Gamma^h_j=\overline{\lan g,h\ran}\Gamma^h_j\Gamma^g_i\,,
\quad i<j\,.
\end{equation}

We extend the endomorphism $\rho$ on $F_Q^{{S}_Q}$ to  $F_Q^{G}$  by $W^g_i \mapsto W^g_{i+1}$, i.e. 
$V^g_i \mapsto U^{-g}_iU^g_{i+1}$ and $ U^{-g}_iU^g_{i+1} \mapsto V^g_{i+1}$. Next, extend $\rho$ to $F_Q$ by 
$U^g_i \mapsto V^g_1 V^g_2 \cdots V^g_i$.  In particular $U^g_1 \mapsto V^g_1$. In the Cuntz algebra description of $\mathcal{O}_G$,
with generators $S_g$, this means that the spectral projection corresponding to $1$ namely $S_0S_0^*$ is taken to  $\Sigma_{h,k} S_h S_k^*/|G|$. So one could try as an ansatz to define an endomorphism $\rho$ on $\mathcal{O}_G$ by $S_0 \rightarrow \pm(\Sigma_h S_h/{\sqrt{|G|}})$, and similarly since $S_g S_g^*$  is taken to  $Ad(U(g)) (\Sigma_{h,k} S_h S_k^*)$ 
so we should take   $\rho(S_g) = \pm Ad(U(g)) (\Sigma_h S_h/{\sqrt{|G|}})$,
on the Cuntz algebra. Here  $g\mapsto U(g) =  \Sigma_h \langle g,h \rangle S_h S_h^*$  is basically the regular unitary representation of $G$ in $\mathcal{O}_G$. Then indeed $\rho$  is a well defined endomorphism on the Cuntz algebra induced by high temperature - low temperature duality. The sign choice $\pm$ determines the two Tambara-Yamagami categories associated to the same group and non-degenerate pairing.
Denote by  $\alpha_g$ the  automorphisms on $\mathcal{O}_G $
defined by  $\alpha_{g}(S_h) = S_{g+h}$. We have recovered the realisation of $\cT\cY_\pm(G,\lan,\ran)$ given in Theorem 4.

\subsection{The double of Tambara-Yamagami} 

Fix any Tambara-Yamagami fusion category $\cT\cY_s(G,\lan,\ran)$.  Let $q$ be any quadratic form on $G$ realising $\lan\cdot,\cdot\ran$ (Lemma 4(a) says that such $q$  always exist).

$\cT\cY_s(G,\lan,\ran)$  is realised as a system of endomorphisms in a Cuntz algebra in Theorem 4 above. From this the modular data of its double can be computed (c.f. section 3 of \cite{iz3}). 
We find there are precisely $4n+n(n-1)/2$ primaries (equivalence classes of simples), which we will parametrise by $\beta^g_i,\rho^g_j,\sigma_{g,h}=\sigma_{h,g}$,
where $g,h\in G$, $g\ne h$, and $i,j\in\{0,1\}$. For each $g\in G$, fix a square-root $\sqrt{q(g)}$. 
Then the modular data for the double is 
$$T_{\beta^g_i,\beta^g_i}=\lan g,g\ran\,,\ T_{\rho^g_i,\rho^g_i}=(-1)^i(\pm x^3)^{-1/2}/\sqrt{q(g)}\,,\ \ T_{\sigma_{g,h},\sigma_{g,h}}=\lan g,h\ran\,,$$
$$S_{\beta^g_i,\beta^h_j}=\frac{\overline{\lan g,h\ran}^2}{2{n}}\,,\ S_{\beta^g_i,\rho^h_j}=(-1)^i\frac{\overline{\lan g,h\ran}}{2\sqrt{n}}\,,\ S_{\beta^g_i,\sigma_{h,k}}=\frac{\overline{\lan g,h+k\ran}}{n} \,,\ S_{\sigma_{g,h},\rho^k_j}=0\,, $$ $$S_{\sigma_{g,h},\sigma_{g',h'}}=\frac{\overline{\lan g,h'\ran\lan h,g'\ran+\lan g,g'\ran\lan h,h'\ran}}{n}\,,\ S_{\rho^g_i,\rho^h_j}=\frac{(-1)^{i+j}x^{-3}}{2{n}{\sqrt{q(g)}\sqrt{q(h)}}}\sum_k\lan k-g-h,k\ran\,.$$ 
where $x$ is defined in section 4.1.

Different choices of $\sqrt{q(g)}$ in the previous paragraph are absorbed into the choice of which $\rho^g_i$ we call $i=0$ or $i=1$.

{Using Gauss sums, it is possible to evaluate the sum in $S_{\rho^g_i,\rho^h_j}$ more explicitly. In particular, write $\langle\cdot,\cdot\rangle$ as a product of the  indecomposables of section 4.1. Then the sum in $S_{\rho^g_i,\rho^h_j}$ will be the product of the corresponding sums for each indecomposable. These sums for each indecomposable are readily computed: e.g. 
for \textbf{type} $p^k_s$ we obtain
$$\sum_{\ell=0}^{p^k-1}\lan \ell-a,\ell\ran=\eps_{p^k}^{-1}s^k\overline{\lan a/2,a/2\ran}\sqrt{p^k}$$
whilst for \textbf{type}  $2^k_m$ we compute
$$\sum_{\ell=0}^{2^k-1}\lan \ell-a,\ell\ran=\left\{\begin{array}{cc}0&\mathrm{if}\ k>1\ \mathrm{and}\ a\ \mathrm{odd}\\ 
(1-i)\eps_m\sqrt{2^k}\left(\frac{2}{m}\right)^k\overline{\lan n/2,n/2\ran}&\mathrm{if}\ k>1\ \mathrm{and}\ a\ \mathrm{even}\\ 
2&\mathrm{if}\ k=1\ \mathrm{and}\ a\ \mathrm{odd}\\ 
0&\mathrm{if}\ k=1\ \mathrm{and}\ a\ \mathrm{even}\end{array}\right.$$
}

We compute (from Verlinde's formula \eqref{Fusl}) the fusions
$$\beta_s^g\cdot \beta_{s'}^h=\beta_{ss'}^{g+h}\,,\ \beta_s^k\cdot\sigma_{g,h}=\sigma_{g+k,h+k}\,,\ \beta_s^g\cdot\rho_{s'}^h=\rho_{ss'}^{h+2g}\,,$$
$$\sigma_{g,h}\cdot\sigma_{g',h'}=\sigma_{g+g',h+h'}+ \sigma_{g+h',h+g'}\,,\ \sigma_{g,h}\cdot
\rho_s^k=\rho_+^{g+h+k}+\rho_-^{g+h+k}\,,$$ 
$$ \rho_s^g\cdot\rho_{s'}^h=\sum_{[k]\ne[(g+h)/2]}\sigma_{k,g+h-k}+\sum_{\ell, 2\ell=g+h}\beta_{ss'}^\ell$$
where the first sum in the $\rho_s^g\cdot\rho_{s'}^h$ fusion product is over the size-2 orbits of $k\leftrightarrow g+h-k$ (to avoid over-counting $\sigma$'s), and the second sum is over the size-1 orbits $\ell$. On the right side of the $\sigma_{g,h}\cdot\sigma_{g',h'}$ fusion product we use the convention $\sigma_{k,k}=\beta_+^k+\beta_-^k$.

Note that the Tambara-Yamagami categories depend on $\langle\cdot,\cdot\rangle$ and not $q$, but we described the modular data of their doubles in terms of $q$. As explained in Lemma 4, different $q$ can correspond to the same $\langle\cdot,\cdot\rangle$. Is it clear that they will give rise to the same modular data? This point is somewhat subtle. If $q,q'$ are equivalent, then so will be the corresponding modular data. However, compare $q$ of \textbf{type} $2^2_1$ with $q'$   of \textbf{type} $2^2_{-3}$ (these are inequivalent but correspond to the same pairing $\langle,\rangle$). The $T$ matrices using $q$ vrs $q'$ are identical, except for the 8 entries for $\rho_i^g$: when $s=+$, using $q$ gives exactly 2 entries $T_{\rho_i^g,\rho_i^g}$ equal to 1, whereas using $q'$ none equal 1. For $s=-$, these are reversed. We find that if $q$ describes the modular data of the double of   $\cT\cY_s(G,\lan,\ran)$, then $q'$ describes that of $\cT\cY_{-s}(G,\lan,\ran)$.

\subsection{The $\bbZ_2$-crossed braiding of Tambara-Yamagami}

There is a simpler way to associate to $\cT\cY_s(G,\lan,\ran)$ a modular tensor category (when $|G|$ is odd).

Given a group $\Gamma$, a $\Gamma$-crossed category $\cC$ (c.f. \cite{Tur,Mug,DGNO}) is a modoidal category with a $\Gamma$-action by automorphisms and a $\Gamma$-valued grading $\partial$ such that $\partial({}^gX)=g\,\partial X\, g^{-1}$. A $\Gamma$-braiding on a $\Gamma$-crossed category is a choice, for  any  objects $\lambda,\mu\in\cC$ where $\lambda$ is homogeneous, of an operator
$\varepsilon_{\lambda,\mu}\in\mathrm{Hom}(\lambda\mu,{^{\partial(\lambda)}\mu}\, \lambda)$ subject to
initial conditions
$$
\varepsilon_{id_A,\mu}=\varepsilon_{\lambda,{id_A}}=\bfe \,,
\label{ini}
$$
and   the naturality and 
 {braiding-fusing} equations (which for later convenience we write in the endomorphism language)
\begin{equation}
\begin{array}{rl}
\varepsilon_{\nu,\kappa}s\lambda(t)&=\,\,{}^{\partial \lambda}t\,{}^{\partial}\mu(s)\varepsilon_{\lambda,\mu}\,, \\ [.4em]  
\varepsilon_{\lambda,\mu\circ\kappa}&=\,\,{}^{\partial\lambda}\mu(\varepsilon_{\lambda,\kappa})\varepsilon_{\lambda,\mu}\,, \\ [.4em]  
\varepsilon_{\lambda\circ\nu,\mu}&=\,\,{\varepsilon_{\lambda,{}^{\partial\nu}\mu}}\lambda(\varepsilon_{\nu,\mu})\,, \\ [.4em]  
\gamma(\varepsilon_{\lambda,\mu})&=\,\,\varepsilon_{\gamma(\lambda),\gamma(\mu)}\,, 
\end{array}
\label{BFE}
\end{equation}
whenever $\lambda,\mu,\nu,\kappa\in\cC$ (with $\lambda,\nu$ homogeneous), $s\in\mathrm{Hom}(\lambda,\nu),t\in\mathrm{Hom}(\mu,\kappa)$, and $\gamma\in\Gamma$.

If $\cC$ carries an action of $\Gamma$, then the $\Gamma$-equivariantization $\cC^\Gamma$   is a tensor category whose simple objects are pairs $(X,\{u_\gamma\}_{\gamma\in\Gamma})$ where $X$ is an orbit of $\Gamma$ and $u_\gamma$ for each $\gamma\in\Gamma$ is an isomorphism in Hom$(g(X),X)$ --- see \cite{Mug} for details.  The point is that if $\cC$ is  $\Gamma$-crossed braided, then the $\Gamma$-equivariantisation $\cC^{\Gamma}$ will be braided \cite{Mug,DGNO}.

Now, if we realise $\cC$ as endomorphisms on a Cuntz algebra $\cO$ (recall Theorem 4), then any C$^*$-algebra automorphism $\gamma$ extends to an automorphism of $\cC$: endomorphism $\rho$ gets mapped to $^\gamma \rho=\gamma\circ\rho\circ\gamma^{-1}$.

For example, {$\cO_n$} has an order-2 automorphism $\varphi$ sending $S_h$ to $S_{-h}$. Then we find $^\varphi\alpha_g=\alpha_{-g}$, $\varphi( U(g))=U(-g)$, and $^\varphi \rho=\rho$: e.g.
$$(^\varphi \rho)(S_h)=\varphi(\rho(\varphi^{-1}(S_h)))=\varphi(U(-h)\frac{1}{\sqrt{|G|}}\sum_{k\in G }S_kU(h)^*)=U(h)\frac{1}{\sqrt{|G|}}\sum_{k\in G} S_{-k}U(-h)^*=\rho(S_h)$$
Hence $\Gamma=\lan\varphi\ran\cong\bbZ_2$ acts on any $\cT\cY_s(G,\lan,\ran)$. The corresponding $\mathbb Z_2$-grading has $\rho$ odd and the $\alpha_g$  even.

\medskip\noindent\textbf{Theorem 5.} \textit{Any Tambara-Yamagami category  $\cT\cY_s(G,\lan,\ran)$ is $\bbZ_2$-crossed braided. This $\bbZ_2$-braiding is non-degenerate iff $|G|$ is odd. In this case, the $\bbZ_2$-equivariantization $\cT\cY_s(G,\lan,\ran)^{\bbZ_2}$ has primaries $\sigma_g=\sigma_{-g}$ ($g\ne 0$), $\beta^\pm$, $\rho^\pm$, and modular data} 
\begin{equation}
\begin{array}{rl}
S_{\sigma_g,\sigma_h}=&2\lambda\,(\langle g,h\rangle^2+ \langle g,h\rangle^{-2})\,,\ 
 S_{\beta^t,\beta^{t'}}=\lambda\,,\ S_{\sigma_g,\rho^t}=S_{\rho^t,\sigma_g}=0\,,\\ [.4em]
S_{\rho^t,\rho^{t'}}=&(-1)^{(|G|^2-1)/2}tt'/2\,,\ 
 S_{\sigma_g,\beta^{t}}=S_{\beta^t,\sigma_g}=2\lambda\,,\ 
S_{\rho^t,\beta^{t'}}=S_{\beta^{t'},\rho^t}=tt'/2\,,\\ [.4em]
 &T_{\sigma_g,\sigma_g}=\langle g,g\rangle\,,\ 
T_{\rho^t,\rho^{t}}=t/\sqrt{sx^3}\,,\ T_{\beta^t,\beta^t}=1\,,\end{array}
\end{equation}
where $\lambda^{-1}=2|G|$.\medskip

\noindent\textit{Proof.} First we must verify the $\bbZ_2$-braiding, for the action and grading defined in the paragraph before Theorem 5. As explained more generally in section 2.2 of \cite{BEK1}, to do this, it suffices to consider equivariance $\gamma(\varepsilon_{\lambda,\mu})=\varepsilon_{\gamma(\lambda),\gamma(\mu)}$ and the braiding-fusing relation
\begin{equation}\label{braidfus} \varepsilon_{\nu,{}^{\partial \kappa}\mu}\,\nu(\varepsilon_{\kappa,\mu})\,s={}^{\partial \lambda}\mu(s)\,\varepsilon_{\lambda,\mu}\end{equation}
and its adjoint,
for all $\gamma\in\Gamma$ and $s\in\mathrm{Hom}(\lambda,\nu\circ \kappa)$, as $\kappa,\lambda,\mu,\nu$ run through 
representatives of the equivalence classes of simples. In our case, we take those representatives to be $\alpha_g$ and $\rho$.

Our starting point is the realisation (Theorem 4) of Tambara-Yamagami as a system of endomorphisms on the Cuntz algebra $\cO_G$. We have $\cO_G$-algebra endomorphisms $\rho,\alpha_g$ and intertwiners $U(g),S_h\in\cO_G$, such that
Hom$(\rho,\alpha_g\rho)=\bbC 1$, Hom$(\alpha_g,\rho^2)=\bbC S_g$, and Hom$(\rho,\rho\alpha_g)=\bbC U(g)$. Taking adjoints, this implies Hom$(\alpha_g\rho,\rho)=\bbC 1$, Hom$(\rho^2,\alpha_g)=\bbC S_g^*$, Hom$(\rho\alpha_g,\rho)=\bbC U(g)^*$. 

The relevant intertwiner spaces are Hom$(\alpha_g \alpha_h , \alpha_h \alpha_g) = \bbC 1$, Hom$(\alpha_g \rho , \rho \alpha_{g}) = \bbC U(g)^*$, Hom$(\rho \alpha_g, \alpha_{-g} \rho ) = \bbC U(-g)$, Hom$(\rho \rho , \rho \rho) =  \sum_g {\mathbb C} S_gS_g^*$, so that  $\mathbb Z_2$-crossed braidings would be $\varepsilon_{\alpha_g, \alpha_{h}} = \epsilon_{g, {h}}1$, $\varepsilon_{\rho, \alpha_g} = \epsilon_{\rho, g} U(g)^*$, $\varepsilon_{\alpha_g, \rho } = \epsilon_{ g, \rho } U(g)$, and
$\varepsilon_{\rho, \rho} = \sum_k \epsilon_{\rho, \rho}^k S_k S_k^*, $ where $\epsilon_{g, {h}}, \epsilon_{\rho, g},  \epsilon_{ g, \rho } , \epsilon_{\rho, \rho}^k \in {\mathbb T}$ (since the category is unitary).
The  braiding-fusing  relations are then:
\begin{eqnarray}
\label{cb1}
& 
 \epsilon_{g,h} \,\epsilon_{g, k}=\epsilon_{g,h+ k}\,,\ \  \epsilon_{g,k} \,\epsilon_{h, k}=\epsilon_{g+h, k}\,,
\\
\label{cb7}
&\epsilon_{g+h,\rho}=\epsilon_{g,\rho} \epsilon_{h,\rho}\overline{\lan g,h\ran} \,,\ \ \epsilon_{\rho,g+h}=\epsilon_{\rho,g} \epsilon_{\rho,h}\overline{\lan g,h\ran} \,,
\\
\label{cb2}
&\lan g,k\ran \epsilon_{g,\rho} \,\epsilon_{\rho,\rho}^{k-g}=\epsilon^k_{\rho, \rho} =\overline{ \lan g,k\ran} \epsilon_{\rho,g} \,\epsilon_{\rho,\rho}^{g+ k} \,,
\\
\label{cb6}&\epsilon_{g,k}=\langle k,g\rangle\\
\label{cb9}
&\epsilon_{\rho,-k}\,\epsilon_{\rho,k}\langle k,g+k\rangle=\varepsilon_{g,k} 
\\ \label{cb10}
&\frac{s}{\sqrt{|G|}}\sum_\ell\epsilon^\ell_{\rho,\rho}\epsilon^k_{\rho,\rho}\langle \ell,k-g\rangle=\langle g,k\rangle\epsilon_{g,\rho}
\end{eqnarray} while equivariance gives \begin{equation} \label{cb3}
\epsilon_{g, h} =  \epsilon_{-g,-h} \,,\ \epsilon_{\rho,g}=\epsilon_{\rho,-g}\,,\  \epsilon_{g,\rho}=\epsilon_{-g,\rho}\,,\ \epsilon_{\rho,\rho}^k=\epsilon_{\rho,\rho}^{-k}\,. \end{equation}
Equation \eqref{cb6} says $\epsilon_{g,h}={\lan g,h\ran}$. Equation \eqref{cb7} and  \eqref{cb3} say $\epsilon_{g,\rho}=a(g)$  for some quadratic form $a$  associated to $\lan\cdot,\cdot\ran$. {Equation \eqref{cb2} and \eqref{cb3} force   $\epsilon_{\rho,h}=\epsilon_{h,\rho}$}. Then \eqref{cb2} gives $\epsilon^k_{\rho,\rho}=\epsilon^0_{\rho,\rho}\overline{a(k)}$. The value of $\epsilon_{\rho,\rho}^0$ is fixed, up to a sign, by \eqref{cb10}: we find $\epsilon^{0\ 2}_{\rho,\rho}=s\,x^{-3}$ where $x$ is defined in section 4.1. It is readily seen that these values satisfy all equations.

Now consider the $\bbZ_2$-equivariantisation $\cT\cY_s(G,\lan,\ran)^{\bbZ_2}$. 
Its simple objects are $\sigma_g=\sigma_{-g}=\alpha_g+\alpha_{-g}$ (for any $g\in G$ with $2g\ne 0$), $\rho^\pm=(\rho,\pm1)$, and $\beta_h^\pm=(\alpha_,\pm 1)$ whenever $2h=0$. It is convenient to write $\sigma_g=\beta_g^++\beta_g^-$ when $2g=0$. Fusion products are (c.f. Definition 2.1 of \cite{Mug}) 
\begin{equation}[\sigma_g][\sigma_h]=[\sigma_{g+h}]+[\sigma_{g-h}]\,,\ [\rho^t][\sigma_g]=2[\rho^t]\,,\ [\rho^t][\rho^{t'}]=\sum_h[\beta_h^{tt'}]+\sum_{[g]}[\sigma_g]\,,\nonumber\end{equation} 
\begin{equation}  [\sigma_g][\beta_h^t]=[\sigma_{g+h}]\,,\ [\rho^t][\beta_h^{t'}]=[\rho^{tt'}]\,,\ [\beta_h^t][\beta_{h'}^{t'}]=[\beta_{h+h'}^{tt'}]\,,\label{TYZ2fus}\end{equation}
 for any $t,t'\in\{\pm\}$ (in the case of greatest interest here, when $|G|$ is odd, the only $\beta^\pm_h$ are $\beta_0^\pm$). Hence the categorical dimensions (which must equal the Perron-Frobenius dimensions, since the category is unitary) of the simple objects $\beta_h^t$, $\sigma_g$ and $\rho^t$ are 1, 2 and $\sqrt{|G|}$, respectively.

The braiding $c_{x,y}$ on the $\bbZ_2$-equivariantisation  is as follows (c.f. Proposition 2.2 of \cite{Mug}). Choosing the obvious bases, we get
\begin{eqnarray}&c_{\sigma_g,\sigma_h}=\mathrm{diag}(\epsilon_{g,h},\epsilon_{-g,h})=\mathrm{diag}(\langle g,h\rangle,\overline{\langle g,h\rangle})\in\mathrm{End}(\sigma_{g+h}+\sigma_{g-h})\\ &c_{\beta^t_h,\beta^{t'}_{h'}}=\langle h,h'\rangle\\ &c_{\rho^t,\rho^{t'}}=\mathrm{diag}_g(t'\,\epsilon^0_{\rho,\rho}\overline{a( g)})\\ &{c_{\sigma_g,\rho^t}=\mathrm{diag}(a(g),a(g))\,,\ c_{\rho^t,\sigma_g}=\left(\begin{array}{cc} 0&t\,a(g)\\ t\,a(g)&0\end{array}\right)}\\ &
c_{\sigma_g,\beta_h^t}=\langle g,h\rangle=c_{\beta_h^t,\sigma_g}\\ &
c_{\rho^t,\beta^{t'}_{h}}= t'\,a(h)\,,\ c_{\beta^{t'}_{h},\rho^t}=a(h)\end{eqnarray}
For example, the difference between $c_{\sigma_g,\rho^t}$ and $c_{\rho^t,\rho^{t'}}$ is due to the fact that the braidings $\varepsilon_{x,y}$ lies in Hom$(x\otimes y,{}^{\partial x}y\otimes x)$.
Thus up to a global normalisation $\lambda$, we obtain the $S$-matrix:
\begin{eqnarray}&S_{\sigma_g,\sigma_k}=\lambda 2\,(\langle g,k\rangle^2+\langle g,k\rangle^{-2})\,,\ \ S_{\beta^t_h,\beta^{t'}_{h'}}=\lambda\,,\\ & S_{\sigma_g,\rho^t_h}=S_{\beta^t_h,\sigma_g}=\lambda 2\,,\ \ S_{\sigma_g,\rho^t}=S_{\rho^t,\sigma_g}=0\,,\ \ S_{\rho^t,\beta^{t'}_h}=S_{\beta^{t'}_h,\rho^t}=0\end{eqnarray} 
For example, the extra $\sqrt{|G|}$ in $S_{\rho^t,\beta^{t'}_h}$ comes from the dimension of $\rho^t$, which enters here  through the trace of the identity endomorphism on $\rho^t$.  Also, note that $\langle g,h\rangle^2=1$ when $2h=0$.
  It is explained after equation (10) in \cite{BEK1} how to obtain the $T$ matrix directly from the braidings. The result is as in the statement of Theorem 5.

Note that if 2 divides $|G|$, then the $S$ matrix is degenerate. The reason is that only the two $\rho^t$ rows have the possibility of distinguishing the $\beta^{t'}_h$ columns. If $G$ is even order, then there will be at least two $h\ne 0$ with $2h=0$, hence at least four such $\beta_h^t$.  The calculation of $S_{\rho^t,\rho^{t'}}$ when $|G|$ is odd requires that we compare $\sum_ga(g)$ to $\sum_ga(g)^2$. Both of these are Gauss sums; their ratio equals the Jacobi symbol $\left(\frac{2}{|G|}\right)$, which equals $(-1)^{(|G|^2-1)/8}$. When $|G|$ is odd, the invertibility of $S$ follows from the calculation of the double given in section 5.2. \textit{QED to Theorem 5}

\medskip For $|G|$ odd, the category will be a modular tensor category. In this case there will be only $\beta^t_0$, which we will abbreviate to $\beta^t$. The tensor unit is $\beta^+$, and $\beta^-$ is a simple-current. This modular tensor category corresponds to a $\bbZ_2$-orbifold of a lattice theory associated to $G$, as we see next subsection.

\subsection{Reconstruction for Tambara-Yamagami}

In this subsection we construct a strongly  rational VOA $\cV$ and conformal net $\cA$ whose representation theories Mod$(\cV)$ and Rep$(\cA)$ are both tensor equivalent to the double of Tambara-Yamagami $\cT\cY_s(G,\langle,\rangle)$ (we also do the same for its $\bbZ_2$-equivariantisation when $|G|$ is odd). Being the double of a fusion category, we should be able to obtain $\cV$ and $\cA$ from some sort of orbifold of a holomorphic theory (i.e. one with trivial representation theory). However, the double of Tambara-Yamagami is manifestly not the double of a group fusion category Vect$_G^\omega$ (its quantum dimensions are not all integers). This implies that we cannot obtain $\cV$ and $\cA$ directly from a holomorphic theory through a group orbifold. However, we do the next best thing.

We obtain $\cV$ and $\cA$  by first performing an orbifold by $G$ of a self-dual lattice theory $\cV(\Lambda)$ (which is necessarily holomorphic), resulting in a different lattice theory $\cV(L)$, and then orbifolding $\cV(L)$ by an involution of $L$. Because the involution doesn't descend to an automorphism of $\cV(\Lambda)$, we cannot combine these two orbifolds into a single group orbifold. 

For a simple  example, first perform a $\bbZ_3$ orbifold of the $E_8$ lattice theory, giving the $A_2\oplus E_6$ lattice theory, followed by an orbifold of either the $A_2$ or $E_6$ piece (which one depending on the choice of $\langle\cdot,\cdot\rangle$) by the $v\mapsto -v$ lattice isometry; the result is the double of $\cT\cY_s(\bbZ_3,\langle,\rangle)$ for one choice of sign $s$. (The construction we give below is necessarily a little more complicated, in order that  it also works for even $|G|$.)

Let $L$ be any even positive-definite lattice of dimension $d$. 
We discussed the automorphisms of  the VOA $\cV(L)$ in section 4.4, and Lemma 6 describes one class of VOA orbifolds, namely by $a_x$ for $x\in\bbQ\otimes_\bbZ L$.

The other class of orbifolds we need are by isometries of $L$. The most familiar of these corresponds to the isometry $\theta:v\mapsto-v$, and and orbifolding by it is denoted
$\cV(L)^+$. The resulting VOA $\cV(L)^+$ is known to be rational, and its irreducible modules are (see \cite{AD}): the simple-currents 
$\beta_{[v]}^\pm$ for each $[v]\in G_L$ with $2[v]=[0]$; the dimension-2 simples $\sigma_{[v]}=\sigma_{[-v]}$ for each $[v]\in G_L$ with  $2[v]\ne 0$; and {$\rho^\chi_\pm$} for each $\chi\in R/2L$, where $R=\{v\in L:v\cdot L\supseteq 2\bbZ\}$ (these  come from twisted $\cV(L)$-modules). We have
$T_{\beta^\pm_{[v]},\beta^\pm_{[v]}}=q_L([v])e^{-\pi i d/12}=T_{\sigma_{[v]},\sigma_{[v]}}$ and $T_{\rho^\chi_\pm,\rho^\chi_\pm}=\mp e^{-\pi i d/8}$. The Frobenius algebra $A$ describing the  (simple-current) extension of $\cV(L)^+$ to $\cV(L)$ is $\beta_{[0]}^++\beta_{[0]}^-$. As always with type 1 systems, there are two equivalent alpha-inductions, so choose one of them. Then it  sends $\beta_{[v]}^\pm$ to the (local simple) $\cV(L)$-module $[g]$, $\sigma_{[v]}$ to the 
local $\cV(L)$-module $[v]+[-v]$, and each $\rho^\chi_\pm$ to a twisted $\cV(L)$-module $\rho^\chi$.
We use this orbifold to reconstruct below the equivariantisation $\cT\cY_s(G,\langle,\rangle)^{\bbZ_2}$. For the double, when $|G|$ is even, another isometry is needed.

\medskip\noindent\textbf{Theorem 6.} \textit{Consider any Tambara-Yamagami category $\cT\cY_s(G,\lan,\ran)$. Then both  the double of $\cT\cY_s(G,\lan,\ran)$ (for any $G$) and the $\bbZ_2$-equivariantization $\cT\cY_s(G,\lan,\ran)^{\bbZ_2}$
(when $|G|$ odd) are realized as Mod$(\cV)$ and Mod$(\cV')$, for certain completely rational VOAs $\cV$ and $\cV'$ which are $\bbZ_2$-orbifolds of lattice VOAs.}\medskip 

\noindent\textit{Proof.} Let $L$ resp.\ $L'$ be any positive-definite even lattices realising the pointed modular tensor category $\cC(G,q)$ resp.\ $\cC(G,\overline{q})$, as promised by Theorem 2. So $G=G_L$ and $q=q_L$. Define $\widehat{L}$ to be the lattice gluing $(L\oplus L\oplus
L'\oplus L')\lan[(g,g,g,g)]_{g\in G}\ran$ (recall the discussion in section 4.2), where we identify each $g\in G$ with the corresponding coset in
$G_L$ and $G_{L'}$. Then $\widehat{L}$ is even and positive-definite, and $G_{\widehat{L}}\cong G\times G$ has elements  $[g,h]:=(g+h,0,g,h)$, for all $g,h\in G$.
We find that $\widehat{L}$ has quadratic form $q_{\hat{L}}([g,h])=\overline{\lan g,h\ran}$.

Now, it is manifest that $\widehat{L}$ has  an isometry $\tau$ sending any $(a,b,c,d)\in \widehat{L}$
to $(a,b,d,c)$. So $\tau$ sends $[g,h]\in G_{\widehat{L}}$ to $[h,g]$. Then (see \cite{BK,DX}) $\tau$ lifts to an
order-2 automorphism of the VOA or conformal net $\cV({\widehat{L}})$ resp. $\cA({\widehat{L}})$ (\textit{a priori} the order could also be 4, a complication due to the 2-cocycle implicit in the construction of $\cV(L)$, but for our $\tau$ the order is 2). We claim that the orbifold $\cA({\widehat{L}})^\tau$ has category of representations Rep$(\cA(\widehat{L})^\tau)$ which is braided tensor equivalent to the double of $\cT\cY_s(G,\lan,\ran)$ for some choice of sign.

To see this, let's first note (using \cite{DX}) that the modules $\beta^g_\pm,\sigma_{g,h},\rho^g_\pm$ for $\cA({\widehat{L}})^\tau$ {match} those for $\cT\cY_s(G,\lan,\ran)$. The nontrivial check are the twisted or solitonic modules $\rho_\pm^g$: the twisted $\cA({\widehat{L}})$-modules are in bijection with $\tau$-invariant classes in $G_{\hat{L}}$, i.e. the $[g,g]$ for $g\in G$, and when restricted to $\cA({\widehat{L}})^\tau$
each of these splits into what we call $\rho_+^g+\rho_-^g$. 

A third conformal net is relevant here. Consider the lattice gluing $\Lambda=\widehat{L}\lan [g,0]\ran_{g\in G}$: it is also even and positive-definite, and is self-dual. Thus the corresponding conformal net $\cA(\Lambda)$ is holomorphic, i.e. its category of representations is Vect. These three conformal nets are related by the conformal inclusions 
$$\cA({\widehat{L}})^\tau\subset \cA({\widehat{L}})\subset \cA(\Lambda)$$
These two inclusions correspond to Frobenius algebras $\theta_{12}=\beta_+^0\oplus\beta_-^0$ and $\theta_{23}=
\oplus_g[g,0]$, respectively. The algebra $\theta_{13}$, governing the extension $\cA(\widehat{L})^\tau\subset  \cA(\Lambda)$,  is determined shortly.
To identify Rep$(\cA({\widehat{L}})^\tau)$ with $\cD(\cT\cY_s(G,\lan,\ran))$ for some sign, we will identify the full system, i.e. Rep$(\cA({\widehat{L}})^\tau)_{\theta_{13}}$ (of course the local modules are 
Rep$(\cA({\widehat{L}})^\tau)_{\theta_{13}}^{loc}=\mathrm{Rep}(\cA({L_{sd}}))=\mathrm{Vect}$).

{As always, there are two different inductions, corresponding to using the braiding or its reverse, but in this type 1 (pure extension) setting they are equivalent, so choose either one.}
First consider alpha-induction from $\cA({\widehat{L}})^\tau$ to $\cA({\widehat{L}})$. Then clearly Ind$\,\beta_\pm^g=[g,g]$, Ind$\,\sigma_{g,h}=[g,h]\oplus[h,g]$ and Ind$\,\rho_\pm^g=\rho^g$, the $\cA({\hat{L}})$-twisted sectors. 

Hence the Frobenius algebra governing the $\cA({\widehat{L}})^\tau\subset\cA(\Lambda)$ extension is $\theta_{13}=\beta_+^0\oplus\beta_-^0\oplus\oplus_{g\ne 0}\sigma_{g,0}$. We compute dim$\,\mathrm{Hom}(\mathrm{Ind}\,\beta_i^g,\mathrm{Ind}\,\beta_j^h)=\mathrm{dim\,Hom}(\beta_i^g\otimes\theta_{13},\beta_j^h)=1$, so Ind$\,\beta_i^g$ is a simple object in the full system which we'll call $\alpha_g$ (these are also the twisted sectors for $\cA(\widehat{L})\subset\cA(\Lambda)$). Because alpha-induction is a tensor functor, we have $\alpha_g\otimes\alpha_h=\alpha_{g+h}$. Likewise,  dim$\,\mathrm{Hom}(\mathrm{Ind}\,\sigma_{g,h},\mathrm{Ind}\,\sigma_{g,h})=2$ and dim$\,\mathrm{Hom}(\mathrm{Ind}\,\beta_i^{k},\mathrm{Ind}\,\sigma_{g,h})=\mathrm{dim\,Hom}(\oplus_{\ell\ne 0}\sigma_{k+\ell,k},\sigma_{g,h})=\delta_{k,h}+\delta_{k,g}$, so Ind$\,\sigma_{g,h}=\alpha_g+\alpha_h$. Finally,
dim$\,\mathrm{Hom}(\mathrm{Ind}\,\rho_i^{g},\mathrm{Ind}\,\rho_j^{h})=\mathrm{dim\,Hom}(\oplus_{k,\pm}\rho_{\pm}^k,\rho_j^{h})=1$, so let $\rho$ denote this common simple object Ind$\,\rho_\pm^g$ in the full sytem. Again, using the fact that alpha-induction respects fusions, we find
$\alpha_g\otimes\rho=\rho=\rho\otimes\alpha_g$ and $\rho\otimes\rho=\oplus_g\alpha_g$. Thus the full system must be $\cT\cY_s(G,\lan,\ran')$ for some  symmetric pairing
$\lan\cdot,\cdot\ran'$ on $G$. By  Corollary 4.8 in \cite{BEK3}, the modular tensor category Rep$(\cA({\widehat{L}})^\tau)$ must be braided tensor equivalent to the double of $\cT\cY_s(G,\lan,\ran')$.  By comparing $T$-eigenvalues, we find that $\lan\cdot,\cdot\ran'=\lan\cdot,\cdot\ran$, as desired. Thus the full system must be $\cT\cY_s(G,\lan,\ran)$ for some sign $s$. 

To get the other sign, consider instead the lattice $\widehat{L}'=\widehat{L}\oplus E_8$ and the automorphism $(\tau,\alpha_v)$ where $v=([1],[1])\in A_1^*\oplus E_7^*\subset\frac{1}{2}E_8$. Then $\cV({\widehat{L}\oplus E_8})^{(\tau,\alpha_v)}$ is a simple-current extension of $\cV(\widehat{L})^\tau\otimes \cV({A_1\oplus E_7})$. Again, the full system is a Tambara-Yamagami category for $G$. We determine the pairing and sign from the modular data: in particular, we find that this $E_8$ trick changes the $T$-eigenvalues for the twisted modules $\rho^g_\pm$ by a factor of $i$, i.e. the sign $s$ has changed.

Now restrict to odd-order $G$. Similar considerations show that the conformal embedding $\cA(L)^+\subset\cA(L)$ has full system $\cT\cY_+(G,\langle,\rangle)$. This means (Corollary 4.8 in \cite{BEK3}) that the double of that Tambara-Yamagami is braided tensor equivalent to the Deligne product Rep$( \cA(L)^+)\stimes \mathrm{Rep}(\cA(L))^{opp}$  (here `opp' means the  category with the opposite braiding  $c^{opp}_{X,Y}=(c_{Y,X})^{-1}$). So
$$\mathrm{Rep}(\cA(L)^+)\stimes\cC(G,\overline{q})\cong \cZ(\cT\cY_+(G,\langle,\rangle))\cong \cT\cY_+(G,\langle,\rangle)^{\bbZ_2}\stimes\cC(G,\overline{q})$$
where the second tensor equivalence comes from the previous subsection. Now,  let $\cC$ resp.\ $\cC'$ be the full subcategory of Rep$(\cA(L)^+)\stimes\cC(G,\overline{q})$ resp.\ $\cT\cY_+(G,\langle,\rangle)^{\bbZ_2}\stimes\cC(G,\overline{q})$ generated by all objects which are invertible of odd order.  Both $\cC$ and $\cC'$  can be identified with the  $\cC(G,\overline{q})$ factors. Then the subcategories commuting with $\cC$ resp.\ $\cC'$ must be tensor equivalent. Thus we obtain
 $\mathrm{Rep}(\cA(L)^+)\cong \cT\cY_+(G,\langle,\rangle)^{\bbZ_2}$, as desired. The other sign is obtained from the $E_8$ trick, as done previously.
  \textit{QED to Theorem 6}\medskip

\subsection{$K$-theory of Tambara-Yamagami and its orbifold}

Until now, the examples we've studied $K$- and $KK$-theoretically are classical, either coming directly from finite groups or Lie theory.
In this subsection we realise  the  fusion category of $\cT\cY_s(G,\langle,\rangle)$ as a fusion category of bundles over a groupoid. This allows us to likewise  identify the modular tensor category of its double and (when $|G|$ is odd) its
$\bbZ_2$-equivariantisation, as modular tensor categories of bundles over groupoids. We explain how to realise  their  module categories using bundles. We hope in the future to find such geometric interpretations of other exotic fusion categories.

Let $G$ be a finite abelian group as before. Fix a non-degenerate symmetric pairing $\langle\cdot,\cdot\rangle$ on $G$. Using it, we can identify $G$ as a group with  $\widehat{G}$, its group of irreps, through $g\mapsto \hat{g}:=\langle g,\star\rangle\in\widehat{G}$. Likewise, we'll write $\hat{\psi}$ for the element in $G$ corresponding via $\langle,\rangle$ to $\psi\in\widehat{G}$.

Our first task is to capture the fusion ring \eqref{TYfusions} using groupoids. Recall the discussion in section 2.5.
Let $X=G\cup \mathrm{pt}$.  Consider the equivariant $K$-group  $K^\star_{G}(X)$, where $G$ acts on $G\subset X$ by left translation, and fixes $\textrm{pt}\in X$. There are two orbits in $X$: $G$ with trivial stabiliser, and pt with stabiliser  $G$. The indecomposable bundles with support pt are in natural bijection with $\psi\in\widehat{G}$ --- call them $a_\psi$ (we can use $\langle,\rangle$ to parametrise these by $g\in G$ if we like). There is only one indecomposable bundle with support $G$  (corresponding to the trivial representation of the trivial group), which we shall call $\rho$. Thus as an additive group, $K^0_{G}(G\cup \mathrm{pt})=K^0_{G}(G^L)\oplus K^0_{G}(\mathrm{pt})=\bbZ\oplus R_G$, whilst $K^1_{G}(G\cup \mathrm{pt})=0$.

Recall the discussion of product of bundles in section 2.5. Consider the map $M:X\times X\rightarrow X$ defined by  $M(g,\mathrm{pt})=M(\mathrm{pt}, g)=g$, $M(\mathrm{pt},\mathrm{pt})=\mathrm{pt}$, and 
$M(g,h)=\mathrm{pt}$ (this choice is necessary in order to recover the Tambara-Yamagami fusions). Then $M$ is $G$-equivariant. We need to modify slightly the generic product defined in section 2.5. The tensor product $V\otimes W$ of bundles carries naturally an action of $G\times G$ as mentioned there; use the pairing to identify this group with $G\times \widehat{G}$, and restrict this action to the diagonal $\Delta'_G=\{(g,\hat{g})\,|\,g\in G\}$. For the union of orbits $Y$  take 
$$Y=\mathrm{pt}\times\mathrm{pt}\,\cup\,\mathrm{pt}\times \widehat{G}\,\cup\, G\times\mathrm{pt}\,\cup\,\Delta'_G$$ which we identify with a subset of $X\times X$ in the obvious way, again using the pairing. 

Consider first the product $\rho\otimes \rho$. It will be a bundle over the orbit $M(0,0)=\mathrm{pt}$, namely Ind$_{\mathbf{1}}^G(\bbC\otimes\bbC)=\oplus_\psi\psi$, the regular representation of $G$, where the superscript $G$ is the stabilizer of pt, where the subscript $\mathbf{1}$ denotes the stabilizer of $0\in G\subset X$ in $G$, and we take the tensor product of the fibres over $0$.   Thus we recover the formula $[\rho]\cdot [\rho]=\sum_g[\alpha_g]$ in \eqref{TYfusions}. Shortly, we will need to understand this more explicitly: we can naturally identify the tensor product $\bbC[G]\otimes_G\bbC[\widehat{G}]$ of bundles with both $\bbC[\widehat{G}]$ and $\bbC[G]$, as we now show. Indeed, the representation $\psi\in\widehat{G}$ is isomorphic to  the submodule $\bbC\sum_g\psi(g)g\otimes \hat{g}$ of $\bbC[G]\otimes_G\bbC[\widehat{G}]$; likewise the basis element $\underline{e}_h\in\bbC[G]$ is identified with $e_h=\sum_g \hat{g}(h)e_g\otimes e_{\hat{g}}\in\bbC[G]\otimes_G\bbC[\widehat{G}]$.

 The product  $\psi\otimes\rho$ will be a bundle over the orbit of $M(0,\mathrm{pt})=0$, namely
$\bbC\otimes\mathrm{Res}^G_{\mathrm{1}}(\psi)=\bbC$. More explicitly, the total space $\bbC^\psi\otimes V_\rho$ has basis $v_\psi\otimes e_g$, which maps $G$-equivariantly to $\psi(g)\underline{e}_g$, a basis for $\bbC[G]$. Thus $\psi\otimes\rho\cong\rho$.

Likewise, the product $\rho\otimes\alpha_\psi$ will also be a bundle over the orbit of $M(0,\mathrm{pt})=0$, namely
$\bbC\otimes\mathrm{Res}^G_{\mathrm{1}}(\psi)=1$. More explicitly, $e_g\otimes v_\psi$ maps equivariantly to $e_{g\hat{\psi}}$. Thus again $\rho\otimes\alpha_g\cong\rho$. It will be  important in the following though that although the bundles $\bbC^\psi\otimes_G V_\rho$ and $V_\rho\otimes_G\bbC^\psi$  are equivalent, they are certainly not equal.

Finally, the product $\alpha_\psi\otimes\alpha_{\phi}$ will be a bundle over the orbit pt, with $\Delta_G'$-module $\bbC^\psi\otimes
\bbC^{\hat{\phi}}$ naturally identified with $G$-module $\bbC^{\psi\phi}$ (or equivalently $\bbC^{\hat{\psi}\hat{\phi}}$). Thus, $\alpha_\psi\otimes\alpha_{\phi}=\alpha_{\psi\phi}$. 

Thus for this choice of $M$ etc, we obtain a ring structure on $K^0_G(X)$ isomorphic to the fusion ring \eqref{TYfusions} of the Tambara-Yamagami category $\cT\cY_s(\widehat{G},\langle,\rangle)$, or equivalently that of $\cT\cY_s(G,\langle,\rangle)$. But much more is true. As we know, this ring structure is independent of the sign and choice of pairing on $G$. The point of introducing the pairing, and of tweaking the product of section 2.5, is to recover the full fusion category of Tambara-Yamagami. To do this, we need to identify the associators.

Consider first the associator $a_{\psi,\rho,\phi}$. First note that $(\bbC^\psi\otimes\rho)\otimes\bbC^\phi$ maps equivariantly to $\rho$ through $(v_\psi\otimes e_g)\otimes v_\phi\mapsto \psi(g)e_{g\hat{\phi}}$, whereas $\bbC^\psi\otimes(\rho\otimes\bbC^\phi)$ maps to $\rho$ through $v_\psi\otimes (e_g\otimes v_\phi)\mapsto \psi(g\hat{\phi})e_{g\hat{\phi}}$. Thus the natural map $(\bbC^\psi\otimes\rho)\otimes\bbC^\phi\to \bbC^\psi\otimes(\rho\otimes\bbC^\phi)$  is multiplication by $\psi(\hat{\phi})=\langle\psi,\phi\rangle$. We choose this to be the associator $a_{\psi,\rho,\phi}$.

For comparison, consider  the associator $a_{\psi,\phi,\rho}$. Then $(\bbC^\psi\otimes\bbC^\phi)\otimes\rho$ maps to $\rho$ through $(v_\psi\otimes v_\phi)\otimes e_g\mapsto \psi(g)\phi(g)e_{g}$, whereas $\bbC^\psi\otimes(\bbC^\phi\otimes\rho)$ maps to $\rho$ through $v_\psi\otimes (v_\phi\otimes e_g)\mapsto \psi(g)\phi(g)e_{g}$. Thus the natural map $(\bbC^\psi\otimes\bbC^\phi)\otimes\rho\to \bbC^\psi\otimes(\bbC^\phi\otimes\rho)$ (hence the associator $a_{\psi,\phi,\rho}$) is multiplication by 1. Likewise, $a_{\psi,\phi,\rho}=1=a_{\phi,\phi',\phi''}$.

Consider next the associator $a_{\rho,\psi,\rho}$. Since $(\rho\otimes\alpha_\psi)\otimes\rho\cong \sum_\phi\alpha_\phi\cong \rho\otimes(\alpha_\psi\otimes\rho)$, $a_{\rho,\psi,\rho}$ is a vector with $|G|$ components, each component indexed by some $\phi\in\widehat{G}$. First note that $e_g\times e_{\hat{g}}\in M^{-1}(\mathrm{pt})\cap Y$ contributes the vector $e_{g\hat{\psi}}\otimes e_{\hat{g}}$ to $(\rho\otimes\alpha_\psi)\otimes\rho$ and $\psi(g)e_g\otimes e_{\hat{g}}$ to $\rho\otimes(\alpha_\psi\otimes\rho)$, hence contributes to the $\phi$-component the terms $\sqrt{|G|}^{-1}\overline{\phi(g\hat{\psi})}e_{g\hat{\psi}}\otimes e_{\hat{g}}$ and $\sqrt{|G|}^{-1}\psi(g)\overline{\phi(g)}e_g\otimes e_{\hat{g}}$, respectively. Thus the $\phi$-component of the associator is multiplication by $\psi(\hat{\phi})=\langle\psi,\phi\rangle$. On the other hand all $|G|$ components of the associators $a_{\rho,\rho,\psi}$ and $a_{\psi,\rho,\rho}$ are 1.

Finally, consider the associator $a_{\rho,\rho,\rho}$. Since $(\rho\otimes\rho)\otimes\rho\cong |G|\rho\cong \rho\otimes(\rho\otimes\rho)$, $a_{\rho,\rho,\rho}$ is a $|G|\times |G|$ matrix with rows and columns indexed by  $\psi\in\widehat{G}$. The associator is the change of basis matrix from the natural basis $b_\psi(g)=\overline{\psi(g)}e_g$ of the multiplicity space of
$(\rho\otimes\rho)\otimes\rho$, to the natural basis $c_\phi(g)=e_{g\hat{\phi}}$ of the multiplicity space of $ \rho\otimes(\rho\otimes\rho)$. We find that $c_{\psi}$ corresponds in the second basis to the Fourier transform $\sqrt{|G|}^{-1}\sum_\phi \overline{\psi(g\hat{\phi})}e_{g\hat{\phi}}$, and thus the change of basis matrix is $a_{\rho,\rho,\rho}=\sqrt{|G|}^{-1}\langle \psi,\phi\rangle^{-1}$.

By comparison with \cite{TY}, we find that these associators match spot on with the associators of  $\cT\cY_+(\widehat{G},\langle,\rangle)$. Thus the category of bundles of $X/\!/G$ form a fusion category tensor equivalent to  $\cT\cY_+(\widehat{G},\langle,\rangle)$. It would be interesting to recover the $s=-1$ analogue as a groupoid. As mentioned earlier, $s=\pm$ corresponds to the twist of the associator coming from $H^3$ of the grading $\bbZ_2$ of Tambara-Yamagami, and so only appears in $a_{\rho,\rho,\rho}$ (as a sign change).

Our groupoid picture also constructs with ease some module categories. 
The module categories for $\cT\cY_s(G,\langle,\rangle)$ were classified by \cite{MeMu}, and come in two classes. The first class is parametrised by a subgroup $H\le G$ and any $\psi\in H^2(H,\bbT)$. Recall the treatment of the module categories for  Rep$(G)$ at the end of section 2.5.  
The irreducible $\psi$-twisted bundles $V$ over $X/\!/_\psi H$, where  $H$ acts on $X=G\cup\mathrm{pt}$ as usual, correspond to a $\rho_{[g]}$ for each coset $[g]\in H\backslash G$ (with support on $[g]$ and trivial stabiliser) and each $\chi\in\mathrm{Irr}_\psi(H)$ (with support on pt). Define $M'=M$ and $Y'=Y$ as before. 
 We compute $\phi\otimes\chi\cong \phi|_H\,\chi$, $\phi\otimes \rho_{[g]}\cong \rho_{[\hat{\phi}g]}$,  $\rho\otimes\chi\cong \sum_{[k]\in G/H}\rho_{[k]}$, and $\rho\otimes\rho_{[g]}$ is  the regular representation $\sum_{\chi\in \widehat{H}}\chi$.

Consider further the case $[\psi]=[1]$ for simplicity. Note that the module categories Bun$(X/\!/H)$ and Bun$(X/\!/H^\perp)$ are equivalent: through the isomorphisms $\widehat{H}\cong G/H^\perp$ and $\widehat{H^\perp}\cong G/H$ given by the pairing $\langle,\rangle$, the indecomposable bundles $\phi\in\widehat{H}$ for $X/\!/H$ correspond to the bundles $\rho_{[k]}$ of $X/\!/H^\perp$, whilst the $\rho_{[g]}$ of the former groupoid correspond to $\chi\in\widehat{H^\perp}$. Generalising this to arbitrary $\psi$ recovers Lemma 30 of \cite{MeMu}.

The other class is much more complicated to describe, and exist only for very special $G$ (e.g.\ $\sqrt{|G|}$ must be integral); it is testament to the $K$-theory method that it provides these with an elegant formulation. Suppose for instance that $G=A\times \widehat{A}$ with pairing $\langle(a,\hat{a}),(b,\hat{b})\rangle = \hat{a}(b)\,\hat{b}(a)$.  Then we can reinterpret the groupoid $X_G/\!/G$ as one copy of $A$ acting on the left and the other copy  acting on the right. The reason for this is that it gives another expression for the product of bundles: as before, first identify $K_{A^L\times \widehat{A}^R}(X)\times K_{A^L\times A^R}(X)$ with
$K_{A^L\times A^R\times A^L\times A^R}(X\times X)$. Then restrict to the subgroup $A^L\times 1\times 1\times A^R$. If we now have a map $M:X\times X\to X$ which is $G$-equivariant in the sense that $M(xa,y)=M(x,ay)$, then the wrong-way map composed with the other two produce a multiplication $K_{A^L\times A^R}(X)\times K_{A^L\times A^R}(X)\to  K_{A^L\times A^R}(X)$. In terms of the notation in section 2.5, choose the same $M$ and $Y$ as before, except that $\Delta'_G$ is replaced with $(A\times 1)\times(1\times \widehat{A})$. In this equivalent formulation of $\cT\cY_+(G,\langle,\rangle)$, the $\alpha_g$ correspond to pairs $(\phi,a)$ where $\phi\in\widehat{A}$ and $a\in A$, with the obvious product.

Given any subgroup $B\le A$, we can interpret $\widehat{A/B}$ as the subgroup of $\psi\in\widehat{A}$ with kernel ker$(\psi)\ge B$. The point is that for any subgroup $B\le A$, the  category Bun$(\mathrm{pt}/\!/B\times \widehat{A/B})$ naturally forms a module category for the fusion category Bun$(X/\!/A\times\widehat{A})\cong\cT\cY_+(G,\langle,\rangle)$, for any subgroup $B\le A$, where now we choose the constant map $M':X\times \mathrm{pt}\to\mathrm{pt}$, and $Y'=\mathrm{pt}\times\mathrm{pt}\cup B\times \widehat{A/B}$. The module category has simples $(\chi,a'B)$ where $\chi\in\widehat{B}$ and $a'B\in A/B$; the nimrep is $(\phi,a).(\chi,a'B)=(\phi|_B\chi,aa'B)$ and $\rho.(\chi,a'B)=\sum_{\chi',a''B}(\chi',a''B)$ where the sum is over all $\chi\in\widehat{B}$ and all $a''B\in A/B$. 

Now turn to the $\bbZ_2$-equivariantisation $\cT\cY_s(G,\langle,\rangle)^{\bbZ_2}$. Here, $G$ must be odd order.  We wish to identify its fusion ring with the equivariant $K$-group $K_{D_G}(X)$, where $D_G$ is the generalised dihedral group $G\sdprod\bbZ_2$  ($\bbZ_2$ acts on $G$ by taking inverse) and $X=G\cup\mathrm{pt}$ as before.  $D_G$ fixes pt. {The additive group structure is again easy: 
$$K_{D_G}(X)=K_{\bbZ_2}(0)\oplus K_{D_G}(\mathrm{pt})=R_{\bbZ_2}\oplus R_{D_G}$$
(More precisely, what should appear  is $R_{D_{\widehat{G}}}$, which we identify with $D_G$ using $\langle,\rangle$ as usual.)
We identify $\rho_\pm$ with $R_{\widehat{\bbZ_2}}$, so these correspond to bundles at $0\in G$. The irreducible representations of $\widehat{D_G}$ consist of two one-dimensional representations, which we identify with $\beta^\pm$, and for each $\psi\in \widehat{G}$, $\psi\ne 1$, we have the two-dimensional representation $\sigma_{\psi}=\sigma_{\overline{\psi}}$. We identify $\sigma_g=\sigma_{-g}$ with $\sigma_{\hat{g}}$.}

To recover the multiplicative structure, use the same $M$ as before.
To compute the products, use Ind$_{\bbZ_2}^{D_G}\pm=\beta^{\pm}+\sum_\psi \sigma_\psi$ and $\sigma_\psi\otimes \sigma_{\psi'}=\sigma_{\psi\psi'}+\sigma_{\psi'\overline{\psi}}$. Otherwise the calculation is as for Tambara-Yamagami. Note that the category  Bun$(X/\!/{D_G})$ naturally possess a fusion category (in fact modular tensor category) structure, by applying the $\bbZ_2$-equivariantisation procedure to the category Bun$(X/\!/G)$.

A large class of module categories for the $\bbZ_2$-equivariantisation are Bun$(X/\!/_\psi H)$, where $H\le G$ and $[\psi]\in H^2_H(\mathrm{pt};\bbT)$. The indecomposable bundles, as we know, are $\psi\in\widehat{H}$ and $\rho_{[kH]}$. The nimrep is $\beta^\pm.\psi=\psi$, $\beta^\pm.\rho_{[kH]}=\rho_{[kH]}$, $\sigma_g.\psi=\hat{g}|_H\psi+\overline{\hat{g}}|_H\psi$, $\sigma_g.\rho_{kH]}=\rho_{[(g+k)H]}+\rho_{[(k-g)H]}$, $\rho^\pm.\psi=\sum_{kH}\rho_{[kH]}$ and $\rho^\pm.\rho_{[kG]}=\sum_\psi\psi$.
The modular invariant is given by the correspondence
$$\begin{tikzpicture}
  \matrix (m) [matrix of math nodes,row sep=2em,column sep=2em,minimum width=2em]
  {
    X/\!/ D_G&&&&&&  X/\!/ D_G \\ &X/\!/G&& \mathrm{pt}/\!/Z &&X/\!/G&\\ &&\mathrm{pt}/\!/G&&\mathrm{pt}/\!/G&&\\ };
  \path[-stealth] (m-1-1) edge node [left] {$\iota\ $} (m-2-2) (m-1-7) edge node [right] {$\ \iota$} (m-2-6)
    (m-2-2) edge node [left] {$\iota'\ $} (m-3-3)  (m-2-6) edge node [right] {$\ \iota'$} (m-3-5) (m-2-4) edge node [left] {$p_+$} (m-3-3)   (m-2-4) edge node [right] {$p_-$} (m-3-5) ;
\end{tikzpicture}$$
where notation is taken from section 4.5 (in particular $Z$ depends on $H,\psi$) or is otherwise clear. These module categories correspond to first extending both sides by the simple-current $\beta^-$, resulting in the lattice theory  of Theorem 6. The type 1 parents are thus lattice theories.
The most important of these choices is  $G=H$ and $\psi=1$, with modular invariant $|\chi_{\beta^+}+\chi_{\beta^-}|^2+2\sum_{\pm g\ne 0}|\chi_{\sigma_g}|^2$. It undoes the $\bbZ_2$-orbifold of Theorem 6.

We also get a groupoid interpretation of the double of $\cT\cY_s(G,\langle,\rangle)$, namely the bundles over $G\cup \{\mathrm{pt}\}/\!/_{(G\times G)\sdprod\bbZ_2}$, where $\bbZ_2$ here permutes $(g,h)\leftrightarrow(h,g)$. $(g,h,0)\in (G\times G)\sdprod\bbZ_2$ sends $\gamma\in G$ to $g+\gamma -h$ whilst $(0,0,1)$ sends $\gamma$ to $-\gamma$. Module categories can be easily constructed as we just did for the $\bbZ_2$-equivariantisation.

\subsection{Tambara-Yamagami and grafting}

In \cite{EG2} the authors introduced the notion of \textit{grafting}, and used it to generalise the modular data of e.g. the double of the Haagerup subfactor.  The key ingredient there is the notion of  $\bbZ_2$-laminated modular data. In \cite{EG2}, the modular data for the Haagerup double was constructed from the loop group LSpin$(13)$ at level 2 together with the double $\cD(S_3)$.

Based on the work of this section, a much more natural approach suggests itself. First, note that any
 $\bbZ_2$-equivariantisation, e.g.\ that of $\cT\cY_s(G,\langle,\rangle)$, is $\bbZ_2$-laminated. Moreover, the double of the Haagerup can be obtained from grafting 
$\cT\cY_s(\bbZ_3\times\bbZ_3,\langle,\rangle)^{\bbZ_2}$ and $\cT\cY_{s'}(\bbZ_{13},\langle,\rangle)^{\bbZ_2}$ for certain choices of signs and symmetric pairings.
This generalises to all known examples $\cD^0\mathrm{Hg}_\nu$ in the Haagerup-Izumi series. 
Namely, replace $\bbZ_3$ with $\bbZ_\nu$ and $\bbZ_{13}$ with $\bbZ_{\nu^2+4}$. Here $\nu$ can be any odd positive integer.
 
Grafting the Haagerup-Izumi series from the $\bbZ_2$-equivariantisations of two Tambara-Yamagami categories seems  much more natural and promising  than the suggestion of  \cite{EG2}, which proposed grafting finite group doubles to the loop group categories at level 2.

Recall that the $\bbZ_2$-equivariantisation of $\cT\cY_s(\bbZ_3\times\bbZ_3,\langle,\rangle)$
 has  simple objects $\sigma_\alpha\,,$ where $\alpha \in \bbZ_3\times\bbZ_3$ and $\alpha \sim -\alpha \ne 0\, $, simple-currents $\beta^{\pm}$ generating $\bbZ_2$, and twisted fields $\rho^{\pm}$ (`twisted' is in reference to the $\bbZ_2$-orbifold of the lattice theory, given in Theorem 6).  Its fusions are given in \eqref{TYZ2fus}.

The fusion rules for the double of the Haagerup can be read off from section 3.2 of \cite{EG2}. They can be written in  a more coherent (and general) form as follows.
The grafted  fusion rules amalgamate the simple-currents $\beta^\pm$ of the $\bbZ_\nu\times\bbZ_\nu$ and $\bbZ_{\nu^2+4}$ theories into simples called $\mathbf{0}$ and $\frak{b}$ in \cite{EG2}, amalgamate the orbifold fields $\sigma_\alpha$ from the $\bbZ_\nu\times\bbZ_\nu$
and twisted fields $\rho^\tau$ from the $\bbZ_{\nu^2+4}$  called $\frak{c}_\alpha$, and amalgamate the orbifold fields $\sigma_\alpha$ from the $\bbZ_{\nu^2+4}$ and the twisted fields $\rho^\tau$ from the
$\bbZ_\nu\times\bbZ_\nu$ called $\frak{d}_a$,  to obtain the elegant fusion rules (first announced in Oberwolfach  in a March 2015 talk by the first author) (we write $\frak{c}_0 = \mathbf{0}+  \frak{b}$ and $\frak{d}_0=\mathbf{0}- \frak{b}$):
\begin{eqnarray}\frak{b}^2 &= &\mathbf{0}+\frak{b}  + \sum \frak{c}_\alpha  + \sum \frak{d}_a =:R\,,\qquad
\frak{c}_\alpha  \frak{d}_b = R - \mathbf{0} =: R_{-}\,,\nonumber\\
\frak{c}_\alpha  \frak{c}_\beta &=& R_{-} + \frak{c}_{\alpha + \beta } + \frak{c}_{\alpha- \beta }\,, \qquad
\frak{d}_a \frak{d}_b = R_{-} - \frak{d}_{a+b} - \frak{d}_{a-b}\,, \nonumber\\
\frak{b} \frak{c}_\alpha &= &R_{-} + \frak{c}_\alpha\,,  \qquad \frak{b} \frak{d}_a = R_{-} - \frak{d}_a\,.\nonumber\end{eqnarray}
Instead of calling $\frak{c}_\alpha$ and $\frak{d}_a$ amalgamations, perhaps we should say the twisted fields of $\bbZ_\nu\times\bbZ_\nu$ and $\bbZ_{\nu^2+4}$ are dropped as being non-local.

\section{Loop groups}

The final important example for us is the loop group $LG$ at level $k\in\bbZ_{>0}$,
where $G$ is a compact, connected, simply-connected Lie group. We write Fus$_k(G)$
for its fusion ring and $P_+^k(G)$ for its primaries $\Phi$. For $G$ of rank $r$, we can
identify $\lambda\in P_+^k(G)$ with the affine highest-weight $\lambda=(\lambda_0;
\lambda_1,\ldots,\lambda_r)\in\bbZ_{\ge 0}^{r+1}$ where $\sum a^\vee_i\lambda_i=k$
for certain positive integers $a_i^\vee$ (the \textit{co-labels}) depending only on $G$. For all $G$,
$a_0^\vee=1$ and $(k;0,\ldots,0)$ denotes the fusion unit $\mathbf{0}$.
Fus$_k(G)$ can be expressed as $R_G/I_k(G)$ for some
ideal $I_k(G)$ of the representation ring $R_G$ called the \textit{fusion ideal}, and primary $\lambda$
is associated to class $[\rho_{\overline{\lambda}}]\in R_G/I_k(G)$ where
$\rho_{\overline{\lambda}}$ is the $G$-irrep with highest-weight 
$\overline{\lambda}=(\lambda_1,\ldots,\lambda_r)$. For example,
for $G=\mathrm{SU}(n)$, we have all $a_i^\vee=1$,
there are $n$ simple-currents (namely the $\lambda$ with some component $\lambda_i=k$), and charge-conjugation $\lambda\mapsto\lambda^*$
is nontrivial iff $n>2$.

The Dixmier-Douady bundles were constructed in \cite{EG1}. Let us review that construction. Recall the construction of bundles for tori in section 4.6. 
Let $G$ be a compact semi-simple Lie group of rank $r$ of 
A-D-E type, e.g.,
$G={\rm SU}(r+1)$. Fix a maximal torus $T\cong \bbR^r/Q^\vee$ of $G$
($Q^\vee$ is the coroot lattice). The orbits of $G$ acting adjointly on
itself are the conjugacy classes of $G$, which are
 parametrised by the Stiefel diagram, which is an affine Weyl chamber.
 More precisely, remove from the Cartan subalgebra
$\bbR^r=\bbR\otimes_\bbZ Q^\vee$ the hyperplanes fixed by a
Weyl reflection $r_\alpha$, as well as the translates of
those hyperplanes by elements of
$Q^\vee$. The Stiefel diagram $S$ is the closure of any
connected component. Any orbit of the adjoint action
intersects $S$ in one and only one point. Points in the
interior of $S$ correspond to generic (so-called regular) elements
of $G$ and have stabiliser $T$, but points on the boundary
have larger stabiliser. 

Consider $G={\rm SU}(2)$ for concreteness. We can 
identify its maximal torus with the circle $\bbR/Q$ where
$Q=\sqrt{2}\bbZ$ is the (co)root lattice and a Stiefel diagram $S$
with half a fundamental domain of $T$, i.e. $0\le x\le 1/\sqrt{2}$. The Hilbert space
is $\cH=L^2(G)$, and the fibres will be the compacts $\cK(\cH)$. 
We want to associate a unitary $U_\gamma$
to any weight $\gamma \in Q^*=\bbZ\Lambda_1$ ($\Lambda_1$ is the fundamental weight).
For any subrepresentation
$\sigma$ in $L^2(G)$, define ``$\sigma\otimes \gamma$'' as follows:
restrict $\sigma$ to $T$ (i.e., write its weight-space
decomposition), and in the Weyl-image $wS\subset T$ act
like the character $\chi_{w\gamma}({\rm e}^{2\pi \i t})
={\rm e}^{2\pi \i \,w\gamma(t)}$. Apply this to the regular representation $\pi$ in $\cH$. 
Then thanks to infinite dimensionality, $\pi\otimes \gamma\simeq\pi$ as
both a representation of $T$ and the Weyl group, so let
$U_\gamma$ be the unitary \hbox{defining} that equivalence. We can
cover $G\simeq S^3$ with two $G$-equivariant patches: $D_1$
about the scalar matrix $I$ and $D_2$ about the scalar
matrix $-I$. The bundle for $G$ on $G$ with level $k\in\bbZ$
 is defined by the following
$G$-equivariant gluing condition: for $x\in T\cap D_1\cap
D_2$ identify $(gxg^{-1},c)$ in $D_1$ with
$(gxg^{-1},{\rm Ad}(\pi_g U_\kappa\pi_g^{-1})c)$ in $D_2$ for any
$g\in G$, $c\in\mathcal{K}$, where $\kappa=k+2$ and we abbreviate  $U_{\kappa\Lambda_1}$ 
to $U_{\kappa}$. We call $U_\kappa$  the twisting
unitary. The consistency condition for these bundles is then that
when $gxg^{-1}=x$, then ${\rm Ad}(\pi_gU_\kappa\pi_g^{-1}U_\kappa^{-1})$
should be the identity, i.e., $\pi_gU_\kappa\pi_g^{-1} =\lambda_g
U_\kappa$ for some character (i.e., one-dimensional
representation) $g\mapsto \lambda_g$ of the stabiliser
$C_G(x)$.

The construction for general $G$ is similar; there $\kappa=k+h^\vee$ where $h^\vee=\sum_{i=0}^ra_i^\vee$ is the dual Coxeter number.
Denote this bundle symbolically by  
$G/\!/_\kappa G$ or just  $\cA_\kappa$ when $G$ is understood. Here, the group $G$  acts on the base $G$ by conjugation,
and $\kappa\in \bbZ\cong H^3_G(G;\bbZ)$ is the twist. The $K$-group ${}^\kappa K_0^G(G)$ is the $K$-theory of the C$^*$-algebra of $\cK$-valued sections of $\cA_\kappa$. The theorem of Freed-Hopkins-Teleman (Theorem 1 of \cite{FHTlg}) is that ${}^\kappa K_0^G(G)\cong {}^\kappa K^{\mathrm{dim}\,G}_G(G)$ is isomorphic as a ring to the fusion ring Fus$_k(G)$.

We expect that (dual) Dirac operators can be used to construct the matrix units for  this 
 bundle $G/\!/_\kappa G$.
More precisely,
\cite{FHTi} explain how to obtain elements of twisted equivariant $K$-theory
using families of Dirac operators. Section 3 of \cite{FHTii} describes how to do this
for ${}^\kappa K^\star_G(G)$. They obtain (for each integrable representation of the centrally extended loop group  $(LG)^\kappa$) a family of Dirac operators parametrised by the
affine space of connections (this has underlying vector space given by the loop algebra $L\mathfrak{g}$).
The connections correspond bijectively to the splittings of $0\rightarrow
L\mathfrak{g}\rightarrow\hat{L}\mathfrak{g}\rightarrow \i\bbR_{rot}\rightarrow 0$
or the linear splittings of the surjection $(L\mathfrak{g})^\kappa\rightarrow L\mathfrak{g}$.
This family of Dirac operators is equivariant with respect to 
$(LG)^\kappa$. This gives an element in the $(LG)^\kappa$-equivariant $K$-theory on the
space of connections, or equivalently on the $\kappa$-twisted $G$-equivariant $K$-theory on $G$. In \cite{FrTe},   the fusion category of the integrable lowest-weight representations of the loop group LG of a compact Lie group G, regarded only as a linear category,
is identified with the twisted, conjugation-equivariant curved Fredholm complexes on the group G: namely, the twisted, equivariant matrix factorizations of a super-potential built from the loop rotation action on LG. 

These Dirac families define elements in the  $KK$-group $KK(\cA_\kappa,\mathrm{pt})$, which can be identified
with $^\kappa K_G^\star(G)$.
We would like to construct (in the spirit of Kasparov \cite{Kas}) the dual Dirac family. This will lie in
$ KK(\mathrm{pt},\cA_\kappa )={}^\kappa K_\star^G(G)$. We expect that
the composition of these will be a matrix unit in $^\kappa KK_G(G,G)$, while the composition
in the other direction will be $\delta_{\lambda,\mu}\in KK(\mathrm{pt},\mathrm{pt})=\bbZ\oplus 0$.
We plan to return to this in a future publication.

The generic module categories for the loop group modular tensor categories $\cC(G,k)$ are now known \cite{Ganlb} to be built in standard ways from symmetries of the extended Dynkin diagram of $G$: for any $G$, there is a bound $K$ (growing cubically with the rank) such that any module category for $\cC(G,k)$, when $k>K$, is generic in that sense. In the remainder of this section we give the $K$-theoretic description of these generic module categories. It would be very interesting to do the same for the known exceptional module categories --- see section 6 in both \cite{EG1} and  \cite{EG3} for work in this direction.

\subsection{Outer automorphism modular invariants}

The most obvious modular invariant $\cZ$, other than the identity $I$, is charge-conjugation $\cZ=S^2$.
In the loop group $LG$ setting, charge-conjugation corresponds to an outer
automorphism of $G$. More generally, see \cite{EG3} where the $K$-theoretic interpretation of the nimrep,
alpha-induction etc for the modular invariants  associated to any
outer automorphism of $G$ is developed and given in more detail.

Let $G$ be as above and write $\kappa=k+h^\vee$ for some  level $k\in\bbZ_{\ge 0}$ as before. The  group of outer automorphisms of $G$
is naturally identified with the group of symmetries of the Dynkin diagram of $G$, and as
a permutation of these vertices also permutes highest weights of $G$ in the usual 
way. More precisely, pick an outer automorphism $\omega$ of $G$, and define $\omega(\la)$ to be the weight whose $i$th component is $\lambda_{\omega{i}}$  ($\omega$
fixes $\lambda_0$).  Through this, $\omega$ permutes the level $k$ primaries $\lambda\in P^k_+(G)$.
For example, an automorphism of $G=\mathrm{SU}(n)$ realising the
charge-conjugation automorphism is taking the complex conjugation of the unitary representations.

The modular invariant $\cZ^\omega$ corresponds to the $KK$-element coming from the 
bundle map $\cA_\kappa\rightarrow \cA_\kappa$
associated to the automorphism $\omega:G\rightarrow G$. The full system is the fusion ring, and alpha-induction is
$\alpha_+=\omega$ and $\alpha_-=\mathrm{id}$.

On the other hand,   inner automorphisms are invisible as they act trivially on $P^k_+(G)$, so from this point of view can be ignored.

\subsection{Simple-current modular invariants}

The other generic source of modular invariants for the loop groups are the simple-current
modular invariants.
These correspond to strings living on non-simply-connected groups
$G/Z$ where $Z$ is some subgroup of the centre $Z(G)$ of $G$.

Recall the discussion of simple-current from section 3.1.
The group of simple-currents for $G\times H$
is the direct product of those for $G$ and for  $H$, so it suffices to
consider simple $G$. All of these simple-currents correspond to extended
Dynkin diagram symmetries, with the single exception Fus$_2(E_8)$ which we will ignore 
 as  it does not yield a modular
 invariant (module category) for $E_8$. For any $G$ and $k$, the simple-currents and outer
automorphisms together generate all symmetries of the extended Dynkin diagram.

For example, the group of
simple-currents for Fus$_k(\mathrm{SU}(n))$ is cyclic of order $n$, generated by
$J=(0;k,0,\ldots,0)$ which permutes $P_+^k(\mathrm{SU}(n))$ through 
 $(\lambda_0;\lambda_1,\ldots,\lambda_{n-1})\mapsto(\lambda_{n-1};\lambda_0,
\lambda_1,\ldots,\lambda_{n-2})$. Then the grading is $\cQ_\lambda(J^d)=\xi_n^{d\sum_{i=1}^{n-1}
i\lambda_i}$ and $T_{J^d,J^d}\overline{T_{\mathbf{0},\mathbf{0}}}=\xi^{kd(n-d)}_{2n}$, where we write $\xi_n=e^{2\pi i/n}$. 

The $K$-theoretic treatment of simple currents  was developed
in \cite{EG3}. Restrict to $G_n={\rm SU}(n)$, the most interesting case.
Fix any divisor $d|n$ and level $k$, and write $n'=n/d$ and $\kappa=k+n$. 
The simple-current modular invariant $\cZ_{\langle J^{n'}\rangle}$ exists (i.e. there are no quaternions in $\langle J^{n'}\rangle$) iff
$n'(n+1)k$ is even. Write $Z_d$ for the order-$d$ subgroup of the centre   
 $Z(G_n)\cong \bbZ_n$. The corresponding nimrep has Grothendieck group
${}^\tau K^{G_n}_0(G_n/Z_d)\simeq {}^\tau K^{G_n\times Z_d}(G_n)$ for some twist
$\tau$, with the module structure coming from the pushforward of the obvious multiplication $G_n\times (G_n/Z_d)\to G_n/Z_d$ of the bases. The type 1 parents have fusion ring ${}^{\tau'}K^{G_n/Z_{d'}}_0(G_n/Z_{d'})$ for some
subgroup $Z_{d'}$ of $Z_d$ and twist $\tau'$. The full system should be
$ R_{Z_{d'}}
\otimes_\bbZ{}^{\tau''}K_0^{G_n^{{adj}}\times Z_{d'}^{L}}(G_n)$, where the product is component-wise, with the second component product coming from the pushforward of multiplication on $G_n$. The appropriate twists
$\tau,\tau',\tau''$ are given in section 5 of \cite{EG3}. These $K$-homology groups all vanish in degree 1.

Let $G$ be any simply-connected,
connected, compact Lie group $G$.
For such $G$, multiplication by the centre $Z(G)$ should correspond naturally 
to the action of the simple-currents in the fusion ring $^\kappa 
K_{G}^{\mathrm{dim}\,G}(G)$, in the following sense. The primaries
$\lambda\in P_+^k(G)$ are identified with certain conjugacy classes ---
this yields a geometric picture of Fus$_k(G)$ dual to
the usual representation ring description $R_G/I_k$.
 Now, $Z(G)$ permutes these conjugacy classes by multiplication,
and this permutation agrees with the simple-current action on primaries. 
For example, for $G=\mathrm{SU}(2)$, the level $k$ primaries are
$\lambda=(k-\lambda_1;\lambda_1)$ for integers $0\le \lambda_1\le k$; this 
corresponds to the conjugacy class intersecting the Stiefel diagram $S$ (half of
the maximal torus $\bbT$) at
diag$(\exp[2\pi\i(\lambda_1+1)/2\kappa],\exp[-2\pi\i(\lambda_1+1)/2\kappa])$ 
for $\kappa=k+2$.
By multiplication, the central element $z=-I$ sends the $\lambda_1$ conjugacy
class  to the $\lambda_1+\kappa$ one,
which the Weyl group identifies with $\kappa-\lambda_1-2$. This matches the 
action of the simple-current.

The relation between the centre of $G$ and the simple-currents is part of CFT folklore.
The generalisation of this calculation to all simple simply connected connected compact $G$ is central to
what follows. To each element $z$ of the centre of $G$, we obtain an invertible bundle map
$\cA_\kappa\rightarrow \cA_\kappa$ given by multiplication of the base space  by $z$, and hence
an invertible element $J_z\in {}^\kappa KK_G(G,G)$.  \cite{EG3} conjectured that
this should correspond directly to
the multiplication of the simple-current in the fusion ring. We can now prove it:

\medskip\noindent\textbf{Proposition 2.} \textit{For each central element $z$ of $G$ and each level $k\in\bbZ_{>0}$, there exists a
simple-current $j_z\in {}^{k+h^\vee}K_G^\star(G)$ such that the map $J_z:{}^{k+h^\vee}K_G^\star(G)\rightarrow
{}^{k+h^\vee}K_G^\star(G)$ corresponds to the fusion product by $j_z$. Moreover, provided $(G,k)\ne (E_8,2)$,
this map $z\mapsto j_z$ is a group isomorphism.}\medskip

The easiest way to reduce this to a familiar calculation from CFT, is perhaps Theorem 4.13 of
\cite{Mein}, which builds ${}^{k+h^\vee}K_0^G(G)$ up from the  conjugacy classes {$cl_\xi$} of prequantised elements
$\exp(\xi)$
at level $k$. If we let $\mu=B^\flat(k\xi)$ denote the associated weight, then the push-forward of the
inclusion $cl_\xi\hookrightarrow G$ sends the fundamental class $[cl_\xi]\in K_0^G(cl_\xi,Cl(Tcl_\xi))$
(where $Cl(T\cC_\xi)$ is the Clifford bundle on the conjugacy class) to the element $[\chi_\mu]\in R_G/I_k
={}^{k+h^\vee}K_0^G(G)$. If we restrict to $\xi$ from the level $k$ alcove, this gives a  natural basis for
${}^{k+h^\vee}K_0^G(G)$. Then
$J_z$ acts on conjugacy classes by sending $cl_{\exp(\xi)}$ to $cl_{z\exp(\xi)}=zcl_{\exp(\xi)}$. Then  $z\exp(\xi)$ is also prequantised at level $k$, and it is known that its weight is $j_z\mu$ for some simple-current $j_z$. 
This defines an injective group homomorphism from $Z(G)$ to the level $k$ simple-currents of $G$.
The latter was computed in \cite{Fuchs}, and it is found that except for $E_8$ at level 2, the group of
simple-currents is the same (finite) size as that of $Z(G)$. The obvious analogue of Proposition 2 holds when
$G$ is no longer simple.

Fix any subgroup $Z$ of the centre $Z(G)$ and write $\overline{G}$ for $G/Z$.
The embedding $Z\hookrightarrow G$ of groups and of the fixed-point pt$\,=1$ into the space $G$,
 yields the $K$-theory map (restriction) $^\kappa K_G^\star(G)\rightarrow
K_Z^\star(\mathrm{pt})=R_{\widehat{Z}}$. This is how we recover here the fact (recall section 3.1) that  
the fusion ring Fus$_k(G)$ carries a grading $a\mapsto\cQ_a$ by
representations of $Z$:  $\mathrm{Fus}_k(G)=\oplus_{\psi\in\widehat{Z}}
\mathrm{Fus}_k(G)^\psi$.

To understand the modular invariants which aren't type 1, note that 
$\overline{G}$ has a well-defined adjoint action on $G$ since the centre of $G$ on $G$ acts trivially. Now, $H^3_{\overline{G}}(G;\bbZ)\cong\bbZ\oplus R_Z$, where $\bbZ=H^3(G;\bbZ)$
and $R_Z=H_{\overline{G}}^3(\mathrm{pt};\bbZ)$, as is
easily verified using the spectral sequence associated to the fibration $G\rightarrow(E_{\overline{G}}\times
G)/\overline{G}\rightarrow B_{\overline{G}}$. To help identify levels etc, we find that $H^3_{
\overline{G}}(\overline{G};\bbZ)=H^3(\overline{G};\bbZ)$, transgression $H^4_{\overline{G}}
(\mathrm{pt};\bbZ)\rightarrow H^3(\overline{G};\bbZ)$ is $\times n$ for SU$(n)$,
and $H^3(\overline{G};\bbZ)\rightarrow H^3(G;\bbZ)$ is $\times 2$ or $\times 1$
depending on $n$ and $d$.

Now,  any representation $\lambda$ of $G$ projects to a projective
representation of $\overline{G}=G/Z$, with multiplier $\psi$ given by the restriction $\cQ_\lambda|_Z$ of the representation
to  the abelian group $Z$ (more precisely, a $G$-irrep $\rho$ restricts to $\rho|_Z$, which
will consist of dim$\,\rho$ copies of the same irrep $\psi_\rho\in\widehat{Z}$).  Indeed, $H^2_{\overline{G}}(\mathrm{pt};\bbT)\cong R_Z$
parametrises the projective equivalence classes of projective representations of $\overline{G}$.
Mimicking the construction of $G/\!/_\kappa G$, 
we can construct a bundle $G/\!/_\kappa\overline{G}$ with components (subbundles) for each $\psi\in\widehat{Z}$.
We denote these subbundles by $G/\!/_{(\kappa,\psi)}\overline{G}$, or simply  $\cA_\kappa^\psi$.

We have proved:

\medskip\noindent\textbf{Proposition 3.} \textit{Let $G$ be compact, connected and simply-connected, of dimension $d$. Fix any subgroup $Z$ of the centre $Z(G)$ and nonzero level
$k$, and write $\kappa=k+h^\vee$. Then the fusion ring $^\kappa K^d_G(G)=\mathrm{Fus}_k(G)$ is graded by restriction to
$Z$: $\mathrm{Fus}_k(G)=\oplus_{\psi\in\widehat{Z}}
\mathrm{Fus}_k(G)^\psi$.   For each level $k$ and $\psi\in\widehat{Z}$,
${}^{(\kappa,\psi)}K^d_{G/Z}(G)\cong \mathrm{Fus}_k(G)^\psi$ (as an abelian group) and ${}^{(\kappa,\psi)}K^{d+1}_{G/Z}(G)=0$.}\medskip

We get something analogous when $G/Z$ acts on $G/Z_0$ for $Z_0\le Z\le Z(G)$  --- in fact this is used at the end of this subsection.

Let's turn now to the modular invariants. Let $G$ be as above.
Recall the parametrisation of simple-current modular invariants in Theorem 1, given by a subgroup $J$ and
group homomorphism $\epsilon:J\rightarrow \widehat{J}$. The   subgroup
$J$ there is our subgroup $Z\le Z(G)$. More precisely, we have the isomorphism $z\mapsto j_z$ of Proposition 2, from the centre $Z(G)$ to the group of simple-currents. For convenience, we'll write $\epsilon_z$ and $q(z)$ in place of $\epsilon_{j_z}$ and $T_{j_z,j_z}\overline{T_{\mathbf{0},\mathbf{0}}}$. When $Z$ is cyclic (which is automatic for all simple $G$ except Spin$(4n)$), $\psi\equiv 1$ and $\epsilon$ is completely determined by the value $q(z)$ for a generator $z\in Z$.

 Consider first a type 1 module category  $(Z,\psi)$ (i.e.\ one of extension type). This requires $q(z)=1$ for all $z\in Z$, and $\psi\equiv 1$. Write $\overline{G}
=G/Z$. As before, we can define a bundle $\overline{G}/\!/_{\kappa'}\overline{G}$ (as usual throughout this subsection we largely suppress the subtleties of the twists, but these are addressed in section 5.2 of \cite{EG3}), and as before this falls into subbundles $\overline{\cA}_{\tau'}^\psi$ parametrised by $\psi\in\widehat{Z}$.  Sigma-restriction is the obvious projection $\pi:\cA_\kappa\to\overline{\cA}_{\kappa'}^1$, which kills any $\lambda\in P_+^k(G)$ unless $\cQ_\lambda|_Z=1$, in which case it sends $\lambda$ to its $Z$-orbit.
This type 1 modular invariant is the correspondence
$$ \begin{tikzpicture}
  \matrix (m) [matrix of math nodes,row sep=2em,column sep=2em,minimum width=2em]
  {
    \cA_\kappa & &\cA_\kappa\\ & \overline{\cA}_{\kappa'} &  \\ };
  \path[-stealth]
    (m-1-1) edge node [left] {$\pi\ $} (m-2-2)
         (m-1-3)   edge node [right] {$\ \pi$} (m-2-2);
\end{tikzpicture}$$
Writing $H=\Delta_G(1\times Z)$, we can pull this back to a  bundle of type $H/\!/H$, formally reminiscent of \eqref{modinvgen1}.

Consider next the  type 2 module categories $(Z,\psi)$ (i.e. those of pure automorphism type). This requires $\epsilon:Z\to\widehat{Z}$ to be an isomorphism; in the case where $Z=\langle z\rangle$ is cyclic, this is equivalent to demanding that $q(z)$ is a primitive $|Z|$-root of unity. The associated $KK$-element is then
\begin{equation}\label{scinvII}\begin{tikzpicture}
  \matrix (m) [matrix of math nodes,row sep=2em,column sep=2em,minimum width=2em]
  {
     \cA_\kappa & & && \cA_\kappa\\
   &  \oplus_{z\in Z}\cA_\kappa^{\epsilon_z}&&\oplus_{z\in Z}\cA_\kappa^{\epsilon_z}&  \\ };
  \path[-stealth]
    (m-1-1) edge node [left] {$\pi$\ } (m-2-2)
    (m-1-5) edge node [right] {\ $\pi$} (m-2-4) 
    (m-2-2) edge node [above] {$\times z$} (m-2-4);
\end{tikzpicture}\end{equation}
where the diagonal maps act on morphisms by projecting $G\rightarrow G/Z$, and $\pi_*$ lifts  projective
 $\overline{G}$-representations to ordinary  $G$-representations. Explicitly, $\pi$ puts $\lambda\in P_+^k(G)$ into the component $\cA_\tau^\phi$ where $\phi=\cQ_\lambda$; unlike the type 1 $\pi$, its kernel is trivial. The horizontal map is multiplication of the subbundle $\cA_\kappa^{\epsilon_z}$ by $z$ (recall Proposition 2). This recovers the modular invariant which Theorem 1(a) associates to $(Z,\psi)$.

Using the preceding two paragraphs, it now is easy to find the $KK$-element corresponding to any $(Z,\psi)$. Write $Z_R$ for the kernel of $\epsilon$, and $Z_L$ for the kernel of the transpose of $\epsilon$, as in the proof of Theorem 1(a). Then $\epsilon$ yields an isomorphism
$\widetilde{\epsilon}:Z/Z_R\rightarrow \widehat{Z/Z_L}$. The two type 1 parents are $\overline{\cA}_{\kappa_L}^1:=(G/Z_L)/\!/_{(\kappa_L,1)}(G/Z_L)$ and $\overline{\cA}_{\kappa_R}^1:=(G/Z_R)/\!/_{(\kappa_R,1)}(G/Z_R)$, and they are linked by $\widetilde{\epsilon}$:
\begin{equation}\label{scinvII}\begin{tikzpicture}
  \matrix (m) [matrix of math nodes,row sep=2em,column sep=2em,minimum width=2em]
  {
     \cA_\kappa & & && &&\cA_\kappa\\  & \overline{\cA}_{\kappa_L}^1& & &&\overline{\cA}_{\kappa_L}^1&\\
   &&  \oplus_{z\in Z/Z_L}\widetilde{\cA}_{\kappa_L}^{\tilde{\epsilon}_z}&&\oplus_{z\in Z/Z_L}\widetilde{\cA}_{\kappa_R}^{\tilde{\epsilon}_z}&&  \\ };
  \path[-stealth]  (m-1-1) edge node [left] {$\pi_L$\ } (m-2-2)
    (m-1-7) edge node [right] {\ $\pi_R$} (m-2-6) 
    (m-2-2) edge node [left] {$\pi_L'$\ } (m-3-3)
    (m-2-6) edge node [right] {\ $\pi'_R$} (m-3-5) 
    (m-3-3) edge node [above] {$\times z$} (m-3-5);
\end{tikzpicture}\end{equation}
where we write $\widetilde{\cA}_{\kappa_L}^\phi$ and $\widetilde{\cA}_{\kappa_R}^\phi$ for the obvious subbundles of $(G/Z)/\!/_{\kappa_L}(G/Z)$ (where $\phi\in\widehat{Z/Z_L}$) and $(G/Z)/\!/_{\kappa_R}(G/Z)$ (where $\phi\in\widehat{Z/Z_R}$).
In particular, $\pi_L$ kills $\lambda\in P_+^k(G)$ unless $\cQ_\lambda|_{Z_L}=1$, in which case $\pi_L'$ sends the $Z_L$-orbit of $\lambda$ to $\widetilde{\cA}_{\kappa_L}^\phi$ with $\phi=\cQ_\lambda|_Z\in\widehat{Z/Z_L}$. The horizontal map is as before.

Combining simple-current modular invariants with outer automorphisms is trivial here: compose (multiply) the correspondence for $(Z,\psi)$ with that for $\omega$.

\section{Modular data reinterpreted}

It would be highly desirable to interpret $K$-theoretically the modular group representation $\rho$ associated
to these fusion rings $^\tau K_G^\star(G)$. For instance, this would allow us in principle to 
understand \textit{intrinsically} (i.e. in the $KK$- or $K$-world) 
why the modular invariants are invariant under the action of
SL$_2(\bbZ)$.
 It suffices to focus on the generators
$S:=\rho\left({0\atop 1}{-1\atop 0}\right)$ and $T:=\rho\left({1\atop 0}{1\atop 1}\right)$. 
However, these matrices $S,T$ in the Verlinde (i.e. primary) basis are complex --- more precisely cyclotomic
--- and never integral. So that means we need to realise them
in the complexification $\bbC\otimes_\bbZ {}^\tau KK_G(G,G)$. This suggests to use
(twisted equivariant) Chern characters, which is a ring isomorphism from complexified
$K$-groups to cohomology rings. 

In this section we consider the Chern character maps. At least for the cases in which we are interested, the target
of these maps (some cohomology space) comes with a preferred basis, and in terms of this basis the Verlinde (fusion)
product is diagonal. This change-of-basis may make the matrices $S,T$ look messier.

\subsection{Finite groups and cyclic homology}

As warm-up, let's begin with the Chern character for the $G/\!/G$ groupoid, where $G$ is finite.

The primaries for $G/\!/G$ are pairs $[g,\phi]$, where $g$ labels conjugacy classes and $\phi\in \mathrm{Irr}(C_G(g))$, and form a basis for the complexified fusion ring $\bbC\otimes_\bbZ\mathrm{Fus}$. Another basis, called the \textit{monomial basis}, is parametrised by equivalence classes $[(g,h)]$ of commuting pairs $g,h\in G$, up to simultaneous conjugation $(g,h)\sim(g^k,h^k)$. In terms of the monomial basis, the SL$_2(\bbZ)$-action is especially simple:
\begin{equation}\left({a\atop c}{b\atop d}\right).[(g,h)]=[(g^ah^b,g^ch^d)]\,.\label{monbas}\end{equation}

The change-of-basis between the primary basis $\{[g,\phi]\}$ and the monomial basis
$\{[(g,h)]\}$,  involves the character tables of the centralisers $C_G(g)$. 
We will show here that this change of basis can be interpreted naturally using equivariant Chern 
characters. Our treatment of Chern characters follows \cite{BauCon}, which describes these for $G$  finite. We will
specialise their discussion to the case where $G$ acts on itself by conjugation ({which we write as a right-action}, to match \cite{BauCon}).

Write $cp({G})=\{(g,\gamma)\in G\times G\,:\,g\gamma=\gamma g\}$ for the set of all commuting
pairs. $G$ acts on $cp({G})$ by simultaneous conjugation: $(g,\gamma)\kappa=(\kappa^{-1}g
\kappa,\kappa^{-1}\gamma\kappa)$. Let $G^\gamma=C_G(\gamma)$. Then
$cp({G})$ can be regarded as the disjoint union $\sqcup_{\gamma\in G} G^\gamma$, where $G^\gamma$ denotes
all $g\in G$ fixed by this action of $\gamma$. The target for the Chern map is
\begin{equation}\nonumber 
HP^0(G/\!/G)\oplus HP^1(G/\!/G)=H_c^\star(cp({G})/G;\bbC)=H^\star_c(cp({G});\bbC)^G=\left(\oplus_{\gamma\in G} H^\star_c(G^\gamma;
\bbC)\right)^G\,,\end{equation}
where $HP^\star$ denotes periodic cyclic cohomology and $H_c$ denotes \v Cech cohomology. Because the base space of $G/\!/G$  is 0-dimensional,
$HP^1(G/\!/G)=H^{od}_c(cp({G});\bbC)^G=0$ while $HP^0(G/\!/G)=H^{ev}_c(cp({G});\bbC)^G$ can be identified with all $\bbC$-valued maps $f$ from the set
of commuting pairs, which are constant on the $G$-orbits: $f(g^\kappa,
\gamma^\kappa)=f(g,\gamma)$.

Take some $G$-equivariant bundle over $G$: $\cE_{(h,\phi)}$ where $h\in G$ is a representative 
of a conjugacy class and $\phi$ is
an irrep of $C_G(h)$. Then $\gamma$ acts on the fibre $(\cE_{(h,\phi)})_{g}$, whenever $g$
commutes with $\gamma$. This fibre is 0-dimensional unless $g$ is conjugate to $h$,
in which case the fibre carries an action of $C_G(g)\cong C_G(h)$ equivalent to $\phi$.
Diagonalise that action by $\gamma$: $\phi(\gamma)\cong V_1\oplus\cdots \oplus V_r$ where
$V_i$ is the eigenspace of dimension $d_i$, with eigenvalue $\lambda_i$. Define
\begin{equation}
\mathrm{ch}_G^\gamma(\cE_{(h,\phi)})=\sum_{i=1}^r\lambda_i\mathrm{ch}(E^i)\,,\end{equation}
where ch$(E^i)$ is the usual (non-equivariant) Chern character of the $d_i$-dimensional
trivial bundle over $C_G(\gamma)$. We can think of ch$(E^i)$ as $d_ie_{[g,\gamma]_G}$ where $[g,\gamma]_G$ denotes the $G$-orbit (acting by simultaneous conjugation),
and $e_{S}$, for $S$ a set, is the characteristic function of $S$ written formally as $e_S=\sum_{s\in S}
e^s$.  Packaging these together,
\begin{equation}
\mathrm{ch}_G(\cE_{(h,\phi)})=\oplus_{\gamma\in G}\mathrm{ch}_G^\gamma(\cE_{(h,\phi)})\,.
\end{equation}

So for finite groups at least, the Chern character maps to functions on commuting pairs.
The ring homomorphism property tells us that the fusion ring structure is the usual
tube algebra one (cf chapter 12 of \cite{EKaw}), so in terms of the cohomology basis, the fusion product has been diagonalised.
Theorem 1.19 of \cite{BauCon} tell us ch$_G:K^0_G(G)\otimes_\bbZ\bbC\rightarrow HP^0(G/\!/G)$ is an
isomorphism of $\bbC$-vector spaces. (Of course in odd degree, both vanish.) 
Any field in which the relevant eigenvalues live, would have worked in place of $\bbC$,
so we could have used e.g. the cyclotomic field $\bbQ[\xi_{|G|}]$ in place of $\bbC$,
if we had wanted. 

Our point is that  not only the fusion product 
has a simplified interpretation in the cohomology basis. We can think of
the SL$_2(\bbZ)$-representation as living in bivariant cohomology $HP^\star(G/\!/G,G/\!/G)$.
The homology group of the torus is of course $H_2(T^2;\bbZ)\cong\bbZ^2$. The set of all group homomorphisms
$f:H_2(T^2;\bbZ)\rightarrow G$ are in natural bijection with pairs $(g,h)\in cp(G)$, namely $f(1,0)=g$ and $f(0,1)=h$.
The group $G$ acts naturally on itself by conjugation, and hence likewise on these maps $f(m,n)$ by conjugation. This recovers the $G$-action on $cp(G)$ given above.   
The modular group SL$_2(\bbZ)$ of the torus acts naturally on homology $H_2(T^2;\bbZ)$ (as change of basis), and hence 
likewise  on the $f(m,n)$, and hence $HP^\star(G/\!/G)$.

This SL$_2(\bbZ)$-action can be worked out explicitly.  For each $\mu=\left({a\atop c}{b\atop d}\right)\in\mathrm{SL}_2(\bbZ)$, construct a correspondence  
\begin{equation}\begin{tikzpicture}
  \matrix (m) [matrix of math nodes,row sep=2em,column sep=2em,minimum width=2em]
  {
     G/\!/G& &{G}/\!/G\\
&    ( G\times G/\!/G\times G, P_\mu ) & \\ };
  \path[-stealth]
    (m-1-1) edge node [left] {$\iota^L\ \ $} (m-2-2)
    (m-1-3) edge node [right] {$\ \iota^R$} (m-2-2);
\end{tikzpicture}\end{equation}
where $P_\mu$ is the element of $HP^0(G\times G/\!/G\times G)$ with fibre 0 everywhere except for $\bbC$ above the point $(g,\gamma;g^{a}\gamma^c,g^b\gamma^d)\in cp(G)^2$. 
Then the  combination of pullback,  multiplying by $P_\mu$ and   pushforward is the map sending $(g,\gamma)$ to $(g^{a}\gamma^c,g^b\gamma^d)$.

The Chern character in the presence of a 3-cocycle $\omega$ is not much different. The target will consist of maps $cp(G)\rightarrow
\bbC$ which are covariant with respect to a $G$-action which twists $\bbC$ by a factor coming from $\omega$. The resulting SL$_2(\bbZ)$-action is \eqref{monbas} with phases thrown in.

\subsection{The torus}

Determining the Chern character for twisted equivariant $K$-theory such as $^\tau K_G^\star(G)$ for $G$ compact is the main  purpose of \cite{FHT}.
We specialise to the case of the torus here and give the more general result next subsection.

Recall a transgressed twist of $T$ acting trivially on itself is given by an even lattice $L\subset \bbR^d$. The Chern character of $^{\tau_L}
K^\star_T(T)$ should localise to conjugacy classes, i.e. elements $\mu+L\in T$ for all $\mu\in L^*$, where $T$ is identified with
$\bbR^d/L$. We obtain $^{\tau_L}K^\star_T(T)\otimes_\bbZ \bbC=\oplus_{[\mu]\in L^*/L}{}^{\tau_L}H^\star_T(T)  $, where
${}^{\tau_L}H^\star_T(T)  $ is $\bbC$ or 0 depending on whether or not $\star=d$. Thus we can regard the vector space
$^{\tau_L}K^\star_T(T)\otimes_\bbZ \bbC$ as having a `Chern' basis parametrised by classes $[\mu]\in L^*/L$, which we
can embed into the $\bbC^d$ algebra as $[\mu]\mapsto(e^{2\pi i \langle \mu,\lambda_1\rangle},\ldots,e^{2\pi i \langle \mu,\lambda_d\rangle})$ where
$\lambda_j$ is some basis of $L^*$. That this is a ring homomorphism, means the Chern character diagonalises 
the fusion product, as with the last subsection. 

A natural guess is that,
at least for the torus, the SL$_2(\bbZ)$-action on the (complexified) fusion ring is related to Mukai's projective 
SL$_2(\bbZ)$-action on the derived
category for abelian varieties. We sketch how this should go  in the following, which is a twisted version of \cite{GJK}.

Suppose that $G$ is a (commutative)  torus $T$.
Then in the standard (non-equivariant) $T$-duality one takes
\begin{equation}\begin{tikzpicture}
  \matrix (m) [matrix of math nodes,row sep=2em,column sep=2em,minimum width=2em]
  {
     T& &T^*\\
&    T\times T^* & \\ };
  \path[-stealth]
    (m-1-1) edge node [left] {} (m-2-2)
    (m-1-3) edge node [right] {} (m-2-2);
\end{tikzpicture}\end{equation}
where $T^*$ is not the Pontryagin group dual of $T$ {(relevant to the usual Fourier transform)} but the dual torus.
In the equivariant case, one takes
\begin{equation}\begin{tikzpicture}
  \matrix (m) [matrix of math nodes,row sep=2em,column sep=2em,minimum width=2em]
  {
     T {/\!/_LT} & & {T^* /\!/_{L^*} T^*}\\
&    ({(T \times T^*)/\!/_{L\times L^*} ( T \times T^*)} , P)& \\ };
  \path[-stealth]
    (m-1-1) edge node [left] {} (m-2-2)
    (m-1-3) edge node [right] {} (m-2-2);
\end{tikzpicture}\end{equation}
where $P$ is the {\it Poincar\'e bundle}, which can be thought of as the equivariant {$T \times {T}^*$} bundle $P$ on {$T \times T^*$} which lives on $(g,\pi)$ with representation
$(\pi, -g)$. {The (equivariant) \textit{Fourier-Mukai functor} $\mathcal{S}:\cD^b_T(T)\rightarrow \cD^b_{T^*}(T^*)$ is a
triangulated category equivalence between those derived categories, built from $P$. This map is defined over $K$-theory
on $\bbZ$. Fix a non-degenerate line bundle
$L$ over $T$ (this will correspond to a choice of `level'). Then we get a projective SL$_2(\bbZ)$-action on the
complexification of the $K$-group of $\cD^b_T(T)$,
where the modular matrix $T_{mod}$ corresponds to multiplying a $K$-class by $L$ and the modular $S$ matrix
corresponds to $\phi_L^*\circ\cS$, for $\phi_L:T\rightarrow T^*$ an isogeny. These matrices can be computed using the pre-quantised points which carry the $K$-theory, as in the proof of Proposition 2 in section 6.2 above. This action is merely projective because $S_L^4$ acts as degree shifts (rather than the identity); put another
way, we are missing the $e^{-\pi ic/12}$  factor (here a 24th root of 1) in the matrix $T_L$.
 This projectivity} is similar to what happens in a modular tensor category, and isn't significant. {Apart from it, this should match the
 SL$_2(\bbZ)$-representation \eqref{WeilT},\eqref{WeilS} at the appropriate level. This action then descends to the Grothendieck
 group $K(\cD^b_T(T))={}^LK_T^0(T)$.}

\subsection{$G$ compact connected simply-connected}

Consider finally the (most interesting) case of the loop group of compact connected simply-connected
$G$. 
Understanding the (twisted equivariant) Chern character here is the main point of \cite{FHT}.
{There,}  $^\tau K_G^\star(G)\otimes_\bbZ\bbC$ is identified with
$\oplus_g {}^\tau H^\star_{Z(g)}(G^g;{}^\tau \cL(g))$, where the
sum is over conjugacy class representatives modulo the Weyl group. All but finitely many of those twisted
cohomology groups are trivial; the only $g$ that contribute are when $g$ lies in the {
 conjugacy classes $\exp(2\pi \i(\lambda+\rho)/(k+h^\vee))$ as $\lambda$ runs through the level $k$ highest weights $P_+^k$}, which each contribute 1-dimension when $\star$ equals the rank$(G)$. For $SU(2)$ these
conjugacy classes have representatives diag$( \xi_{2(k+2)}^l,\xi_{2(k+2)}^{-l})$
where $l=1,\ldots,k+1$ and $\xi_n=e^{2\pi\i /n}$.
We should think of the direct sum of those twisted cohomology groups as the space spanned
by the characters of $G$ evaluated at those special $g$. That is how
to think of the Chern characters here: as a vector-valued map (one component for
each special conjugacy class  $g$). As it must be, this is a ring homomorphism.
One can think of this as the map associating primary $\lambda$ to
the vector $(S_{\lambda,\mu}/S_{\mathbf{0},\mu})_{\mu\in {P_+^k}}$. Although this is a ring homomorphism,
it is  simpler to drop the denominators, and regard this as a map $\lambda\mapsto
(S_{\lambda,\mu})_{\mu\in P_+}$. In either case, 
SL$_2(\bbZ)$ doesn't act any simpler: $\left({0\atop -1}{1\atop 0}\right)$ again corresponds to $S$
and $\left({1\atop 0}{1\atop 1}\right)$ goes to
$STS^*=T^*S^*T^*$. So the image of  $\left({1\atop 0}{1\atop 1}\right)$ is no longer diagonal. However, the matrix corresponding
to $\left({1\atop 1}{0\atop 1}\right)$ will be diagonal, and correspond to
$T$. So although we haven't simplified the SL$_2(\bbZ)$ action, we are now 
in a complex vector space, and we have diagonalised the fusion products.

The expectation here is that we can derive the SL$_2(\bbZ)$ action from that of the maximal
torus, equivariantised over the (finite) Weyl group $W$. In particular, Atiyah  \cite{Ati}
proved that the 
restriction map from $K_G^\star (X)\rightarrow K_T^\star(X)$  has a natural left inverse, so
$K_G^\star(X)$ is a direct summand of $K_T^\star(X)$, contained in the $W$-invariant part.
This $W$-action on $K_T^\star(X)$ can be realised by $KK$-elements. Presumably this extends to
$Z^3$-twists. We're interested in
$X=G$, and then  we project $K_T^\star(G)\rightarrow {}^LK_T^\star(T)$, where $L=(k+h^\vee)Q^\vee$, through the inclusion of the
space $T$ into $G$. The latter is the fusion algebra of the torus, with a basis parametrised
(if we like) by theta functions (which carry the SL$_2(\bbZ)$ action corresponding
to its modular data), and because of Kac-Peterson we know we can identify the characters of the loop group with alternating sums over $W$ of those theta functions, divided by some
anti-symmetric denominator. It all sounds like we should be able to recover not only
the loop group modular data, but in fact the characters themselves by reducing it
to the maximal torus and (anti)-symmetrising over $W$ in the appropriate way.

\newcommand\biba[7]   {\bibitem{#1} {#2:} {\sl #3.} {\rm #4} {\bf #5,}
                    {#6 } {#7}}
                    \newcommand\bibx[4]   {\bibitem{#1} {#2:} {\sl #3} {\rm #4}}

\def\ASENS            {Ann. Sci. \'Ec. Norm. Sup.}
\def\AM   {Acta Math.}
   \def\AnM              {Ann. Math.}
   \def\CMP              {Commun.\ Math.\ Phys.}
   \def\IJM              {Internat.\ J. Math.}
   \def\JAMS             {J. Amer. Math. Soc.}
\def\JFA              {J.\ Funct.\ Anal.}
\def\JMP              {J.\ Math.\ Phys.}
\def\JRA              {J. Reine Angew. Math.}
\def\JSP              {J.\ Stat.\ Physics}
\def\LMP              {Lett.\ Math.\ Phys.}
\def\RMP              {Rev.\ Math.\ Phys.}
\def\RNM              {Res.\ Notes\ Math.}
\def\RIMS             {Publ.\ RIMS.\ Kyoto Univ.}
\def\Inv              {Invent.\ Math.}
\def\npbp             {Nucl.\ Phys.\ {\bf B} (Proc.\ Suppl.)}
\def\nupb             {Nucl.\ Phys.\ {\bf B}}
\def\nup              {Nucl.\ Phys. }
\def\nupp             {Nucl.\ Phys.\ (Proc.\ Suppl.) }
\def\adma             {Adv.\ Math.}
\def\coma             {Con\-temp.\ Math.}
\def\PAMS             {Proc. Amer. Math. Soc.}
\def\PJM              {Pacific J. Math.}
\def\ijmp             {Int.\ J.\ Mod.\ Phys.\ {\bf A}}
\def\jpa              {J.\ Phys.\ {\bf A}}
\def\PLB              {Phys.\ Lett.\ {\bf B}}
\def\RIMS             {Publ.\ RIMS, Kyoto Univ.}
\def\Top               {Topology}
\def\TAMS             {Trans.\ Amer.\ Math.\ Soc.}

\def\Duke              {Duke Math.\ J.}
\def\K                 {K-theory}
\def\JOP               {J.\ Oper.\ Theory}

\vspace{0.2cm}\addtolength{\baselineskip}{-2pt}
\begin{footnotesize}
\noindent{\it Acknowledgement.}

The authors thank the University of Alberta Mathematics Dept, Cardiff School of Mathematics,  Neuadd Gregynog (University of Wales), University of Warwick 
Mathematics Institute,  Swansea University Dept of Computer Science, the  Isaac Newton Institute for Mathematical Sciences, Cambridge (during the programme Operator Algebras: Subfactors and their applications), and the Mathematical Sciences Research Institute Berkeley (NSF Grant No.\ DMS-1440140) during the Quantum Symmetries programme) for generous hospitality while researching this
paper. They also benefitted greatly from Research-in-Pairs held at Oberwolfach and BIRS,
and the von Neumann algebra trimester at the Hausdorff Institute. 
DEE thanks  WWU M\" unster and Wilhelm Winter whilst TG thanks Karlstads Universitet and J\"urgen Fuchs for stimulating environments while part of this paper was written. 
Their research was supported in part by  EPSRC grant nos EP/K032208/1 and EP/N022432/1, NSERC and SFB 878.

\end{footnotesize}

\end{document}